
\magnification=1200
\hsize=13.50cm
\vsize=18cm
\parindent=12pt   \parskip=0pt
\pageno=1

\ifnum\mag=\magstep1
\hoffset=-2mm
\voffset=.8cm
\fi

\pretolerance=500 \tolerance=1000  \brokenpenalty=5000

\catcode`\@=11

\font\eightrm=cmr8         \font\eighti=cmmi8
\font\eightsy=cmsy8        \font\eightbf=cmbx8
\font\eighttt=cmtt8        \font\eightit=cmti8
\font\eightsl=cmsl8        \font\sixrm=cmr6
\font\sixi=cmmi6           \font\sixsy=cmsy6
\font\sixbf=cmbx6

\font\tengoth=eufm10
\font\tenbboard=msbm10
\font\eightgoth=eufm10 at 8pt
\font\eightbboard=msbm10 at 8 pt
\font\sevengoth=eufm7
\font\sevenbboard=msbm7
\font\sixgoth=eufm7 at 6 pt
\font\fivegoth=eufm5

 \font\tencyr=wncyr10

\font\eightcyr=wncyr10 at 8 pt

\font\sevencyr=wncyr10 at 7 pt

\font\sixcyr=wncyr10 at 6 pt

\skewchar\eighti='177 \skewchar\sixi='177
\skewchar\eightsy='60 \skewchar\sixsy='60

\newfam\gothfam           \newfam\bboardfam
\newfam\cyrfam

\def\tenpoint{%
  \textfont0=\tenrm \scriptfont0=\sevenrm
  \scriptscriptfont0=\fiverm
  \def\rm{\fam\z@\tenrm}%
  \textfont1=\teni  \scriptfont1=\seveni
  \scriptscriptfont1=\fivei
  \def\oldstyle{\fam\@ne\teni}\let\old=\oldstyle
  \textfont2=\tensy \scriptfont2=\sevensy
  \scriptscriptfont2=\fivesy
  \textfont\gothfam=\tengoth \scriptfont\gothfam=\sevengoth
  \scriptscriptfont\gothfam=\fivegoth
  \def\goth{\fam\gothfam\tengoth}%
  \textfont\bboardfam=\tenbboard
  \scriptfont\bboardfam=\sevenbboard
  \scriptscriptfont\bboardfam=\sevenbboard
  \def\bb{\fam\bboardfam\tenbboard}%
  \textfont\cyrfam=\tencyr \scriptfont\cyrfam=\sevencyr
  \scriptscriptfont\cyrfam=\sixcyr
  \def\cyr{\fam\cyrfam\tencyr}%
  \textfont\itfam=\tenit
  \def\it{\fam\itfam\tenit}%
  \textfont\slfam=\tensl
  \def\sl{\fam\slfam\tensl}%
  \textfont\bffam=\tenbf \scriptfont\bffam=\sevenbf
  \scriptscriptfont\bffam=\fivebf
  \def\bf{\fam\bffam\tenbf}%
  \textfont\ttfam=\tentt
  \def\tt{\fam\ttfam\tentt}%
  \abovedisplayskip=12pt plus 3pt minus 9pt
  \belowdisplayskip=\abovedisplayskip
  \abovedisplayshortskip=0pt plus 3pt
  \belowdisplayshortskip=4pt plus 3pt
  \smallskipamount=3pt plus 1pt minus 1pt
  \medskipamount=6pt plus 2pt minus 2pt
  \bigskipamount=12pt plus 4pt minus 4pt
  \normalbaselineskip=12pt
  \setbox\strutbox=\hbox{\vrule height8.5pt depth3.5pt width0pt}%
  \let\bigf@nt=\tenrm       \let\smallf@nt=\sevenrm
  \normalbaselines\rm}

\def\eightpoint{%
  \textfont0=\eightrm \scriptfont0=\sixrm
  \scriptscriptfont0=\fiverm
  \def\rm{\fam\z@\eightrm}%
  \textfont1=\eighti  \scriptfont1=\sixi
  \scriptscriptfont1=\fivei
  \def\oldstyle{\fam\@ne\eighti}\let\old=\oldstyle
  \textfont2=\eightsy \scriptfont2=\sixsy
  \scriptscriptfont2=\fivesy
  \textfont\gothfam=\eightgoth \scriptfont\gothfam=\sixgoth
  \scriptscriptfont\gothfam=\fivegoth
  \def\goth{\fam\gothfam\eightgoth}%
  \textfont\cyrfam=\eightcyr \scriptfont\cyrfam=\sixcyr
  \scriptscriptfont\cyrfam=\sixcyr
  \def\cyr{\fam\cyrfam\eightcyr}%
  \textfont\bboardfam=\eightbboard
  \scriptfont\bboardfam=\sevenbboard
  \scriptscriptfont\bboardfam=\sevenbboard
  \def\bb{\fam\bboardfam}%
  \textfont\itfam=\eightit
  \def\it{\fam\itfam\eightit}%
  \textfont\slfam=\eightsl
  \def\sl{\fam\slfam\eightsl}%
  \textfont\bffam=\eightbf \scriptfont\bffam=\sixbf
  \scriptscriptfont\bffam=\fivebf
  \def\bf{\fam\bffam\eightbf}%
  \textfont\ttfam=\eighttt
  \def\tt{\fam\ttfam\eighttt}%
  \abovedisplayskip=9pt plus 3pt minus 9pt
  \belowdisplayskip=\abovedisplayskip
  \abovedisplayshortskip=0pt plus 3pt
  \belowdisplayshortskip=3pt plus 3pt
  \smallskipamount=2pt plus 1pt minus 1pt
  \medskipamount=4pt plus 2pt minus 1pt
  \bigskipamount=9pt plus 3pt minus 3pt
  \normalbaselineskip=9pt
  \setbox\strutbox=\hbox{\vrule height7pt depth2pt width0pt}%
  \let\bigf@nt=\eightrm     \let\smallf@nt=\sixrm
  \normalbaselines\rm}

\tenpoint

\def\pc#1{\bigf@nt#1\smallf@nt}         \def\pd#1 {{\pc#1} }

\catcode`\;=\active
\def;{\relax\ifhmode\ifdim\lastskip>\z@\unskip\fi
\kern\fontdimen2  -1.2 \fontdimen3 \string;}

\catcode`\:=\active
\def:{\relax\ifhmode\ifdim\lastskip>\z@
\unskip\fi\penalty\@M\ \fi\string:}

\catcode`\!=\active
\def!{\relax\ifhmode\ifdim\lastskip>\z@
\unskip\fi\kern\fontdimen2  -1.1
\fontdimen3 \string!}

\catcode`\?=\active
\def?{\relax\ifhmode\ifdim\lastskip>\z@
\unskip\fi\kern\fontdimen2  -1.1
\fontdimen3 \string?}

\def\^#1{\if#1i{\accent"5E\i}\else{\accent"5E #1}\fi}
\def\"#1{\if#1i{\accent"7F\i}\else{\accent"7F #1}\fi}

\frenchspacing

\newtoks\auteurcourant      \auteurcourant={\hfil}
\newtoks\titrecourant       \titrecourant={\hfil}

\newtoks\hautpagetitre      \hautpagetitre={\hfil}
\newtoks\baspagetitre       \baspagetitre={\hfil}

\newtoks\hautpagegauche
\hautpagegauche={\eightpoint\rlap{\folio}\hfil\the\auteurcourant
\hfil}
\newtoks\hautpagedroite
\hautpagedroite={\eightpoint\hfil\the\titrecourant\hfil
\llap{\folio}}

\newtoks\baspagegauche      \baspagegauche={\hfil}
\newtoks\baspagedroite      \baspagedroite={\hfil}

\newif\ifpagetitre          \pagetitretrue

\headline={\ifpagetitre\the\hautpagetitre
\else\ifodd\pageno\the\hautpagedroite\else\the\hautpagegauche
\fi\fi}

\footline={\ifpagetitre\the\baspagetitre\else
\ifodd\pageno\the\baspagedroite\else\the\baspagegauche\fi\fi
\global\pagetitrefalse}

\def\raggedbottom{\topskip 10pt plus 36pt\r@ggedbottomtrue}

\def\pointir{\unskip . --- \ignorespaces}

\def\Medbreak{\vskip-\lastskip\medbreak}

\def\rem#1\endrem{%
\Medbreak
{\it#1\unskip} : }

\long\def\th#1 #2\enonce#3\endth{%
   \Medbreak
   {\pc#1} {#2\unskip}\pointir{\it #3}\medskip}

\long\def\tha#1 #2\enonce#3\endth{%
   \Medbreak
   {\pc#1} {#2\unskip}\par\nobreak{\it #3}\medskip}

\def\decale#1{\smallbreak\hskip 28pt\llap{#1}\kern 5pt}
\def\decaledecale#1{\smallbreak\hskip 34pt\llap{#1}\kern 5pt}

\let\@ldmessage=\message

\def\message#1{{\def\pc{\string\pc\space}%
                \def\'{\string'}\def\`{\string`}%
                \def\^{\string^}\def\"{\string"}%
                \@ldmessage{#1}}}

\def\qed{\raise -2pt\hbox{\vrule\vbox to 10pt{\hrule width 4pt
                 \vfill\hrule}\vrule}}

\def\cqfd{\unskip\penalty 500\quad\qed\medbreak}

\catcode`\@=12

\showboxbreadth=-1  \showboxdepth=-1

\expandafter\ifx\csname amssym.def\endcsname\relax \else
\endinput\fi
%
\expandafter\edef\csname amssym.def\endcsname{%
       \catcode`\noexpand\@=\the\catcode`\@\space}
\catcode`\@=11
%

\def\undefine#1{\let#1\undefined}
\def\newsymbol#1#2#3#4#5{\let\next@\relax
 \ifnum#2=\@ne\let\next@\msafam@\else
 \ifnum#2=\tw@\let\next@\msbfam@\fi\fi
 \mathchardef#1="#3\next@#4#5}
\def\mathhexbox@#1#2#3{\relax
 \ifmmode\mathpalette{}{\m@th\mathchar"#1#2#3}%
 \else\leavevmode\hbox{$\m@th\mathchar"#1#2#3$}\fi}
\def\hexnumber@#1{\ifcase#1 0\or 1\or 2\or 3\or 4\or 5\or 6\or 7\or 8\or
 9\or A\or B\or C\or D\or E\or F\fi}

\font\tenmsa=msam10
\font\sevenmsa=msam7
\font\fivemsa=msam5
\newfam\msafam
\textfont\msafam=\tenmsa
\scriptfont\msafam=\sevenmsa
\scriptscriptfont\msafam=\fivemsa
\edef\msafam@{\hexnumber@\msafam}
\mathchardef\dabar@"0\msafam@39
\def\dashrightarrow{\mathrel{\dabar@\dabar@\mathchar"0\msafam@4B}}
\def\dashleftarrow{\mathrel{\mathchar"0\msafam@4C\dabar@\dabar@}}

\def\ulcorner{\delimiter"4\msafam@70\msafam@70 }
\def\urcorner{\delimiter"5\msafam@71\msafam@71 }
\def\llcorner{\delimiter"4\msafam@78\msafam@78 }
\def\lrcorner{\delimiter"5\msafam@79\msafam@79 }
\def\yen{{\mathhexbox@\msafam@55}}
\def\checkmark{{\mathhexbox@\msafam@58}}
\def\circledR{{\mathhexbox@\msafam@72}}
\def\maltese{{\mathhexbox@\msafam@7A}}

\font\tenmsb=msbm10
\font\sevenmsb=msbm7
\font\fivemsb=msbm5
\newfam\msbfam
\textfont\msbfam=\tenmsb
\scriptfont\msbfam=\sevenmsb
\scriptscriptfont\msbfam=\fivemsb
\edef\msbfam@{\hexnumber@\msbfam}

\def\widehat#1{\setbox\z@\hbox{$\m@th#1$}%
 \ifdim\wd\z@>\tw@ em\mathaccent"0\msbfam@5B{#1}%
 \else\mathaccent"0362{#1}\fi}
\def\widetilde#1{\setbox\z@\hbox{$\m@th#1$}%
 \ifdim\wd\z@>\tw@ em\mathaccent"0\msbfam@5D{#1}%
 \else\mathaccent"0365{#1}\fi}
\font\teneufm=eufm10
\font\seveneufm=eufm7
\font\fiveeufm=eufm5
\newfam\eufmfam
\textfont\eufmfam=\teneufm
\scriptfont\eufmfam=\seveneufm
\scriptscriptfont\eufmfam=\fiveeufm

\let\goth\frak

\csname amssym.def\endcsname
\expandafter\ifx\csname pre amssym.tex at\endcsname\relax \else
\endinput\fi
\expandafter\chardef\csname pre amssym.tex at\endcsname=\the
\catcode`\@
\catcode`\@=11
\begingroup\ifx\undefined\newsymbol \else\def\input#1 {\endgroup}\fi
\input amssym.def \relax
\newsymbol\boxdot 1200
\newsymbol\boxplus 1201
\newsymbol\boxtimes 1202
\newsymbol\square 1003
\newsymbol\blacksquare 1004
\newsymbol\centerdot 1205
\newsymbol\lozenge 1006
\newsymbol\blacklozenge 1007
\newsymbol\circlearrowright 1308
\newsymbol\circlearrowleft 1309
\undefine\rightleftharpoons
\newsymbol\rightleftharpoons 130A
\newsymbol\leftrightharpoons 130B
\newsymbol\boxminus 120C
\newsymbol\Vdash 130D
\newsymbol\Vvdash 130E
\newsymbol\vDash 130F
\newsymbol\twoheadrightarrow 1310
\newsymbol\twoheadleftarrow 1311
\newsymbol\leftleftarrows 1312
\newsymbol\rightrightarrows 1313
\newsymbol\upuparrows 1314
\newsymbol\downdownarrows 1315
\newsymbol\upharpoonright 1316
 
\newsymbol\downharpoonright 1317
\newsymbol\upharpoonleft 1318
\newsymbol\downharpoonleft 1319
\newsymbol\rightarrowtail 131A
\newsymbol\leftarrowtail 131B
\newsymbol\leftrightarrows 131C
\newsymbol\rightleftarrows 131D
\newsymbol\Lsh 131E
\newsymbol\Rsh 131F
\newsymbol\rightsquigarrow 1320
\newsymbol\leftrightsquigarrow 1321
\newsymbol\looparrowleft 1322
\newsymbol\looparrowright 1323
\newsymbol\circeq 1324
\newsymbol\succsim 1325
\newsymbol\gtrsim 1326
\newsymbol\gtrapprox 1327
\newsymbol\multimap 1328
\newsymbol\therefore 1329
\newsymbol\because 132A
\newsymbol\doteqdot 132B
 
\newsymbol\triangleq 132C
\newsymbol\precsim 132D
\newsymbol\lesssim 132E
\newsymbol\lessapprox 132F
\newsymbol\eqslantless 1330
\newsymbol\eqslantgtr 1331
\newsymbol\curlyeqprec 1332
\newsymbol\curlyeqsucc 1333
\newsymbol\preccurlyeq 1334
\newsymbol\leqq 1335
\newsymbol\leqslant 1336
\newsymbol\lessgtr 1337
\newsymbol\backprime 1038
\newsymbol\risingdotseq 133A
\newsymbol\fallingdotseq 133B
\newsymbol\succcurlyeq 133C
\newsymbol\geqq 133D
\newsymbol\geqslant 133E
\newsymbol\gtrless 133F
\newsymbol\sqsubset 1340
\newsymbol\sqsupset 1341
\newsymbol\vartriangleright 1342
\newsymbol\vartriangleleft 1343
\newsymbol\trianglerighteq 1344
\newsymbol\trianglelefteq 1345
\newsymbol\bigstar 1046
\newsymbol\between 1347
\newsymbol\blacktriangledown 1048
\newsymbol\blacktriangleright 1349
\newsymbol\blacktriangleleft 134A
\newsymbol\vartriangle 134D
\newsymbol\blacktriangle 104E
\newsymbol\triangledown 104F
\newsymbol\eqcirc 1350
\newsymbol\lesseqgtr 1351
\newsymbol\gtreqless 1352
\newsymbol\lesseqqgtr 1353
\newsymbol\gtreqqless 1354
\newsymbol\Rrightarrow 1356
\newsymbol\Lleftarrow 1357
\newsymbol\veebar 1259
\newsymbol\barwedge 125A
\newsymbol\doublebarwedge 125B
\undefine\angle
\newsymbol\angle 105C
\newsymbol\measuredangle 105D
\newsymbol\sphericalangle 105E
\newsymbol\varpropto 135F
\newsymbol\smallsmile 1360
\newsymbol\smallfrown 1361
\newsymbol\Subset 1362
\newsymbol\Supset 1363
\newsymbol\Cup 1264
 
\newsymbol\Cap 1265
 
\newsymbol\curlywedge 1266
\newsymbol\curlyvee 1267
\newsymbol\leftthreetimes 1268
\newsymbol\rightthreetimes 1269
\newsymbol\subseteqq 136A
\newsymbol\supseteqq 136B
\newsymbol\bumpeq 136C
\newsymbol\Bumpeq 136D
\newsymbol\lll 136E
 
\newsymbol\ggg 136F
 
\newsymbol\circledS 1073
\newsymbol\pitchfork 1374
\newsymbol\dotplus 1275
\newsymbol\backsim 1376
\newsymbol\backsimeq 1377
\newsymbol\complement 107B
\newsymbol\intercal 127C
\newsymbol\circledcirc 127D
\newsymbol\circledast 127E
\newsymbol\circleddash 127F
\newsymbol\lvertneqq 2300
\newsymbol\gvertneqq 2301
\newsymbol\nleq 2302
\newsymbol\ngeq 2303
\newsymbol\nless 2304
\newsymbol\ngtr 2305
\newsymbol\nprec 2306
\newsymbol\nsucc 2307
\newsymbol\lneqq 2308
\newsymbol\gneqq 2309
\newsymbol\nleqslant 230A
\newsymbol\ngeqslant 230B
\newsymbol\lneq 230C
\newsymbol\gneq 230D
\newsymbol\npreceq 230E
\newsymbol\nsucceq 230F
\newsymbol\precnsim 2310
\newsymbol\succnsim 2311
\newsymbol\lnsim 2312
\newsymbol\gnsim 2313
\newsymbol\nleqq 2314
\newsymbol\ngeqq 2315
\newsymbol\precneqq 2316
\newsymbol\succneqq 2317
\newsymbol\precnapprox 2318
\newsymbol\succnapprox 2319
\newsymbol\lnapprox 231A
\newsymbol\gnapprox 231B
\newsymbol\nsim 231C
\newsymbol\ncong 231D
\newsymbol\diagup 201E
\newsymbol\diagdown 201F
\newsymbol\varsubsetneq 2320
\newsymbol\varsupsetneq 2321
\newsymbol\nsubseteqq 2322
\newsymbol\nsupseteqq 2323
\newsymbol\subsetneqq 2324
\newsymbol\supsetneqq 2325
\newsymbol\varsubsetneqq 2326
\newsymbol\varsupsetneqq 2327
\newsymbol\subsetneq 2328
\newsymbol\supsetneq 2329
\newsymbol\nsubseteq 232A
\newsymbol\nsupseteq 232B
\newsymbol\nparallel 232C
\newsymbol\nmid 232D
\newsymbol\nshortmid 232E
\newsymbol\nshortparallel 232F
\newsymbol\nvdash 2330
\newsymbol\nVdash 2331
\newsymbol\nvDash 2332
\newsymbol\nVDash 2333
\newsymbol\ntrianglerighteq 2334
\newsymbol\ntrianglelefteq 2335
\newsymbol\ntriangleleft 2336
\newsymbol\ntriangleright 2337
\newsymbol\nleftarrow 2338
\newsymbol\nrightarrow 2339
\newsymbol\nLeftarrow 233A
\newsymbol\nRightarrow 233B
\newsymbol\nLeftrightarrow 233C
\newsymbol\nleftrightarrow 233D
\newsymbol\divideontimes 223E
\newsymbol\varnothing 203F
\newsymbol\nexists 2040
\newsymbol\Finv 2060
\newsymbol\Game 2061
\newsymbol\mho 2066
\newsymbol\eth 2067
\newsymbol\eqsim 2368
\newsymbol\beth 2069
\newsymbol\gimel 206A
\newsymbol\daleth 206B
\newsymbol\lessdot 236C
\newsymbol\gtrdot 236D
\newsymbol\ltimes 226E
\newsymbol\rtimes 226F
\newsymbol\shortmid 2370
\newsymbol\shortparallel 2371
\newsymbol\smallsetminus 2272
\newsymbol\thicksim 2373
\newsymbol\thickapprox 2374
\newsymbol\approxeq 2375
\newsymbol\succapprox 2376
\newsymbol\precapprox 2377
\newsymbol\curvearrowleft 2378
\newsymbol\curvearrowright 2379
\newsymbol\digamma 207A
\newsymbol\varkappa 207B
\newsymbol\Bbbk 207C
\newsymbol\hslash 207D
\undefine\hbar
\newsymbol\hbar 207E
\newsymbol\backepsilon 237F
\catcode`\@=\csname pre amssym.tex at\endcsname

\font\quatorzebf=cmbx10 at 14pt
\def\({{\rm (}}
\def\){{\rm )}}
\def\Lotimes{\smash{\mathop{\otimes}\limits^L}}
\def\Lboxtimes{\smash{\mathop{\boxtimes}\limits^L}}
\def\Lstar{\smash{\mathop{*}\limits^L}}

\def\maprightover#1{\smash{\mathop{\longrightarrow}
\limits^{#1}}}

\def\maprightunder#1{\smash{\mathop{\longrightarrow}
\limits_{#1}}}

\def\mapdownleft#1{\llap{$\vcenter
{\hbox{$\scriptstyle#1$}}$}\Big\downarrow}
\def\mapdownright#1{\Big\downarrow
\rlap{$\vcenter{\hbox{$\scriptstyle#1$}}$}}


\auteurcourant={G. LAUMON}
\titrecourant={TRANSFORMATION DE FOURIER}

\centerline{\quatorzebf Transformation de Fourier
g\'en\'eralis\'ee}
\vskip 5mm
\centerline{G\'erard Laumon}
\vskip 15mm

{\bf 0. Introduction}
\vskip 5mm

Soit $G$ un groupe alg\'ebrique ab\'elien connexe, extension  d'une
vari\'et\'e ab\'elienne par le produit d'un tore et d'un vectoriel.
Grothendieck a d\'efini  la notion de $\natural$-extensions de $G$
par ${\bb G}_{\rm m}$. Il s'agit de ${\cal O}_G$-Modules inversibles
\`a connexion int\'egrable $({\cal L},\nabla )$ satisfaisant le
th\'eor\`eme  du carr\'e, i.e. tels que
$$
\mu^*({\cal L},\nabla )\otimes {\rm pr}_1^*({\cal L},\nabla )^{-1}
\otimes {\rm pr}_2^*({\cal L},\nabla )^{-1}\cong ({\cal O}_G,d),
$$
o\`u $\mu :G\times G\rightarrow G$ est la loi de groupe et
${\rm pr}_1,{\rm pr}_2: G\times G\rightarrow G$ sont les deux
projections canoniques.

Notons $G^\natural$ l'ensemble de ces $\natural$-extensions et
${\cal D}_G$ l'Anneau des op\'erateurs diff\'erentiels sur $G$. On
s'int\'eresse au probl\`eme suivant: \'etant donn\'e un ${\cal
D}_G$-Module (ou un complexe de ${\cal D}_G$-Modules) ${\cal M}$,
dans quelle mesure la donn\'ee des complexes
$$
R\Gamma\bigl(G,{\rm DR}_G(({\cal L},\nabla )\otimes_{{\cal O}_G}
{\cal M})\bigr) \qquad (({\cal L},\nabla )\in G^\natural )
$$
d\'etermine-t-elle ${\cal M}$?

Nous r\'esolvons ce probl\`eme en munissant $G^\natural$ d'une
structure alg\'ebrique, puis en montrant que, pour $({\cal L},\nabla
)$ parcourant $G^\natural$, les complexes
$$
R\Gamma\bigl(G,{\rm DR}_G(({\cal L},\nabla )\otimes_{{\cal
O}_G}{\cal M})\bigr)
$$
s'organisent en un complexe de faisceaux ${\cal F}({\cal M})$ sur
$G^\natural$ et enfin en construisant une transformation inverse de
la ``transformation de Fourier'' ${\cal M}\mapsto {\cal F}({\cal M})$.
\vskip 2mm

Le plan de l'article est le suivant.

Apr\`es des rappels sur la transformation de Fourier-Mukai  pour
les vari\'et\'es ab\'eliennes (section 1) et sur l'extension
vectorielle universelle d'une vari\'et\'e ab\'elienne (section 2),
nous d\'efinissons une transformation de Fourier pour les
${\cal D}$-Modules sur les vari\'et\'es ab\'eliennes qui prolonge la
transformation de Fourier-Mukai et nous en \'etudions les
propri\'et\'es principales (section 3).

Le reste de l'article est consacr\'e \`a la construction d'une
transformation de Fourier qui englobe \`a la fois la transformation
de Fourier-Mukai, son prolongement aux ${\cal D}$-Modules, la
transformation de Fourier pour les ${\cal D}$-Modules sur les
vectoriels et la transformation de Mellin pour les ${\cal
D}$-Modules sur les tores. Pour cela, nous g\'en\'eralisons la notion
de $1$-motifs de Deligne et la dualit\'e de Cartier pour ces
$1$-motifs de mani\`ere \`a inclure les vectoriels (sections 4 et
5). Nous introduisons ensuite les notions de Modules
quasi-coh\'erents et coh\'erents sur ces $1$-motifs
g\'en\'eralis\'es et \'etudions les fonctorialit\'es naturelles de ces
objets (section 6). Enfin, nous d\'efinissons la transformation de
Fourier pour les Modules quasi-coh\'erents sur les $1$-motifs
g\'en\'eralis\'es (section 7).
\vskip 2mm

Cet article est une version remani\'ee d'une pr\'epublication de
l'auteur intitul\'ee ``Transformation de Fourier g\'eom\'etrique''
(IHES, septembre 1985). R\'ecemment, la transformation de  Fourier
g\'en\'eralis\'ee introduite dans cette  pr\'epublication a \'et\'e
utilis\'ee par Beilinson et Drinfeld dans leur programme de
Langlands g\'eom\'etrique. Cette transformation de Fourier joue un
r\^ole dans l'\'etude des faisceaux automorphes associ\'es aux
syst\`emes locaux de rang un sur une surface de Riemann compacte
connexe priv\'ee d'un ensemble fini de points. Les r\'esultats de la
section 3 ont aussi \'et\'e obtenus par M. Rothstein ([Ro]).
\vskip 2mm

Ce travail a \'et\'e commenc\'e lors d'un s\'ejour au D\'epartement
de Math\'ematiques de l'Universit\'e de  Harvard que je remercie
pour son hospitalit\'e. J'ai eu des discussions fructueuses avec L.
Breen, P. Deligne, D. Kazdhan, D. Mumford et R. Thomason durant la
pr\'eparation de cet article et J. Bernstein a simplifi\'e certains de
mes arguments originaux.

\vskip 5mm
{\bf 1. Transformation de Fourier-Mukai: rappels.}
\vskip 5mm

(1.1) Soient $S$ un sch\'ema localement noeth\'erien et $A$ un
$S$-sch\'ema ab\'elien purement de dimension relative $g$ ($A$
est donc un $S$-sch\'ema en groupes ab\'eliens dont le
$S$-sch\'ema sous-jacent est propre, lisse, purement de dimension
relative $g$ et \`a fibres g\'eom\'etriquement connexes). On  note
$\pi :A\rightarrow S$ le morphisme structural, $\mu :
A\times_SA\rightarrow A$ la loi de groupe, $\epsilon :
S\rightarrow A$ la section nulle, $\langle -1\rangle :A\rightarrow
A$ l'application inverse pour la loi de groupe et ${\rm pr}_1, {\rm
pr}_2: A\times_SA\rightarrow A$ les projections canoniques.

Un ${\cal O}_A$-Module inversible ${\cal L}$ est dit {\it rigidifi\'e}
s'il est muni d'un isomorphisme de ${\cal O}_S$-Modules ${\cal O}_S
\buildrel\sim\over\longrightarrow\epsilon^*{\cal L}$. On note
${\rm Pic}(A/S)$ le groupe ab\'elien (pour le produit tensoriel) des
classes d'isomorphie de ${\cal O}_A$-Modules inversibles
rigidifi\'es.

Pour tout ${\cal O}_A$-Module inversible ${\cal L}$, on peut former
le ${\cal O}_{A\times_SA}$-Module inversible
$$
{\cal D}_2({\cal L})=\mu^*{\cal L}\otimes {\rm pr}_1^*{\cal
L}^{\otimes -1}\otimes {\rm pr}_2^*{\cal L}^{\otimes -1}\otimes
(\pi\times\pi )^*\epsilon^*{\cal L}.
$$
Par construction, la restriction de ${\cal D}_2({\cal L})$ \`a
$(\epsilon (S)\times_SA)\cup (A\times_S\epsilon (S))$ est
canoniquement triviale. On dit qu'un ${\cal O}_A$-Module inversible
${\cal L}$ {\it satisfait le th\'eor\`eme du carr\'e} si la
trivialisation  canonique de ${\cal D}_2({\cal L})$ sur
$(\epsilon (S)\times_SA)\cup (A\times_S\epsilon (S))$ se prolonge
\`a $A\times_SA$ tout entier. On note ${\rm Pic}^0(A/S) \subset
{\rm Pic}(A/S)$ le sous-groupe form\'e des classes d'isomorphie de
${\cal O}_A$-Modules inversibles rigidifi\'es qui satisfont le
th\'eor\`eme du carr\'e.

Les applications $T\mapsto {\rm Pic}(A\times_ST/T)$ et
$T\mapsto {\rm Pic}^0(A\times_ST/T)$ se prolongent de mani\`ere
\'evidente en des foncteurs de la cat\'egorie des $S$-sch\'emas
localement noeth\'eriens dans la cat\'egorie des groupes ab\'eliens.

\th TH\'EOR\`EME (1.1.1) (Grothendieck dans le cas o\`u $A$ est
localement projectif sur $S$; Artin et Raynaud en g\'en\'eral)
\enonce
Le premier foncteur est rep\'esentable par un $S$-sch\'ema en
groupes ab\'eliens localement de type fini, not\'e
${\rm Pic}_{A/S}$.

Le deuxi\`eme foncteur est repr\'esentable par la composante
neutre ${\rm Pic}_{A/S}^0$ de ${\rm Pic}_{A/S}$ et cette
composante neutre est un $S$-sch\'ema ab\'elien purement de
dimension relative $g$.
\endth

\rem Preuve
\endrem
Voir l'expos\'e de Grothendieck au s\'eminaire Bourbaki [Gr] et la
discussion dans la section 1 du chapitre I du livre Faltings et Chai
([Fa-Ch]). Voir aussi [Mum] Ch. III, \S 13, pour le cas o\`u $S$ est le
spectre d'un corps.

\hfill\hfill\cqfd

Le $S$-sch\'ema ab\'elien ${\rm Pic}_{A/S}^0$ est appel\'e le
$S$-sch\'ema ab\'elien {\it dual} de $A$ et sera not\'e simplement
$A'$ dans la suite. On notera $\pi'$, $\mu'$, ...  les morphismes
$\pi$, $\mu$, ... relatifs \`a $A'$ et $\sigma :A'\times_SA
\rightarrow A\times_SA'$ l'application de permutation des deux
facteurs.

Par d\'efinition de $A'$, on dispose d'un ${\cal O}_{A'\times_S
A}$-Module inversible universel ${\cal P}$, dit {\it de Poincar\'e},
trivialis\'e le long de $A'\times_S\epsilon (S)$ et satisfaisant le
th\'eor\`eme du carr\'e pour le $A'$-sch\'ema ab\'elien
$A'\times_SA$. En fait, ${\cal P}$ est aussi trivialis\'e le long de
$\epsilon'(S)\times_SA$ et ces trivialisations co\"incident sur
$\epsilon'(S)\times_S\epsilon (S)$.

Soient $A''$ le $S$-sch\'ema ab\'elien dual de $A'$ et ${\cal P}'$ le
${\cal O}_{A''\times_SA'}$-Module de Poincar\'e correspondant.  On
dispose d'un morphisme de $S$-sch\'emas en groupes
$$
\iota :A\rightarrow A''
$$
qui, pour tout $S$-sch\'ema localement noeth\'erien $T$, envoie
$a\in A(T)$ sur la classe de $({\rm id}_{A'}\times a)^* {\cal P}$
dans $A''(T)$.

\th TH\'EOR\`EME (1.1.2)
\enonce
Le morphisme $\iota$ ci-dessus est un isomorphisme et
$\sigma^*(\iota\times {\rm id}_{A'})^*{\cal P}'$ est isomorphe \`a
${\cal P}$ en tant que ${\cal O}_{A'\times_SA}$-Module inversible
trivialis\'e le long de $(A'\times_S\epsilon (S))\cup
(\epsilon'(S)\times_SA)$.
\endth

\rem Preuve
\endrem
Voir [Mum] Ch. III, \S 13, pour le cas o\`u $S$ est le spectre d'un
corps.

\hfill\hfill\cqfd

En particulier, ${\cal P}$ satisfait aussi le th\'eor\`eme du carr\'e
pour le $A$-sch\'ema ab\'elien $A'\times_SA$. Dans la suite, on
identifiera $A''$ \`a $A$ et $\sigma^*{\cal P}'$ \`a ${\cal P}$ par
les  isomorphismes ci-dessus.

\th LEMME (1.1.3)
\enonce
Les ${\cal O}_S$-Modules inversibles
$$
\pi_*\Omega_{A/S}^g\cong\epsilon^*\Omega_{A/S}^g
$$
et
$$
\pi_*'\Omega_{A'/S}^g\cong\epsilon'^*\Omega_{A'/S}^g
$$
sont canoniquement isomorphes.
\endth

Dans la suite, on identifiera ces Modules inversibles et on les
notera simplement $\omega_{A/S}$.

\rem Preuve
\endrem
La th\'eorie des d\'eformations fournit un isomorphisme canonique
de ${\cal O}_S$-Modules
$$
R^1\pi_*{\cal O}_A\buildrel\sim\over\longrightarrow
\epsilon'^*{\cal T}_{A'/S},
$$
o\`u ${\cal T}_{A'/S}=\underline{\rm Hom}_{{\cal O}_{A'}}
(\Omega_{A'/S}^1,{\cal O}_{A'})$ est le ${\cal O}_{A'}$-Module
tangent relatif. Par suite, en prenant la puissance ext\'erieure
$g$-i\`eme de cet isomorphisme et en tenant compte de
l'isomorphisme canonique de ${\cal O}_S$-Modules $\bigwedge^g
R^1\pi_*{\cal O}_A\buildrel\sim\over \longrightarrow
R^g\pi_*{\cal O}_A$, on obtient un isomorphisme de ${\cal
O}_S$-Modules
$$
R^g\pi_*{\cal O}_A\buildrel\sim\over\longrightarrow
\epsilon'^*\bigl(\bigwedge^g{\cal T}_{A'/S}\bigr)
$$
et donc, par dualit\'e de Grothendieck, l'isomorphisme
$$
\epsilon'^*\Omega_{A'/S}^g\buildrel\sim\over\longrightarrow
\pi_*\Omega_{A/S}^g
$$
cherch\'e.

\hfill\hfill\cqfd
\vskip 5mm

(1.2) Notons ${\rm pr}':A'\times_SA\rightarrow A'$ et ${\rm pr}:
A'\times_SA\rightarrow A$ les projections canoniques. La {\it
transformation de Fourier-Mukai} pour $A$ est le foncteur
$$
{\cal F}:D_{\rm qcoh}^{\rm b}({\cal O}_A)\rightarrow D_{\rm
qcoh}^{\rm b}({\cal O}_{A'})
$$
d\'efini par
$$
{\cal F}(\cdot )=R{\rm pr}_*'({\cal P}
\,\Lotimes_{{\cal O}_{A'\times_SA}}\, L{\rm pr}^*(\cdot ))\leqno
(1.2.1)
$$
(cf. [Muk]). Par construction ${\cal F}$ est un foncteur triangul\'e.

Il r\'esulte du th\'eor\`eme de changement de base que:

\th LEMME (1.2.2)
\enonce
Pout tout morphisme $T\buildrel f\over\longrightarrow S$ de
sch\'emas localement noeth\'eriens on a un isomorphisme
canonique
$$
L({\rm id}_{A'}\times f)^*\circ {\cal F}\buildrel\sim\over
\longrightarrow {\cal F}_T\circ L({\rm id}_A\times f)^*
$$
de foncteurs de $D_{\rm qcoh}^{\rm b}({\cal O}_A)$ dans
$D_{\rm qcoh}^{\rm b}({\cal O}_{A'\times_ST})$, o\`u ${\cal F}_T$
est la transformation de Fourier-Mukai pour le $T$-sch\'ema
ab\'elien $A\times_ST$ \(de dual identifi\'e \`a $A'\times_ST$\).

\hfill\hfill\cqfd
\endth

Il r\'esulte du th\'eor\`eme de finitude de Grothendieck pour les
images directes par les morphismes propres que:

\th PROPOSITION (1.2.3)
\enonce
La transformation de Fourier-Mukai ${\cal F}$ envoie la
sous-cat\'egorie strictement pleine $D_{\rm coh}^{\rm b}({\cal
O}_A)$ de $D_{\rm qcoh}^{\rm b}({\cal O}_A)$ dans la
sous-cat\'egorie strictement pleine $D_{\rm coh}^{\rm b}({\cal
O}_{A'})$ de $D_{\rm qcoh}^{\rm b}({\cal O}_{A'})$.

\hfill\hfill\cqfd
\endth

La propri\'et\'e la plus important de la transformation de
Fourier-Mukai est l'involutivit\'e.

\th TH\'EOR\`EME (1.2.4) (Mukai)
\enonce
Si l'on note
$$
{\cal F}':D_{\rm qcoh}^{\rm b}({\cal O}_{A'})\rightarrow D_{\rm
qcoh}^{\rm b}({\cal O}_A)
$$
la transformation de Fourier-Mukai pour le $S$-sch\'ema ab\'elien
$A'$, les foncteurs compos\'es ${\cal F}'\circ {\cal F}$ et
${\cal F}\circ {\cal F}'$ sont canoniquement isomorphes aux
foncteurs
$$
\langle -1\rangle^*(\cdot )\,\Lotimes_{{\cal O}_A}\,
\pi^*\omega_{A/S}^{\otimes -1}[-g]: D_{\rm qcoh}^{\rm b}({\cal
O}_A)\rightarrow D_{\rm qcoh}^{\rm b}({\cal O}_A)
$$
et
$$
\langle -1\rangle'^*(\cdot )\,\Lotimes_{{\cal O}_{A'}}\,
\pi'^*\omega_{A/S}^{\otimes -1}[-g]: D_{\rm qcoh}^{\rm b}({\cal
O}_{A'})\rightarrow D_{\rm qcoh}^{\rm b}({\cal O}_{A'})
$$
respectivement. En particulier, ${\cal F}$ est une \'equivalence de
cat\'egories de quasi-inverse $\langle -1\rangle^*{\cal F}'(\cdot )
\,\Lotimes_{{\cal O}_A}\,\pi^*\omega_{A/S} [g]$.
\endth

\rem Preuve
\endrem
Compte tenu du th\'eor\`eme de changement de base et de la
formule des projections, ce th\'eor\`eme est une cons\'equence
formelle du lemme ci-dessous.

\hfill\hfill\cqfd

\th LEMME (1.2.5)
\enonce
Le complexe $R{\rm pr}_*{\cal P}$ \(resp. $R{\rm pr}_*'{\cal P}$\)
est canoniquement isomorphe \`a
$\epsilon_*(\omega_{A/S}^{\otimes -1})[-g]$ \(resp.
$\epsilon_*'(\omega_{A/S}^{\otimes -1})[-g]$\) dans
$D_{\rm qcoh}^{\rm b}({\cal O}_A)$ \(resp. $D_{\rm qcoh}^{\rm
b}({\cal O}_{A'})$\).
\endth

\rem Preuve
\endrem
Nous rappellerons seulement la construction du morphisme
canonique $R{\rm pr}_*{\cal P}\rightarrow
\epsilon_*(\omega_{A/S}^{\otimes -1})[-g]$. Par changement de
base on a un isomorphisme
$$
\epsilon'^*R^g{\rm pr}_*{\cal P}\buildrel\sim\over\longrightarrow
R^g\pi_*'{\cal O}_{A'}
$$
et, par dualit\'e de Gothendieck, on a un isomorphisme
$$
R^g\pi_*'{\cal O}_{A'}\buildrel\sim\over\longrightarrow
\omega_{A/S}^{\otimes -1},
$$
d'o\`u, par adjunction, un morphisme
$$
R{\rm pr}_*{\cal P}\rightarrow R^g{\rm pr}_*{\cal P}[-g]
\rightarrow\epsilon_*(\omega_{A/S}^{\otimes -1})[-g]
$$
qui est le morphisme cherch\'e.

\hfill\hfill\cqfd
\vskip 5mm

(1.3) La transformation de Fourier-Mukai poss\`ede les
propri\'et\'es de fonctorialit\'e suivantes.

\th PROPOSITION (1.3.1)
\enonce
Soient $f:A_1\rightarrow A_2$ un morphisme de $S$-sch\'emas
ab\'eliens et $f':A_2'\rightarrow A_1'$ le morphisme  transpos\'e.
Alors on a des isomorphismes canoniques de  foncteurs
$$
{\cal F}_2\circ Rf_*\cong Lf'^*\circ {\cal F}_1
$$
et
$$
{\cal F}_1\circ Rf^!\cong Rf_*'\circ {\cal F}_2,
$$
o\`u ${\cal F}_1$ et ${\cal F}_2$ sont les transformations de
Fourier-Mukai pour $A_1$ et $A_2$ respectivement.

\hfill\hfill\cqfd
\endth

\th PROPOSITION (1.3.2)
\enonce
Soient $A_1$ et $A_2$ deux $S$-sch\'emas ab\'eliens et
$A=A_1\times_SA_2$. Alors on a un isomorphisme canonique de
foncteurs
$$
{\cal F}\bigl((\cdot )\,\Lboxtimes_S\,(\cdot )\bigr)\cong {\cal
F}_1(\cdot )\,\Lboxtimes_S\,{\cal F}_2(\cdot ),
$$
o\`u ${\cal F}$, ${\cal F}_1$ et ${\cal F}_2$ sont les
transformations de Fourier-Mukai pour $A$, $A_1$ et $A_2$
respectivement \(on a bien entendu identifi\'e $A'$ \`a
$A_1'\times_SA_2'$\).

\hfill\hfill\cqfd
\endth

On appelle {\it produit de convolution} pour $A$ et on note
$$
(\cdot )\,\Lstar\,(\cdot ):D_{\rm qcoh}^{\rm b}({\cal O}_A)\times
D_{\rm qcoh}^{\rm b}({\cal O}_A)\rightarrow D_{\rm qcoh}^{\rm
b}({\cal O}_A)
$$
le foncteur d\'efini par
$$
(\cdot )\,\Lstar\,(\cdot )=R\mu_*\bigl((\cdot )
\,\Lboxtimes_S\,(\cdot )\bigr).
$$

\th COROLLAIRE (1.3.3)
\enonce
On a des isomorphismes de foncteurs
$$
{\cal F}\bigl((\cdot )\,\Lstar\,(\cdot )\bigr)\cong {\cal F}(\cdot
)\,\Lotimes_{{\cal O}_A}\,{\cal F}(\cdot )
$$
et
$$
{\cal F}\bigl((\cdot )\,\Lotimes_{{\cal O}_A}\,(\cdot )\bigr)
\cong \bigl({\cal F}(\cdot )\,\Lstar\,{\cal F}(\cdot )\bigr)
\,\Lotimes_{{\cal O}_A}\,\pi^*\omega_{A/S} [g].
$$

\hfill\hfill\cqfd
\endth

Enfin la transformation de Fourier-Mukai commute \`a la dualit\'e:

\th PROPOSITION (1.3.4)
\enonce
Soient
$$
D(\cdot )\buildrel{\rm dfn}\over{=\!=} R\underline{\rm Hom}_{{\cal
O}_A}(\,\cdot\,,\pi^*\omega_{A/S})[g]: D_{\rm coh}^{\rm b}({\cal
O}_A)^{\rm opp}\rightarrow D_{\rm coh}^{\rm b}({\cal O}_A)
$$
et
$$
D'(\cdot )\buildrel{\rm dfn}\over{=\!=} R\underline{\rm
Hom}_{{\cal O}_{A'}}(\,\cdot\,,\pi'^*\omega_{A/S})[g]: D_{\rm
coh}^{\rm b}({\cal O}_{A'})^{\rm opp}\rightarrow D_{\rm coh}^{\rm
b}({\cal O}_{A'})
$$
les foncteurs de dualit\'e. Alors, on a un isomorphisme canonique de
foncteurs
$$
D'\circ {\cal F}\cong \langle -1\rangle'^*({\cal F} ( D(\cdot
)))\,\Lotimes_{{\cal O}_A}\,\pi'^*\omega_{A/S} [g].
$$

\hfill\hfill\cqfd
\endth
\vskip 5mm

{\bf 2. L'extension vectorielle universelle: rappels et
compl\'ements.}
\vskip 5mm

(2.1) Soit $A$ un sch\'ema ab\'elien purement de dimension relative
$g$ sur un sch\'ema localement noeth\'erien $S$.

Un ${\cal O}_A$-Module inversible ${\cal L}$ muni d'une connexion
int\'egrable (relative \`a $S$)
$$
\nabla :{\cal L}\rightarrow \Omega_{A/S}^1\otimes_{{\cal O}_A}
{\cal L}
$$
est dit {\it rigidifi\'e} si ${\cal L}$ l'est.

Pour tout ${\cal O}_A$-Module inversible \`a connexion int\'egrable
(relative \`a $S$) $({\cal L},\nabla )$, on peut former le
${\cal O}_{A\times_SA}$-Module \`a connexion int\'egrable
(relative \`a $S$)
$$
{\cal D}_2({\cal L},\nabla )=\mu^*({\cal L},\nabla )\otimes {\rm
pr}_1^*({\cal L},\nabla )^{\otimes -1}\otimes {\rm pr}_2^*({\cal
L},\nabla )^{\otimes -1}\otimes (\pi\times \pi )^*\epsilon^*({\cal
L},\nabla ).
$$
Par construction, la restriction de ${\cal D}_2({\cal L},\nabla )$
\`a $(\epsilon (S)\times_SA)\cup (A\times_S\epsilon (S))$ est
canoniquement triviale. On dit qu'un ${\cal O}_A$-Module inversible
\`a connexion int\'egrable (relative \`a $S$) $({\cal L},\nabla )$ {\it
satisfait le th\'eor\`eme du carr\'e} si la trivialisation canonique
de ${\cal D}_2({\cal L},\nabla )$ sur $(\epsilon (S)\times_SA)\cup
(A\times_S\epsilon (S))$ se prolonge \`a $A\times_SA$ tout entier,
de sorte que
$$
{\cal D}_2({\cal L},\nabla )\cong ({\cal O}_{A\times_SA},d)
$$
o\`u $d: {\cal O}_{A\times_SA} \rightarrow
\Omega_{A\times_SA/S}^1$ est la diff\'erentielle.

\th LEMME (2.1.1)
\enonce
Soit $({\cal L},\nabla )$ un ${\cal O}_A$-Module inversible \`a
connexion int\'egrable \(relative \`a $S$\). Alors les conditions
suivantes sont \'equivalentes:

\decale{\rm (i)} $({\cal L},\nabla )$ satisfait le th\'eor\`eme du
carr\'e,

\decale{\rm (ii)} ${\cal L}$ satisfait le th\'eor\`eme du carr\'e.
\endth

\rem Preuve
\endrem
L'implication (i)$\Rightarrow$(ii) est imm\'ediate.

R\'eciproquement, si ${\cal L}$ satisfait le th\'eor\`eme du carr\'e,
on a
$$
{\cal D}_2({\cal L},\nabla )\cong ({\cal O}_{A\times_SA},d+\omega )
$$
pour une forme diff\'erentielle $\omega\in H^0(A\times_SA,
\Omega_{A\times_SA}^1)$ telle que
$$
(\epsilon\times {\rm id}_A)^*\omega = ({\rm id}_A\times\epsilon
)^*\omega=0.
$$
Mais une telle forme est automatiquement nulle puisque
$A\times_SA$ est un $S$-sch\'ema ab\'elien.

\hfill\hfill\cqfd

On note ${\rm Pic}^\natural (A/S)$ le groupe (pour le produit
tensoriel) des classes d'isomorphie de ${\cal O}_A$-Modules
inversibles \`a connexion int\'egrable (relative \`a $S$) rigidifi\'es
qui satisfont le th\'eor\`eme du carr\'e. On a une suite exacte de
groupes ab\'eliens
$$
{\rm E}(A/S)=(0\rightarrow H^0(A,\Omega_A^1)\buildrel
i\over\longrightarrow {\rm Pic}^\natural (A/S)\buildrel
p\over\longrightarrow {\rm Pic}^0(A/S)),
$$
o\`u $i(\omega )=({\cal O}_A,d+\omega )$ et $p({\cal L},\nabla )
={\cal L}$. De plus, si $S$ est affine, l'homomorphisme $p$ est
surjectif (cf. [Ma-Me] Prop. (3.2.3)(a) et  section (2.6)).

L'application $T\mapsto {\rm E}(A\times_ST/T)$ se prolonge de
mani\`ere \'evidente en un foncteur de la cat\'egorie des
$S$-sch\'emas localement noeth\'eriens dans la cat\'egorie des
suites exactes de groupes ab\'eliens.

\th TH\'EOR\`EME (2.1.2) (Mazur et Messing)
\enonce
Ce foncteur est repr\'esentable par une suite exacte de
$S$-sch\'emas en groupes ab\'eliens
$$
0\rightarrow {\bb V}(\epsilon^*{\cal T}_{A/S})\buildrel
i\over\longrightarrow {\rm Pic}_{A/S}^\natural\buildrel
p\over\longrightarrow {\rm Pic}_{A/S}^0\rightarrow 0.
$$

En particulier, le foncteur $T\mapsto {\rm Pic}^\natural
(A\times_ST/T)$ est repr\'esentable par un $S$-sch\'ema en
groupes ab\'eliens ${\rm Pic}_{A/S}^\natural $ dont le
$S$-sch\'ema sous-jacent est de type fini, lisse, purement de
dimension relative $2g$ et \`a fibres g\'eom\'etriquement connexes.
\endth

\rem Preuve
\endrem
Voir [Ma-Me] Ch. I, (2.6) et (3.2.3).

\hfill\hfill\cqfd

On notera encore
$$
{\rm E}_{A/S}=(0\rightarrow\Omega_{A/S}\buildrel
i\over\longrightarrow A^\natural\buildrel p\over\longrightarrow
A'\rightarrow 0)
$$
la suite exacte ci-dessus et on notera $\pi^\natural$ ,
$\mu^\natural$, ... les morphismes $\pi$, $\mu$, ... relatifs \`a
$A^\natural$. L'image r\'eciproque
$$
\widetilde{\cal P}\buildrel{\rm dfn}\over {=\!=} (p\times {\rm
id}_A)^*{\cal P}
$$
du Module de Poincar\'e est un ${\cal O}_{A^\natural
\times_SA}$-Module inversible muni d'une connexion int\'egrable
(relative \`a $A^\natural$) universelle
$$
\widetilde\nabla :\widetilde{\cal P}\rightarrow
\Omega_{A^\natural \times_SA/A^\natural}^1
\otimes_{{\cal O}_{A^\natural\times_SA}}\widetilde{\cal P}.
$$
\vskip 5mm

(2.2) On identifiera de la mani\`ere habituelle la cat\'egorie des
$S$-sch\'emas \`a une sous-cat\'egorie pleine de la cat\'egorie des
faisceaux sur le site fppf de $S$. Pour tout ${\cal O}_S$-Module
localement libre de rang fini ${\cal E}$, de dual $\underline{\rm
Hom}_{{\cal O}_S}({\cal E},{\cal O}_S)$ not\'e simplement ${\cal
E}^\vee$, on identifiera donc ${\bb V}({\cal E}^\vee)$ \`a ${\cal E}$.
Par exemple, on identifiera $\Omega_{A/S}={\bb V}(\epsilon^*{\cal
T}_{A/S})$ \`a $\epsilon^*\Omega_{A/S}^1$.

L'extension $A^\natural$ de $A'$ par le $S$-sch\'ema en vectoriels
$\Omega_{A/S}$ est universelle pour les extensions (dans la
cat\'egorie des $S$-sch\'emas en groupes ab\'eliens) de $A'$ par
les $S$-sch\'emas en vectoriels:

\th TH\'EOR\`EME (2.2.1) (Mazur et Messing)
\enonce
Si ${\cal M}$ est un ${\cal O}_S$-Module quasi-coh\'erent, toute
extension de $A'$ par ${\cal M}$ dans la cat\'egorie des faisceaux
en groupes ab\'eliens sur le site fppf de $S$ est le ``push-out'' de
l'extension ${\rm E}_{A/S}$ ci-dessus par un et un seul morphisme
de faisceaux en groupes $\Omega_{A/S}\rightarrow {\cal M}$. En
d'autres termes, on a un isomorphisme, fonctoriel en ${\cal M}$,
$$
{\rm Hom}_{{\cal O}_S}(\epsilon^*\Omega_{A/S}^1, {\cal M})
\buildrel\sim\over\longrightarrow {\rm Ext}_{S_{\rm fppf}}^1
(A',{\cal M}).
$$

En particulier, pour tout ${\cal O}_S$-Module localement libre de
rang fini ${\cal E}$ et toute extension de $A'$ par ${\bb V}({\cal
E}^\vee)$ dans la cat\'egorie des $S$-sch\'emas en groupes
ab\'eliens, il existe un et un seul morphisme de $S$-sch\'emas en
groupes $\Omega_{A/S}\rightarrow {\bb V}({\cal E}^\vee )$ tel que
cette extension soit le ``push-out'' de l'extension ${\rm E}_{A/S}$
par ce morphisme.
\endth

\rem Preuve
\endrem
Voir [Ma-Me] Ch. I, (1.9) et (1.10).

\hfill\hfill\cqfd

D'apr\`es Mazur et Messing (cf. [Ma-Me] Ch. I, \S 4), l'extension de
${\cal O}_S$-Alg\`ebres de Lie (commutatives) associ\'ee
\`a l'extension ${\rm E}_{A/S}$ n'est autre que l'extension
$$
0\rightarrow\epsilon^*\Omega_{A/S}^1\rightarrow {\cal H}_{\rm
dR}^1(A/S)\rightarrow R^1\pi_*{\cal O}_A\rightarrow 0
$$
donn\'ee par la suite spectrale de la cohomologie de Hodge vers la
cohomologie de de Rham pour $A/S$, suite spectrale  qui
d\'eg\'en\`ere en $E_1$.
\vskip 5mm

(2.3) Soient $\pi :X\rightarrow S$ un $S$-sch\'ema et ${\cal E}$ un
${\cal O}_S$-Module localement libre de rang fini. On consid\`ere
les trois cat\'egories suivantes:

\decale{\rm (1)} la cat\'egorie ${\cal C}_1$ des torseurs
$Y\buildrel p\over\longrightarrow X$ sous le $X$-sch\'ema en
vectoriels ${\bb V}(\pi^*{\cal E}^\vee )=X\times_S {\bb V}({\cal
E}^\vee )$ (il s'agit a priori de torseurs dans la cat\'egorie des
faisceaux sur le site fppf de $X$ mais tout tel torseur est
repr\'esentable par un $X$-sch\'ema affine),

\decale{\rm (2)} la cat\'egorie ${\cal C}_2$ des
${\cal O}_X$-Alg\`ebres (commutatives et avec unit\'e)
quasi-coh\'erentes munies d'une filtration exhaustive
$$
(0)={\cal B}_{-1}\subset {\cal B}_0\subset {\cal B}_1\subset {\cal
B}_2\subset\cdots\subset {\cal B}
$$
par des sous-${\cal O}_X$-Modules localement libres de rang fini
faisant de ${\cal B}$ une ${\cal O}_X$-Alg\`ebre filtr\'ee, i.e. telle
que
$$
{\cal B}_i\cdot {\cal B}_j\subset {\cal B}_{i+j}\qquad (\forall
i,j\in {\bb N}),
$$
et munies d'un isomorphisme de ${\cal O}_X$-Alg\`ebres gradu\'ees
$$
\varphi_\bullet : {\rm gr}_\bullet {\cal B}\buildrel\sim
\over\longrightarrow {\rm Sym}_\bullet^{{\cal O}_X} (\pi^*{\cal
E}^\vee)=\pi^*{\rm Sym}_\bullet^{{\cal O}_S} ({\cal E}^\vee),
$$

\decale{\rm (3)} la cat\'egorie ${\cal C}_3$ des extensions de
$\pi^*{\cal E}^\vee$ par ${\cal O}_X$ dans la cat\'egorie ab\'elienne
des ${\cal O}_X$-Modules quasi-coh\'erents.
\vskip 2mm
On d\'efinit des foncteurs
$$
{\cal C}_1\rightarrow {\cal C}_2\rightarrow {\cal C}_3
\rightarrow {\cal C}_1
$$
de la fa\c con suivante.

\decale{$\bullet$} Si $Y\buildrel p\over\longrightarrow X$ est un
objet de ${\cal C}_1$, ${\cal B}=p_*{\cal O}_Y$ est une
${\cal O}_X$-Alg\`ebre quasi-coh\'erente que l'on peut filtrer par
les sous-${\cal O}_X$-Modules localement libres de rang fini
$$
{\cal B}_i=(\mu^*)^{-1}\bigl({\cal B}\otimes_{{\cal O}_S}
\bigoplus_{j=0}^i {\rm Sym}_j^{{\cal O}_S}({\cal E}^\vee)\bigr),
$$
o\`u $\mu :Y\times_S{\bb V}({\cal E}^\vee)\rightarrow Y$ est
l'action de $X\times_S{\bb V}({\cal E}^\vee)$ sur $Y$ et
$\mu^*:{\cal B}\rightarrow {\cal B}\otimes_{{\cal O}_S} {\rm
Sym}_\bullet^{{\cal O}_S}({\cal E}^\vee)$ est l'homomorphisme
d'Alg\`ebres correspondant; pour tout entier $i\geq 0$ et toute
section  $b$ de ${\cal B}_i$, la projection de $\mu^*b$ sur ${\cal
B}\otimes_{{\cal O}_S} {\rm Sym}_i^{{\cal O}_S}({\cal E}^\vee)$ est
en fait contenue dans
$$
{\cal O}_X\otimes_{{\cal O}_S}{\rm Sym}_i^{{\cal O}_S} ({\cal
E}^\vee)=\pi^*{\rm Sym}_i^{{\cal O}_S}({\cal E}^\vee)\hbox{;}
$$
si on note $\varphi_i(b+{\cal B}_{i-1})$ cette projection, on
v\'erifie facilement que $({\cal B},{\cal B}_\bullet ,\varphi_\bullet
)$ est un objet de ${\cal C}_2$.

\decale{$\bullet$} Si $({\cal B},{\cal B}_\bullet ,\varphi_\bullet )$
est un objet de ${\cal C}_2$, on a une extension
$$
0\rightarrow {\cal B}_0\rightarrow {\cal B}_1\rightarrow  {\cal
B}_1/{\cal B}_0\rightarrow 0\hbox{;}
$$
or on peut identifier ${\cal B}_0$ \`a ${\cal O}_X$ via $\varphi_0$
et ${\cal B}_1/{\cal B}_0$ \`a $\pi^*{\cal E}^\vee$ via $\varphi_1$
et on obtient ainsi une extension de $\pi^*{\cal E}^\vee$ par ${\cal
O}_X$, i.e. un objet de ${\cal C}_3$.

\decale{$\bullet$} Enfin, si
$$
0\rightarrow {\cal O}_X\rightarrow {\cal F}\rightarrow
\pi^*{\cal E}^\vee\rightarrow 0,
$$
est un objet de ${\cal C}_3$, on consid\`ere l'ouvert
$$
{\bb P}({\cal F})-{\bb P}(\pi^*{\cal E}^\vee)=Y\rightarrow X
$$
du fibr\'e projectif (au sens de Grothendieck) ${\bb P}({\cal F})$ sur
$X$; $Y$ est le $X$-sch\'ema des scindages de l'extension ci-dessus
et est donc naturellement muni d'une structure de ${\bb
V}(\pi^*{\cal E}^\vee)$-torseur; c'est l'objet de ${\cal C}_1$ voulu.

\th LEMME (2.3.1)
\enonce
{\rm (i)} Les trois foncteurs que l'on vient de d\'efinir sont des
\'equivalences de cat\'egories. Plus pr\'ecis\'ement, pour chaque
$\alpha\in {\bb Z}/3{\bb Z}$, le foncteur compos\'e
$$
{\cal C}_\alpha\rightarrow {\cal C}_{\alpha +1}\rightarrow {\cal
C}_{\alpha +2}\rightarrow {\cal C}_{\alpha +3}= {\cal C}_\alpha
$$
est isomorphe au foncteur identit\'e.

\decale{\rm (ii)} La classe d'un objet $Y\buildrel p\over
\longrightarrow X$ de ${\cal C}_1$ dans $H^1(X,\pi^*{\cal E})$ n'est
autre que la classe de l'objet correspondant de ${\cal C}_3$ dans
${\rm Ext}_{{\cal O}_X}^1(\pi^*{\cal E}^\vee ,{\cal O}_X)$, compte
tenu de l'identification naturelle entre $H^1(X,\pi^*{\cal E})$ et
${\rm Ext}_{{\cal O}_X}^1(\pi^*{\cal E}^\vee ,{\cal O}_X)$.
\endth

\rem Preuve
\endrem
Il suffit de remarquer que, pour tout objet $Y\buildrel p\over
\longrightarrow X$ de ${\cal C}_1$, le $X$-sch\'ema $Y$ n'est autre
que ${\rm Spec}({\cal B})$ o\`u ${\cal B}=p_*{\cal O}_Y$ et que,
pour tout objet $({\cal B},{\cal B}_\bullet ,\varphi_\bullet )$ de
${\cal C}_2$, ${\cal B}$ est isomorphe \`a
$$
{\rm Sym}_\bullet^{{\cal O}_X}({\cal B}_1)/(1-\xi )
$$
o\`u $1\in {\cal O}_X={\rm Sym}_0^{{\cal O}_X}({\cal B}_1)$ et o\`u
$\xi$ est l'image de $1\in {\cal O}_X={\cal B}_0$ par l'inclusion
${\cal B}_0\hookrightarrow {\cal B}_1= {\rm Sym}_1^{{\cal
O}_X}({\cal B}_1)$.

\hfill\hfill\cqfd

On peut appliquer les consid\'erations pr\'ec\'edentes au
$(A'\times_S\Omega_{A/S})$-torseur $A^\natural\buildrel p\over
\longrightarrow A'$ sous-jacent \`a l'extension vectorielle
universelle. La ${\cal O}_{A'}$-Alg\`ebre quasi-coh\'erente
$$
{\cal A}^\natural =p_*{\cal O}_{A^\natural}
$$
est donc munie d'une filtration
$$
(0)={\cal A}^\natural_{-1}\subset {\cal A}^\natural_0\subset {\cal
A}^\natural_1\subset {\cal A}^\natural_2\subset\cdots\subset
{\cal A}^\natural
$$
par des sous-${\cal O}_{A'}$-Modules localement libres de rang fini
qui en fait une ${\cal O}_{A'}$-Alg\`ebre filtr\'ee et est munie d'un
isomorphisme de ${\cal O}_{A'}$-Alg\`ebres gradu\'ees
$$
\varphi_\bullet :{\rm gr}_\bullet {\cal A}^\natural\buildrel
\sim\over\longrightarrow\pi'^*{\rm Sym}_\bullet^{{\cal O}_S}
(\epsilon^*{\cal T}_{A/S}).
$$
De plus, l'extension
$$
0\rightarrow {\cal O}_{A'}\rightarrow {\cal A}_1^\natural
\rightarrow\pi'^*\epsilon^*{\cal T}_{A/S}\rightarrow 0\leqno
(2.3.2)
$$
d\'etermine le $(A'\times_S\Omega_{A/S})$-torseur $A^\natural
\buildrel p\over\longrightarrow A'$ \`a isomorphisme unique pr\`es.

On va maintenant d\'eterminer la classe de cette extension  dans
$$
{\rm Ext}_{{\cal O}_{A'}}^1(\pi'^*\epsilon^*{\cal T}_{A/S}, {\cal
O}_{A'})\cong H^1(A',\pi'^*\epsilon^*\Omega_{A/S}^1).
$$

La suite exacte des termes de bas degr\'es de la suite spectrale de
Leray pour $(\pi',\pi'^*\epsilon^*\Omega_{A/S}^1)$ s'\'ecrit
$$
0\rightarrow H^1(S,\epsilon^*\Omega_{A/S}^1)\buildrel
\alpha\over\longrightarrow H^1(A',\pi'^*\epsilon^*\Omega_{A/S}^1)
\buildrel\beta\over\longrightarrow H^0(S,R^1\pi_*'\pi'^*
\epsilon^*\Omega_{A/S}^1)
$$
puisque $\pi_*'{\cal O}_{A'}={\cal O}_S$ et, comme $\pi'$ admet la
section $\epsilon'$, la fl\`eche $\beta$ est surjective et la fl\`eche
$\alpha$ admet la r\'etraction $\epsilon'^*$. On a donc une
d\'ecomposition en somme directe
$$
H^1(A',\pi'^*\epsilon^*\Omega_{A/S}^1)\cong
H^1(S,\epsilon^*\Omega_{A/S}^1)\oplus
H^0(S,R^1\pi_*'\pi'^*\epsilon^*\Omega_{A/S}^1).
$$
La restriction de l'extension (2.3.2) \`a la section nulle
$\epsilon'(S)$ est canoniquement scind\'ee puisque $A^\natural
\buildrel p\over\longrightarrow A'$ est non seulement un
$(A'\times_S\Omega_{A/S})$-torseur mais encore une extension de
$A'$ par $\Omega_{A/S}$. Par suite, la classe de l'extension (2.3.2)
appartient au facteur direct $H^0(S,R^1\pi_*'\pi'^*\epsilon^*
\Omega_{A/S}^1)$ de ${\rm Ext}_{{\cal O}_{A'}}^1
(\pi'^*\epsilon^*{\cal T}_{A/S}, {\cal O}_{A'})$. Mais, comme
$\epsilon^*\Omega_{A/S}^1$ est localement libre et que
$R^1\pi_*'{\cal O}_{A'}$ est isomorphe \`a
$\epsilon^*{\cal T}_{A/S}$ par bi-dualit\'e et th\'eorie de la
d\'eformation, on a un isomorphisme canonique
$$
\epsilon^*{\cal T}_{A/S}\otimes_{{\cal O}_S}
\epsilon^*\Omega_{A/S}^1\buildrel\sim\over
\longrightarrow R^1\pi_*'\pi'^*\epsilon^*\Omega_{A/S}^1,
$$
de sorte que la classe de l'extension (2.3.2) peut \^etre vue comme
un \'el\'ement de
$$
H^0(S,\epsilon^*{\cal T}_{A/S}\otimes_{{\cal O}_S}
\epsilon^*\Omega_{A/S}^1)\cong {\rm Hom}_{{\cal O}_S}
(\epsilon^*\Omega_{A/S}^1,\epsilon^*\Omega_{A/S}^1).
$$

\th PROPOSITION (2.3.3)
\enonce
La classe de l'extension $(2.3.2)$ dans
$$
{\rm Hom}_{{\cal O}_S}
(\epsilon^*\Omega_{A/S}^1,\epsilon^*\Omega_{A/S}^1)
$$
n'est autre que l'identit\'e de $\epsilon^*\Omega_{A/S}^1$.
\endth

\rem Preuve
\endrem
D'apr\`es [Se] Ch. VII, \S 15, on a la suite exacte
$$
0\rightarrow {\rm Ext}_{S_{\rm fppf}}^1(A', \epsilon^*
\Omega_{A/S}^1) \rightarrow H^1(A',\pi'^*\epsilon^*
\Omega_{A/S}^1)\buildrel\epsilon'^*\over\longrightarrow
H^1(S,\epsilon^*\Omega_{A/S}^1)
$$
(cf. [Ma-Me] (1.10)). La proposition r\'esulte donc de  l'isomorphisme
$$
{\rm Hom}_{{\cal O}_S}(\epsilon^*\Omega_{A/S}^1,
\epsilon^*\Omega_{A/S}^1)\buildrel\sim\over\longrightarrow {\rm
Ext}_{S_{\rm fppf}}^1(A',\epsilon^*\Omega_{A/S}^1).
$$
du th\'eor\`eme (2.2.1) et du lemme (2.3.1).

\hfill\hfill\cqfd
\vskip 5mm

(2.4) On suppose dans cette section que $S$ {\it est de
caract\'eristique nulle}. Sous cette hypoth\`ese, on a

\th TH\'EOR\`EME (2.4.1)
\enonce
La fl\`eche d'adjonction
$$
{\cal O}_S\rightarrow R\pi_*^\natural {\cal O}_{A^\natural}
$$
est un isomorphisme dans $D_{\rm qcoh}^{\rm b}({\cal O}_S)$, i.e.
la fl\`eche d'adjonction ${\cal O}_S\rightarrow \pi_*^\natural
{\cal O}_{A^\natural}$ est un isomorphisme et $R^n\pi_*^\natural
{\cal O}_{A^\natural}=(0)$ pour tout entier $n\not= 0$.
\endth

\rem Preuve
\endrem
Consid\'erons la suite spectrale
$$
E_1^{pq}=R^{p+q}\pi_*'{\rm gr}_{-p}{\cal A}^\natural
\Rightarrow R^{p+q}\pi_*'{\cal A}^\natural = R^{p+q}\pi_*^\natural
{\cal O}_{A^\natural}
$$
associ\'ee \`a l'objet filtr\'e $R\pi_*'({\cal A}^\natural , {\cal
A}_\bullet^\natural )$ de la cat\'egorie d\'eriv\'ee filtr\'ee $D^{\rm
b}F_{\rm qcoh}({\cal O}_S)$. On va montrer que
$$
E_2^{pq}=\left\{\matrix{{\cal O}_S\hfill &\hbox{si }p=q=0,\hfill\cr
\noalign{\smallskip} (0)\hfill &\hbox{sinon,}\hfill\cr}\right.
$$
ce qui bien entendu impliquera le th\'eor\`eme.

Il revient au m\^eme de montrer que le complexe de
${\cal O}_S$-Modules
$$
0\rightarrow E_1^{p,-p}\buildrel d_1\over\longrightarrow
E_1^{p+1,-p}\buildrel d_1\over\longrightarrow\cdots
\buildrel d_1\over\longrightarrow E_1^{p',-p}\rightarrow 0,
$$
o\`u $p'={\rm Inf}(0,p+g)$, est acyclique pour tout entier
$p<0$ (on a trivialement $E_1^{pq}=(0)$ si $p+q\notin [0,g]$ ou si
$p>0$). Or, d'apr\`es la formule des projections, l'isomorphisme
$\varphi_{-p}$ induit un isomorphisme
$$
E_1^{pq}\cong (R^{p+q}\pi_*'{\cal O}_{A'})\otimes_{{\cal O}_S} {\rm
Sym}_{-p}^{{\cal O}_S}(\epsilon^*{\cal T}_{A/S})
$$
i.e. un isomorphisme
$$
E_1^{pq}\cong\bigl(\bigwedge_{{\cal O}_S}^{p+q}\epsilon^* {\cal
T}_{A/S}\bigr)\otimes_{{\cal O}_S} {\rm Sym}_{-p}^{{\cal O}_S}
(\epsilon^*{\cal T}_{A/S})
$$
puisque le cup-produit induit un isomorphisme
$$
\bigwedge_{{\cal O}_S}^{p+q}R^1\pi_*'{\cal O}_{A'}
\buildrel\sim\over\longrightarrow R^{p+q}\pi_*'{\cal O}_{A'}
$$
et que
$$
R^1\pi_*'{\cal O}_{A'}\cong \epsilon^*{\cal T}_{A/S}
$$
par bi-dualit\'e et th\'eorie des d\'eformations. De plus, il est
facile de voir, compte tenu des r\'esultats de la section (2.3), que
la diff\'erentielle $d_1:E_1^{pq}\rightarrow E_1^{p+1,q}$ est
donn\'ee par la formule
$$\displaylines{
(2.4.2)\qquad d_1\bigl((v_1\wedge\cdots\wedge v_{p+q})\otimes
(w_1\cdots w_{-p})\bigr)
\hfill\cr\hfill
\sum_{i=1}^{-p}(v_1\wedge\cdots\wedge v_{p+q}\wedge w_i)
\otimes (w_1\cdots w_{i-1}w_{i+1}\cdots w_{-p})\qquad}
$$
pour toutes sections locales $v_1,\ldots ,v_{p+q},w_1,\ldots
,w_{-p}$ de $\epsilon^*{\cal T}_{A/S}$. Le th\'eor\`eme r\'esulte
donc du lemme bien connu suivant.

\th LEMME (2.4.3) (Koszul)
\enonce
Soit $V$ un espace vectoriel de dimension finie sur un corps $k$ de
caract\'eristique nulle. Pour tout entier $n>0$, le complexe de
$k$-espaces vectoriels
$$
0\rightarrow {\rm S}^nV\rightarrow V\otimes {\rm
S}^{n-1}V\rightarrow\bigwedge^2V\otimes {\rm
S}^{n-2}V\rightarrow\cdots \rightarrow \bigwedge^{n-1}V\otimes
V\rightarrow\bigwedge^nV\rightarrow 0,
$$
o\`u ${\rm S}^iV$ est une notation abr\'eg\'ee pour
${\rm Sym}_i^k(V)$ et o\`u la diff\'erentielle $\bigwedge^iV
\otimes {\rm S}^{n-i}V\rightarrow\bigwedge^{i+1} V\otimes {\rm
S}^{n-i-1}V$ est d\'efinie par la formule $(2.4.2)$ pour $p+q=i$ et
$-p=n$, est acyclique.

\hfill\hfill\cqfd
\endth
\vskip 5mm

{\bf 3. Une nouvelle transformation de Fourier pour les sch\'emas
ab\'eliens en caract\'erisitique nulle.}
\vskip 5mm

(3.1) Soient $S$ un sch\'ema localement noeth\'erien {\it de
caract\'eristique nulle} et $A$ un $S$-sch\'ema ab\'elien  purement
de dimension relative $g$.

On note ${\cal D}_{A/S}$ la ${\cal O}_A$-Alg\`ebre des op\'erateurs
diff\'erentiels sur $A$ relatifs \`a $S$. Rappelons qu'un ${\cal
D}_{A/S}$-Module (\`a gauche) quasi-coh\'erent n'est rien d'autre
qu'un ${\cal O}_A$-Module quasi-coh\'erent muni d'une connexion
int\'egrable relative \`a $S$.

On note $\widetilde{\rm pr}^\natural :A^\natural\times_SA
\rightarrow A^\natural$ et $\widetilde{\rm pr}:A^\natural
\times_SA\rightarrow A$ les deux projections canoniques.

Si ${\cal M}$ est un ${\cal D}_{A/S}$-Module (\`a gauche), le ${\cal
O}_{A^\natural\times_SA}$-Module
$$
(\widetilde{\cal P},\widetilde\nabla )\,\Lotimes_{{\cal
O}_{A^\natural \times_SA}}\,\widetilde{\rm pr}^*{\cal M}
$$
sera muni de la connexion int\'egrable  (relative \`a $A^\natural$)
produit tensoriel de $\widetilde\nabla$ par la connexion sur
$\widetilde{\rm pr}^*{\cal M}$ image inverse de celle de ${\cal M}$.

Si $({\cal M},\nabla )$ est un ${\cal O}_{A^\natural\times_S
A}$-Module \`a connexion int\'egrable (relative \`a $A^\natural$),
on notera ${\rm DR}_{A^\natural\times_SA/A^\natural} ({\cal
M},\nabla )$ son complexe de de Rham
$$
[{\cal M}\buildrel\nabla\over\longrightarrow \Omega_{A^\natural
\times_SA/A^\natural}^1\otimes_{{\cal O}_{A^\natural
\times_SA}}{\cal M}\buildrel \nabla \over \longrightarrow
\cdots\buildrel\nabla\over\longrightarrow\Omega_{A^\natural
\times_SA/A^\natural}^g\otimes_{{\cal O}_{A^\natural
\times_SA}}{\cal M}]\hbox{;}
$$
c'est un complexe de $((\widetilde{\rm pr}^\natural )^{-1} {\cal
O}_{A^\natural})$-Modules concentr\'e en degr\'es compris entre
$-g$ et $0$.

Si ${\cal M}^\natural$ est un ${\cal O}_{A^\natural}$-Module, le
${\cal O}_{A^\natural\times_SA}$-Module
$$
(\widetilde{\cal P},\widetilde\nabla ) \otimes_{{\cal
O}_{A^\natural\times_SA}}\widetilde{\rm pr}^{\natural *}{\cal
M}^\natural
$$
sera muni de la connexion int\'egrable (relative \`a $A^\natural$)
induite par $\widetilde\nabla$.

On d\'efinit alors des foncteurs triangul\'es
$$
\widetilde{\cal F}:D_{\rm qcoh}^{\rm b}({\cal D}_{A/S}) \rightarrow
D_{\rm qcoh}^{\rm b}({\cal O}_{A^\natural})
$$
et
$$
\widetilde{\cal F}^\natural : D_{\rm qcoh}^{\rm b}({\cal
O}_{A^\natural})\rightarrow D_{\rm qcoh}^{\rm b}({\cal D}_{A/S})
$$
par
$$
\widetilde{\cal F}(\cdot )=R\widetilde{\rm pr}_*^\natural {\rm
DR}_{A^\natural\times_SA/A^\natural}\bigl((\widetilde {\cal
P},\widetilde\nabla )\,\Lotimes_{{\cal O}_{A^\natural
\times_SA}}\, L\widetilde{\rm pr}^*(\cdot )\bigr)
$$
et
$$
\widetilde{\cal F}^\natural(\cdot )=R\widetilde{\rm pr}_*
\bigl((\widetilde{\cal P},\widetilde\nabla )
\,\Lotimes_{{\cal O}_{A^\natural\times_SA}}\, L\widetilde{\rm
pr}^{\natural *}(\cdot )\bigr).
$$

Il r\'esulte du th\'eor\`eme de changement de base que:

\th LEMME (3.1.1)
\enonce
Pout tout $S$-sch\'ema localement noeth\'erien $T\buildrel
f\over\longrightarrow S$ on a des isomorphismes canoniques de
foncteurs
$$
L({\rm id}_{A^\natural}\times f)^*\circ\widetilde{\cal F}
\buildrel\sim\over\longrightarrow\widetilde{\cal F}_T
\circ L({\rm id}_A\times f)^*
$$
de $D_{\rm qcoh}^{\rm b}({\cal D}_{A/S})$ dans $D_{\rm qcoh}^{\rm
b}({\cal O}_{A^\natural\times_ST})$ et
$$
L({\rm id}_A\times f)^*\circ\widetilde{\cal F}^\natural
\buildrel\sim\over\longrightarrow\widetilde{\cal F}_T^\natural
\circ L({\rm id}_{A^\natural}\times f)^*
$$
de $D_{\rm qcoh}^{\rm b}({\cal O}_{A^\natural})$ dans $D_{\rm
qcoh}^{\rm b}({\cal D}_{A\times_ST/T})$, o\`u
$\widetilde{\cal F}_T$ et $\widetilde{\cal F}_T^\natural$ sont les
foncteurs $\widetilde{\cal F}$ et $\widetilde{\cal F}^\natural$ pour
le $T$-sch\'ema ab\'elien $A\times_ST$ \(on a identifi\'e $(A
\times_ST)^\natural$ \`a $A^\natural\times_ST$\).

\hfill\hfill\cqfd
\endth

La transformation de Fourier $\widetilde{\cal F}$ est reli\'ee
\`a la transformation de Fourier-Mukai de la fa\c con suivante

\th PROPOSITION (3.1.2)
\enonce
Le carr\'e de foncteurs
$$
\matrix{
D_{\rm qcoh}^{\rm b}({\cal O}_A)&\maprightover{{\cal F}}&D_{\rm
qcoh}^{\rm b}({\cal O}_{A'})\cr
\noalign{\medskip}
\mapdownleft{{\cal D}_{A/S}\,\Lotimes_{{\cal O}_A}\,(\cdot )}
&&\mapdownright{Lp^*(\cdot )\,\Lotimes_{{\cal O}_{A^\natural}}\,
\pi^{\natural *}\omega_{A/S}}\cr
\noalign{\medskip}
D_{\rm qcoh}^{\rm b}({\cal D}_{A/S})&\maprightunder{\widetilde
{\cal F}}& D_{\rm qcoh}^{\rm b}({\cal O}_{A^\natural})\cr}
$$
est commutatif \(\`a un isomorphisme canonique pr\`es\).
\endth

\th COROLLAIRE (3.1.3)
\enonce
Le foncteur $\widetilde{\cal F}$ envoie la sous-cat\'egorie
strictement pleine $D_{\rm coh}^{\rm b}({\cal D}_{A/S})$ de
$D_{\rm qcoh}^{\rm b}({\cal D}_{A/S})$ dans la sous-cat\'egorie
strictement pleine $D_{\rm coh}^{\rm b}({\cal O}_{A^\natural})$ de
$D_{\rm qcoh}^{\rm b}({\cal O}_{A^\natural})$.
\endth

\rem Preuve du corollaire
\endrem
Localement sur $S$, tout ${\cal D}_{A/S}$-Module coh\'erent admet
une r\'esolution par des ${\cal D}_{A/S}$-Modules coh\'erents
induits, i.e. de la forme ${\cal D}_{A/S}\otimes_{{\cal O}_A}{\cal
E}$ avec ${\cal E}$ un ${\cal O}_A$-Module coh\'erent. Or, d'apr\`es
les propositions $(3.1.2)$ et $(1.2.3)$, $\widetilde{\cal F}$ envoie
un tel ${\cal D}_{A/S}$-Module coh\'erent induit dans $D_{\rm
coh}^{\rm b}({\cal O}_{A^\natural})$, d'o\`u le corollaire.

\hfill\hfill\cqfd

Pour d\'emontrer la proposition nous aurons besoin du r\'esultat
suivant.

\th LEMME (3.1.4)
\enonce
Soient $T$ un $S$-sch\'ema et $X$ un $T$-sch\'ema lisse. Soit
${\cal E}$ un ${\cal O}_X$-Module localement libre de rang fini
muni d'une connexion int\'egrable \(relative \`a $T$\) $\nabla :{\cal
E}\rightarrow\Omega_{X/T}^1\otimes_{{\cal O}_X}{\cal E}$ et soit
${\cal F}$ un ${\cal O}_X$-Module. On peut former d'une part le
${\cal D}_{X/T}$-Module \`a gauche
$$
({\cal E},\nabla )\otimes_{{\cal O}_X}({\cal D}_{X/T} \otimes_{{\cal
O}_X}{\cal F}),
$$
produit tensoriel des ${\cal D}_{X/T}$-Modules \`a gauche $({\cal
E},\nabla )$ et ${\cal D}_{X/T}\otimes_{{\cal O}_X}{\cal F}$, et
d'autre part le ${\cal D}_{X/T}$-Module \`a gauche
$$
{\cal D}_{X/T}\otimes_{{\cal O}_X}({\cal E}\otimes_{{\cal
O}_X}{\cal F}),
$$
qui lui ne d\'epend pas de $\nabla$. Alors la fl\`eche
$$\matrix{
{\cal D}_{X/T}\otimes_{{\cal O}_X}({\cal E}\otimes_{{\cal
O}_X}{\cal F})&\rightarrow & ({\cal E},\nabla )\otimes_{{\cal
O}_X}({\cal D}_{X/T}\otimes_{{\cal O}_X}{\cal F})\cr
\noalign{\medskip}
\hfill P\otimes (e\otimes f)&\mapsto &P\cdot (e\otimes (1\otimes
f))\hfill\cr}
$$
est un isomorphisme entre ces ${\cal D}_{X/T}$-Modules.
\endth

\rem Preuve
\endrem
Si $P\in {\cal D}_{X/T,i}$ ($P$ est de degr\'e $\leq i$), on a
$$
P\cdot (e\otimes (1\otimes f))\in {\cal E}\otimes_{{\cal O}_X}
({\cal D}_{X/T,i}\otimes_{{\cal O}_X}{\cal F})
$$
et
$$
P\cdot (e\otimes (1\otimes f))-e\otimes (P\otimes f)\in {\cal
E}\otimes_{{\cal O}_X}({\cal D}_{X/T,i-1}\otimes_{{\cal O}_X}{\cal
F}).
$$
Par suite, la fl\`eche du lemme est filtr\'ee, i.e. envoie ${\cal
D}_{X/T,i}\otimes_{{\cal O}_X}({\cal E}\otimes_{{\cal O}_X}{\cal
F})$ dans ${\cal E}\otimes_{{\cal O}_X}({\cal D}_{X/T,i}
\otimes_{{\cal O}_X}{\cal F})$, et induit une fl\`eche entre les
gradu\'es correspondants qui n'est autre que le morphisme de ${\cal
O}_X$-Modules
$$\matrix{
{\rm gr}_\bullet {\cal D}_{X/T}\otimes_{{\cal O}_X}({\cal E}
\otimes_{{\cal O}_X}{\cal F})&\rightarrow&{\cal E}
\otimes_{{\cal O}_X}({\rm gr}_\bullet {\cal D}_{X/T}
\otimes_{{\cal O}_X}{\cal F}).\cr
\hfill \xi\otimes (e\otimes f)&\mapsto &
e\otimes (\xi\otimes f)\hfill\cr}
$$
Mais ce dernier morphisme est clairement un isomorphisme, d'o\`u
la conclusion.

\hfill\hfill\cqfd

\rem Preuve de la proposition $(3.1.2)$
\endrem
D'apr\`es le lemme ci-dessus, on a canoniquement
$$
\widetilde{\cal F}({\cal D}_{A/S}\,\Lotimes_{{\cal O}_A}\, (\cdot
)) =R\widetilde{\rm pr}_*^\natural {\rm DR}_{A^\natural
\times_SA/A^\natural}\bigl({\cal D}_{A^\natural
\times_SA/A^\natural}\,\Lotimes_{{\cal O}_{A^\natural
\times_SA}}\, (\widetilde{\cal P}\,\Lotimes_{{\cal O}_{A^\natural
\times_SA}}\, L\widetilde{\rm pr}^*(\cdot ))\bigr).
$$
Or, pour tout $S$-sch\'ema $T$ et tout $T$-sch\'ema $X$, lisse et
purement de dimension relative $d$, on a
$$
{\rm DR}_{X/T}({\cal D}_{X/T}\,\Lotimes_{{\cal O}_X}\,(\cdot ))
\buildrel\sim\over\longrightarrow\Omega_{X/T}^d\,\Lotimes_{{\cal
O}_X}\,(\cdot ).
$$
Par suite, on a canoniquement
$$
\widetilde{\cal F}({\cal D}_{A/S}\,\Lotimes_{{\cal O}_A}\, (\cdot
)) =R\widetilde{\rm pr}_*^\natural\bigl((\pi^\natural \times\pi
)^*\omega_{A/S}\,\Lotimes_{{\cal O}_{A^\natural\times_SA}}\,
(\widetilde{\cal P}\,\Lotimes_{{\cal O}_{A^\natural\times_SA}}\,
L\widetilde{\rm pr}^*(\cdot ))\bigr)
$$
et la proposition s'en suit compte tenu du th\'eor\`eme de
changement de base pour ${\rm pr}'$.

\hfill\hfill\cqfd
\vskip 5mm

(3.2) Comme pour la transformation de Fourier-Mukai, la
propri\'et\'e la plus importante du couple $(\widetilde{\cal F},
\widetilde{\cal F}^\natural )$ est le th\'eor\`eme suivant.

\th TH\'EOR\`EME (3.2.1)
\enonce
Les foncteurs compos\'es $\widetilde{\cal F}^\natural\circ
\widetilde{\cal F}$ et $\widetilde{\cal F}\circ\widetilde{\cal
F}^\natural$ sont canoniquement isomorphes aux foncteurs
$$
\langle -1\rangle^*(\cdot )[-g]:D_{\rm qcoh}^{\rm b}({\cal D}_{A/S})
\rightarrow D_{\rm qcoh}^{\rm b}({\cal D}_{A/S})
$$
et
$$
\langle -1\rangle^{\natural *}(\cdot )[-g]: D_{\rm qcoh}^{\rm b}
({\cal O}_{A^\natural})\rightarrow D_{\rm qcoh}^{\rm b}({\cal
O}_{A^\natural})
$$
respectivement. En d'autres termes, $\widetilde{\cal F}$ est une
\'equivalence de cat\'egories de quasi-inverse $\langle
-1\rangle^*\widetilde{\cal F}^\natural (\cdot )[g]$.
\endth

De nouveau, compte tenu du th\'eor\`eme de changement de  base et
de la formule des projections, ce th\'eor\`eme est une
cons\'equence formelle du lemme suivant.

\th LEMME (3.2.2)
\enonce
Le complexe
$$
R\widetilde{\rm pr}_*^\natural {\rm DR}_{A^\natural
\times_SA/A^\natural} (\widetilde{\cal P},\widetilde\nabla
)\quad\hbox{\(resp. } R\widetilde{\rm pr}_*(\widetilde{\cal P},
\widetilde\nabla )\,)
$$
est canoniquement isomorphe \`a $\epsilon_*^\natural {\cal
O}_S[-g]$ \(resp. $\epsilon_+{\cal O}_S[-g]$\) dans $D_{\rm
qcoh}^{\rm b} ({\cal O}_{A^\natural})$ \(resp. $D_{\rm qcoh}^{\rm b}
({\cal D}_{A/S})$\).
\endth

\vskip 2mm
{\pc RAPPEL} (3.2.3)\pointir
Le ${\cal D}_{A/S}$-Module \`a gauche $\epsilon_+{\cal O}_S$ est
d\'efini par
$$
\epsilon_+{\cal O}_S=\epsilon_*{\cal D}_{(A\buildrel\epsilon
\over\longleftarrow S)/S},
$$
o\`u
$$
{\cal D}_{(A\buildrel\epsilon\over\longleftarrow S)/S}=
\epsilon^{-1}\bigr({\cal D}_{A/S}\otimes_{{\cal O}_A}
(\Omega_{A/S}^g)^{\otimes -1}\bigl)\otimes_{\epsilon^{-1} {\cal
O}_A}{\cal O}_S
$$
est muni de sa structure de $(\epsilon^{-1}{\cal D}_{A/S}, {\cal
O}_S)$-bi-Module gauche-droite naturelle (la structure de
$(\epsilon^{-1}{\cal D}_{A/S})$-Module \`a gauche vient de celle de
${\cal D}_{A/S}$-Module \`a gauche de ${\cal D}_{A/S}
\otimes_{{\cal O}_A} (\Omega_{A/S}^g)^{\otimes -1}$, qui elle
m\^eme vient des structures de ${\cal D}_{A/S}$-Module \`a droite
de ${\cal D}_{A/S}$ et de $\Omega_{A/S}^g$, et la structure de
${\cal O}_S$-Module \`a droite est la structure \'evidente). Voir [Bo] \S 7
pour plus de d\'etails.

\hfill\hfill\cqfd
\vskip 2mm

\rem Preuve du lemme
\endrem
Par changement de base, le complexe
$$
L\epsilon^{\natural *}R\widetilde{\rm pr}_*^\natural {\rm
DR}_{A^\natural\times_SA/A^\natural} (\widetilde{\cal
P},\widetilde\nabla )
$$
est canoniquement isomorphe \`a
$$
R\pi_*{\rm DR}_{A/S}({\cal O}_A,d)
$$
et, comme le ${\cal O}_S$-Module
$$
{\cal H}_{\rm dR}^{2g}(A/S)=R^g\pi_*{\rm DR}_{A/S} ({\cal O}_A,d)
$$
est canoniquement isomorphe \`a ${\cal O}_S$, on en d\'eduit une
fl\`eche
$$
L\epsilon^{\natural *}R\widetilde{\rm pr}_*^\natural {\rm
DR}_{A^\natural\times_SA/A^\natural} (\widetilde{\cal P},
\widetilde\nabla )\rightarrow  {\cal O}_S[-g]
$$
et donc par adjonction une fl\`eche
$$
R\widetilde{\rm pr}_*^\natural {\rm DR}_{A^\natural
\times_SA/A^\natural} (\widetilde{\cal P},\widetilde\nabla )
\rightarrow \epsilon_*^\natural {\cal O}_S[-g].
$$
On va montrer que cette derni\`ere fl\`eche est un isomorphisme.

Remarquons tout d'abord que $R\widetilde{\rm pr}_*^\natural {\rm
DR}_{A^\natural\times_SA/A^\natural} (\widetilde{\cal P},
\widetilde\nabla )$ est un objet de $D_{\rm coh}^{\rm b} ({\cal
O}_{A^\natural})$ puisqu'il en est ainsi de
$$
R\widetilde{\rm pr}_*^\natural (\Omega_{A^\natural
\times_SA/A^\natural}^i\otimes_{A^\natural
\times_SA}\widetilde{\cal P}) =Lp^*{\cal F}(\Omega_{A/S}^i)
$$
pour tout entier $i$. D'apr\`es le lemme de Nakayama, on peut donc
supposer que $S$ est le spectre d'un corps de caract\'eristique
nulle et il suffit de montrer que, pour  tout ${\cal O}_A$-Module
inversible \`a connexion int\'egrable $({\cal L},\nabla )$
satisfaisant le th\'eor\`eme du carr\'e et tout entier $i$, le
$k$-espace vectoriel
$$
H_{\rm dR}^{i+g}({\cal L},\nabla )=H^i(A,{\rm DR}_{A/k}({\cal
L},\nabla ))
$$
est nul si $({\cal L},\nabla )$ n'est pas isomorphe \`a $({\cal
O}_A,d)$ et est canoniquement isomorphe \`a
$$
L_{g-i}\epsilon^{\natural *}\epsilon_*^\natural k
$$
si $({\cal L},\nabla )=({\cal O}_A,d)$.

La preuve de cette assertion est en tout point similaire \`a [Mum]
Ch. III, \S 16, PROOFS (2), et nous ne ferons que l'esquisser.

Si $({\cal L},\nabla)$ n'est pas isomorphe \`a $({\cal O}_A,d)$, on a
$H_{\rm dR}^0({\cal L},\nabla )=(0)$. Raisonnons par l'absurde en
supposant qu'il existe des entiers $n>0$ tels que $H_{\rm dR}^n
({\cal L},\nabla )\not=(0)$ et notons $n_0$ le plus petit de ces
entiers. Comme $({\cal L},\nabla )$ satisfait le th\'eor\`eme du
carr\'e, on a un isomorphisme de $k$-espaces vectoriels
$$
\bigoplus_{i+j=n_0}H_{\rm dR}^i({\cal L},\nabla )\otimes_k H_{\rm
dR}^j({\cal L},\nabla )\cong\bigoplus_{i+j=n_0} H_{\rm dR}^i({\cal
L},\nabla )\otimes_k H_{\rm dR}^j(A/k)
$$
(formule de K\" unneth et formule de projection). Or, $H_{\rm
dR}^0(A/k)=k$ n'est pas nul, d'o\`u une contradiction.

Si $({\cal L},\nabla )=({\cal O}_A,d)$, pour tout entier $i$,
$H_{\rm dR}^i({\cal L},\nabla )=H_{\rm dR}^i(A/k)$ est
canoniquement isomorphe \`a $\bigwedge^iH_{\rm dR}^1(A/k)$.
Mais, la dualit\'e locale fournit, pour tout entier $i$, un
isomorphisme canonique
$$
\epsilon_*^\natural L_{g-i}\epsilon^{\natural *}
\epsilon_*^\natural k\cong\underline{\rm Ext}_{{\cal
O}_{A^\natural }}^{g+i} (\epsilon_*^\natural k,\epsilon_*^\natural k)
$$
(le ${\cal O}_{A^\natural}$-Module inversible $\Omega_{A^\natural
/k}^{2g}$ est canoniquement trivial).  De plus,
$$
(\epsilon^\natural )^{-1}\underline{\rm Ext}_{{\cal O}_{A^\natural
}}^1 (\epsilon_*^\natural k,\epsilon_*^\natural k)
$$
n'est autre que le $k$-espace vectoriel tangent \`a l'origine de
$A^\natural$ et, pour tout entier $i$,
$$
(\epsilon^\natural )^{-1}\underline{\rm Ext}_{{\cal O}_{A^\natural
}}^i (\epsilon_*^\natural k,\epsilon_*^\natural k)
$$
est canoniquement isomorphe \`a
$$
\bigwedge^i(\epsilon^\natural )^{-1}\underline{\rm Ext}_{{\cal
O}_{A^\natural }}^1 (\epsilon_*^\natural k,\epsilon_*^\natural k)
$$
(consid\'erer une r\'esolution de Koszul de $\epsilon_*^\natural k$).
Par suite, on est ramen\'e \`a identifier $H_{\rm dR}^1(A/k)$ \`a
l'espace tangent \`a l'origine de $A^\natural$, ce qui est fait dans
[Ma-Me] (2.2).

Passons maintenant \`a l'assertion ``resp.'' du lemme. Pour montrer
que $R\widetilde{\rm pr}_*(\widetilde{\cal P},\widetilde \nabla )$
est canoniquement isomorphe \`a $\epsilon_+{\cal O}_S$ dans
$D_{\rm qcoh}^{\rm b} ({\cal D}_{A/S})$ il suffit, d'apr\`es un
th\'eor\`eme de Kashiwara (cf. [Bo] Ch. VI, Theorem 7.11), de
montrer que $R\widetilde{\rm pr}_*(\widetilde{\cal P},
\widetilde\nabla )$ est support\'e par $\epsilon (S)$ et que
$L\epsilon^* R\widetilde{\rm pr}_*\widetilde{\cal P}[-g]$ est
canoniquement isomorphe \`a ${\cal O}_S[-g]$. La seconde assertion
r\'esulte aussit\^ot du th\'eor\`eme de changement de base et du
th\'eor\`eme (2.4.1). Pour la premi\`ere assertion, on remarque que
$R\widetilde{\rm pr}_*\widetilde{\cal P}$ n'est rien d'autre que
${\cal F}'({\cal A}^\natural )$ et donc cette assertion est une
cons\'equence imm\'ediate de l'isomorphisme
$$
{\rm gr}_\bullet {\cal A}^\natural\buildrel\sim\over
\longrightarrow \pi'^*{\rm Sym}_\bullet^{{\cal O}_S}
(\epsilon^*{\cal T}_{A/S}).
$$
et du lemme (1.2.5) puisque ${\cal F}'({\cal O}_{A'})= R{\rm
pr}_*{\cal P}$.

\hfill\hfill\cqfd

\th COROLLAIRE (3.2.4)
\enonce
Le carr\'e de foncteurs
$$
\matrix{
D_{\rm qcoh}^{\rm b}({\cal O}_{A'})&\maprightover{{\cal F}'}
&D_{\rm qcoh}^{\rm b}({\cal O}_A)\cr
\noalign{\medskip}
\mapdownleft{Lp^*(\cdot )} &&\mapdownright{{\cal D}_{A/S}
\,\Lotimes_{{\cal O}_A}(\cdot )}\,\cr
\noalign{\medskip}
D_{\rm qcoh}^{\rm b}({\cal O}_{A^\natural})&
\maprightunder{\widetilde {\cal F}^\natural}& D_{\rm qcoh}^{\rm
b}({\cal D}_{A/S})\cr}
$$
est commutif \(\`a un isomorphisme canonique pr\`es\).

\hfill\hfill\cqfd
\endth

\th COROLLAIRE (3.2.5)
\enonce
Le foncteur $\widetilde{\cal F}^\natural$ envoie la
sous-cat\'egorie strictement pleine $D_{\rm coh}^{\rm b} ({\cal
O}_{A^\natural})$ de $D_{\rm qcoh}^{\rm b} ({\cal O}_{A^\natural})$
dans la sous-cat\'egorie strictement pleine $D_{\rm coh}^{\rm
b}({\cal D}_{A/S})$ de $D_{\rm qcoh}^{\rm b}({\cal D}_{A/S})$.

\hfill\hfill\cqfd
\endth
\vskip 5mm

(3.3) Tout comme la transformation de Fourier-Mukai, les
transformations $\widetilde{\cal F}$ et $\widetilde{\cal
F}^\natural$ poss\`edent des propri\'et\'es de fonctorialit\'e.
\vskip 2mm

{\pc RAPPELS} (3.3.1)\pointir Soient $T$ un $S$-sch\'ema lisse
purement de dimension relative $d_T$ et $X$ et $Y$ des
$T$-sch\'emas lisses et purement de dimensions relatives
$d_{X/T}$ et $d_{Y/T}$ respectivement.

Si $X\buildrel f\over\longrightarrow Y$ est un $T$-morphisme, on
a des foncteurs exacts
$$
f_+:D_{\rm qcoh}^{\rm b}({\cal D}_{X/T})\rightarrow D_{\rm
qcoh}^{\rm b}({\cal D}_{Y/T})
$$
et
$$
f^!:D_{\rm qcoh}^{\rm b}({\cal D}_{Y/T})\rightarrow D_{\rm
qcoh}^{\rm b}({\cal D}_{X/T})
$$
d\'efinis par
$$
f_+(\cdot )=Rf_*({\cal D}_{(Y\leftarrow X)/T}\,\Lotimes_{{\cal
D}_{X/T}}\,(\cdot ))
$$
et
$$
f^!(\cdot )={\cal D}_{(X\rightarrow Y)/T}\,\Lotimes_{f^{-1}{\cal
D}_{Y/T}}\,f^{-1}(\cdot ) [d_{X/T}-d_{Y/T}],
$$
o\`u ${\cal D}_{(Y\leftarrow X)/T}$ est le $(f^{-1}{\cal
D}_{Y/T},{\cal D}_{X/T})$-bi-Module
$$
f^{-1}({\cal D}_{Y/T}\otimes_{{\cal O}_Y}
(\Omega_{Y/T}^{d_{Y/T}})^{\otimes -1})\otimes_{f^{-1}{\cal O}_Y}
\Omega_{X/T}^{d_{X/T}}
$$
et o\`u ${\cal D}_{(X\rightarrow Y)/T}$ est le $({\cal D}_{X/T},
f^{-1}{\cal D}_{Y/T})$-bi-Module
$$
{\cal O}_X\otimes_{f^{-1}{\cal O}_Y}f^{-1}{\cal D}_{Y/T}
$$
(cf. [Bo] Ch. VI, \S 5.1 et \S 4.2). Si $f$ est lisse, on a encore
$$
f_+(\cdot )=Rf_*{\rm DR}_{X/Y}(\cdot ),
$$
o\`u le complexe de de Rham ${\rm DR}_{X/Y}({\cal M})$ d'un ${\cal
D}_{X/S}$-Module ${\cal M}$ est le complexe
$$
[{\cal M}\buildrel\nabla\over\longrightarrow \Omega_{X/Y}^1
\otimes_{{\cal O}_X}{\cal M}\buildrel\nabla\over\longrightarrow
\cdots\buildrel\nabla\over\longrightarrow\Omega_{X/Y}^{d_{X/T}
-d_{Y/T}}\otimes_{{\cal O}_X}{\cal M}]
$$
concentr\'e en degr\'es compris entre $d_{Y/T}-d_{X/T}$ et $0$ (cf.
[Bo] Ch. VI, 5.3.2).

Si $T'$ est un autre $S$-sch\'ema lisse purement de dimension
relative $d_{T'}$, si $T'\buildrel\alpha\over\longrightarrow T$
 est un morphisme de $S$-sch\'emas et si $X'$ est le
$T'$-sch\'ema d\'eduit du $T$-sch\'ema $X$ par le changement  de
base $\alpha$, on a aussi des foncteurs
$$
(\beta ,\alpha )_+:D_{\rm qcoh}^{\rm b}({\cal D}_{X'/T'})
\rightarrow D_{\rm qcoh}^{\rm b}({\cal D}_{X/T})
$$
et
$$
(\beta ,\alpha )^!:D_{\rm qcoh}^{\rm b}({\cal D}_{X/T})
\rightarrow D_{\rm qcoh}^{\rm b}({\cal D}_{X'/T'})
$$
d\'efinis par
$$
(\beta ,\alpha )_+=R\beta_*
$$
et
$$
(\beta ,\alpha )^!(\cdot )=R\beta^!(\cdot )=L\beta^*(\cdot )
\,\Lotimes_{\pi'^{-1}{\cal O}_{T'}}\,\pi'^{-1}R\alpha^! {\cal O}_T,
$$
o\`u $\beta :X'=T'\times_TX\rightarrow X$ et $\pi' :
X'=T'\times_TX\rightarrow T'$ sont les projections canoniques.

On a un produit tensoriel externe
$$
\boxtimes : D_{\rm qcoh}^{\rm b}({\cal D}_{X/T})\times D_{\rm
qcoh}^{\rm b}({\cal D}_{Y/T})\rightarrow D_{\rm qcoh}^{\rm
b}({\cal D}_{X\times_SY/T\times_ST})
$$
et un produit tensoriel interne
$$
\otimes^!: D_{\rm qcoh}^{\rm b}({\cal D}_{X/T})\times D_{\rm
qcoh}^{\rm b}({\cal D}_{X/T})\rightarrow D_{\rm qcoh}^{\rm
b}({\cal D}_{X/T})
$$
d\'efini par
$$
(\cdot )\otimes^!(\cdot )=\Delta_{X/T}^! (\widetilde
\Delta_T,\Delta_T)^!((\cdot )\boxtimes (\cdot ))=((\cdot)\,
\Lotimes_{{\cal O}_X}\,(\cdot ))\,\Lotimes_{\pi^{-1}({\cal O}_T)}\,
\pi^{-1}(\Omega_{T/S}^{d_T})^{\otimes -1}[-d_{X/T}-d_T],
$$
o\`u $\Delta_T:T \hookrightarrow T\times_ST$, $\widetilde
\Delta_T:X\times_TX \hookrightarrow X\times_SX$ et
$\Delta_{X/T}:X\hookrightarrow X\times_TX$ sont les immersions
diagonales (comparer \`a [Bo] Ch. VIII, \S 14.3).

On a enfin un foncteur de dualit\'e
$$
D:D_{\rm coh}^{\rm b}({\cal D}_{X/T})^{\rm opp}\rightarrow D_{\rm
coh}^{\rm b}({\cal D}_{X/T})
$$
d\'efini par
$$
D(\cdot )=R\underline{\rm Hom}_{{\cal D}_{X/T}}(\,\cdot\,,{\cal
D}_{X/T}\otimes_{{\cal O}_X} (\Omega_{X/T}^{d_{X/T}})^{\otimes
-1})\otimes_{\pi^{-1}({\cal O}_T)}\pi^{-1}\Omega_{T/S}^{d_T}
[d_{X/T}+d_T]
$$
(cf. [Bo] Ch. VI, \S 3.6).

\hfill\hfill\cqfd
\vskip 2mm

Compte tenu de ces rappels, on a
$$
\widetilde{\cal F}(\cdot )=\widetilde{\rm pr}_+^\natural
\bigl((\widetilde{\cal P},\widetilde{\nabla})\otimes^!
(\widetilde{\rm pr},\pi^\natural )^! (\cdot )\bigr)[g]
$$
et
$$
\widetilde{\cal F}^\natural (\cdot )=\bigl(\widetilde{\rm
pr},\pi^\natural )_+ ((\widetilde{\cal P},\widetilde{\nabla})
\otimes^!\widetilde{\rm pr}^{\natural !}(\cdot )\bigr)[2g].
$$

\th PROPOSITION (3.3.2)
\enonce
Soit $f:A_1\rightarrow A_2$ un morphisme $S$-sch\'emas
ab\'eliens, avec $A_i$ purement de dimension relative $g_i$
\($i=1,2$\), et soit $f^\natural :A_2^\natural\rightarrow
A_1^\natural$ le morphisme transpos\'e. Alors on a les couples
d'isomorphismes canoniques de foncteurs
$$
\left\{\matrix{
\widetilde{\cal F}_2\circ f_+\cong Rf^{\natural
!}\circ\widetilde{\cal F}_1(\cdot ) [2(g-1-g_2)],\hfill\cr
\noalign{\smallskip}
f_+\circ\widetilde{\cal F}_1^\natural\cong\widetilde {\cal
F}_2^\natural\circ Rf^{\natural !}(\cdot ) [g-1-g_2],\hfill\cr}\right.
$$
avec $Rf^{\natural !}(\cdot )=Lf^{\natural *}(\cdot ) [2(g_2-g_1)]$,
et
$$
\left\{\matrix{
\widetilde{\cal F}_1^\natural\circ Rf_*^\natural (\cdot )
[g-1-g_2]\cong f^!\circ\widetilde{\cal F}_2^\natural,\hfill\cr
\noalign{\smallskip}
Rf_*^\natural\circ\widetilde{\cal F}_2\cong\widetilde{\cal
F}_1\circ f^!,\hfill\cr}\right.
$$
o\`u $\widetilde{\cal F}_i$ et $\widetilde{\cal F}_i^\natural$ sont
les foncteurs $\widetilde{\cal F}$ et $\widetilde{\cal F}^\natural$
pour le $S$-sch\'ema ab\'elien $A_i$ \($i=1,2$\).
\endth

\rem Preuve
\endrem
Les premiers isomorphismes de chaque couple r\'esultent
formellement du th\'eor\`eme de changement de base, de la formule
des projections et de l'isomorphisme canonique
$$
(f^\natural\times {\rm id}_{A_1},f^\natural )^! (\widetilde{\cal
P}_1,\widetilde\nabla_1)[3(g-1-g_2)]\cong ({\rm
id}_{A_2^\natural}\times f)^! (\widetilde{\cal P}_2,
\widetilde\nabla_2).
$$
Quant aux seconds, ils se d\'eduisent des premiers par le
th\'eor\`eme (3.2.1).

\hfill\hfill\cqfd

On appelle {\it produits de convolution} et on note par
$$
(\cdot )*_{\cal D}(\cdot ):D_{\rm qcoh}^{\rm b}({\cal D}_{A/S})
\times D_{\rm qcoh}^{\rm b}({\cal D}_{A/S})\rightarrow D_{\rm
qcoh}^{\rm b}({\cal D}_{A/S})
$$
et
$$
(\cdot )*^\natural (\cdot ):D_{\rm qcoh}^{\rm b}({\cal
O}_{A^\natural})\times D_{\rm qcoh}^{\rm b}({\cal
O}_{A^\natural})\rightarrow D_{\rm qcoh}^{\rm b}({\cal
O}_{A^\natural})
$$
les foncteurs d\'efinis par
$$
(\cdot )*_{\cal D}(\cdot )=\mu_+((\cdot )\boxtimes (\cdot ))
$$
et
$$
(\cdot )*^\natural (\cdot )=R\mu_*^\natural ((\cdot )
\boxtimes (\cdot ))
$$
(on rappelle que $\mu :A\times_SA\rightarrow A$ et
$\mu^\natural :A^\natural \times_SA^\natural \rightarrow
A^\natural $ sont les lois de groupe).

\th COROLLAIRE (3.3.3)
\enonce
On a des couples d'isomorphismes canoniques de foncteurs
$$
\left\{\matrix{
\widetilde{\cal F}((\cdot )*_{\cal D}(\cdot ))\cong\widetilde{\cal
F}(\cdot )\,\Lotimes_{{\cal O}_{A^\natural}}\,\widetilde{\cal
F}(\cdot ),\hfill\cr
\noalign{\medskip}
\widetilde{\cal F}^\natural(\cdot )*_{\cal D}\widetilde{\cal
F}^\natural(\cdot )\cong\widetilde{\cal F}^\natural ((\cdot )
\,\Lotimes_{{\cal O}_{A^\natural}}\,(\cdot ))[-g],\hfill\cr}\right.
$$
et
$$
\left\{\matrix{
\widetilde{\cal F}^\natural ((\cdot )*^\natural (\cdot ))\cong
\widetilde{\cal F}^\natural (\cdot )[g]\otimes_S^! \widetilde{\cal
F}^\natural (\cdot ),\hfill\cr
\noalign{\medskip}
\widetilde{\cal F}(\cdot )*^\natural\widetilde{\cal F}(\cdot) \cong
\widetilde{\cal F}((\cdot )\otimes_S^!(\cdot )).\hfill\cr}\right.
$$

\hfill\hfill\cqfd
\endth

\th PROPOSITION (3.3.4)
\enonce
Consid\'erons les  foncteurs de dualit\'e $D:D_{\rm coh}^{\rm
b}({\cal D}_{A/S})^{\rm opp}\rightarrow D_{\rm coh}^{\rm b}({\cal
D}_{A/S})$ \(pour $(X,T)=(A,S)$\) et $D^\natural : D_{\rm coh}^{\rm
b}({\cal O}_{A^\natural})^{\rm opp}\rightarrow D_{\rm coh}^{\rm
b}({\cal O}_{A^\natural})$ \(pour $(X,T)=(A^\natural ,A^\natural )$\)
d\'efinis en $(3.3.1)$. Alors, on a des isomorphismes canoniques de
foncteurs
$$
D^\natural\circ\widetilde{\cal F}\cong\langle -1 \rangle^{\natural
*}\circ\widetilde{\cal F}\circ D(\cdot )[2g]
$$
et
$$
D\circ\widetilde{\cal F}^\natural\cong\langle
-1\rangle^*\circ\widetilde{\cal F}^\natural\circ D^\natural .
$$
\endth

\rem Preuve
\endrem
Les seconds isomorphismes se d\'eduisent des premiers
\`a l'aide du th\'eor\`eme (3.2.1). Quant aux premiers, ils r\'esultent
des isomorphismes suivants: avec les notations de (3.3.1), on a
$$
D\circ f_+\cong f_+\circ D
$$
si $f$ est propre, on a
$$
D(({\cal L},\nabla )\otimes^!(\cdot ))\cong ({\cal L},\nabla
)^{\otimes -1}\otimes^!D(\cdot )\otimes_{\pi^{-1}{\cal O}_T}\pi^{-1}
(\Omega_{T/S}^{d_T})^{\otimes 2}[2d_{X/T}+2d_T]
$$
si $({\cal L},\nabla )$ est un ${\cal O}_X$-Module inversible \`a
connexion int\'egrable relative \`a $T$ et on a
$$
(\beta ,\alpha )^!(D(\cdot ))\otimes_{{\cal O}_{X'}}
(\Omega_{X'/X}^n)^{\otimes -1}[-n]\cong D'\circ (\beta ,\alpha
)^!(\cdot )
$$
si $\alpha$ est lisse, purement de dimension relative $n$, o\`u
$D'$ est le foncteur de dualit\'e pour $X'/T'$.

\hfill\hfill\cqfd

\vskip 5mm
{\bf 4. Groupes alg\'ebriques et groupes formels sur un corps de
caract\'eristique nulle: rappels.}
\vskip 5mm

(4.1) Dans toute la suite de cet article, on fixe un corps $k$ de
caract\'eristique nulle et une cl\^oture alg\'ebrique $\overline k$ de
$k$.

On appellera simplement $k$-{\it groupes alg\'egriques} les
$k$-sch\'emas en groupes ab\'eliens qui sont s\'epar\'es et de type
fini en tant que $k$-sch\'emas. Puisque $k$ est de caract\'eristique
nulle, le sch\'ema sous-jacent  \`a tout $k$-groupe alg\'ebrique est
lisse et de dimension pure sur $k$ (dimension not\'ee dans la suite
$d_G$) et on a:

\th PROPOSITION (4.1.1)
\enonce
{\rm (i)} La cat\'egorie des $k$-groupes alg\'ebriques est
ab\'elienne.

\decale{\rm (ii)} Tout $k$-groupe alg\'ebrique $G$ admet le
d\'evissage canonique
$$
\matrix{
1\rightarrow G^0\rightarrow G\rightarrow G/G^0\rightarrow 1,\cr
\noalign{\medskip}
1\rightarrow L\rightarrow G^0\rightarrow A\rightarrow 0,\cr
\noalign{\medskip}
1\rightarrow T\rightarrow L\rightarrow V\rightarrow 0,\cr}
$$
o\`u:
\decale{$\bullet$} $G^0$ est la composante neutre de $G$ et
$G/G^0$ est fini et \'etale sur $k$,
\decale{$\bullet$} $L$ est le plus petit sous-groupe alg\'ebrique de
$G^0$ tel que $G^0/L$ soit propre sur $k$, $L$ est affine et $A$ est
un $k$-sch\'ema ab\'elien,
\decale{$\bullet$} $T$ est le plus grand $k$-tore contenu dans $L$
et $V$ est un $k$-vectoriel.
\vskip 2mm
De plus, $L$ est aussi le plus grand sous-$k$-groupe alg\'ebrique
affine et connexe de $G^0$ et l'extension de $V$ par $T$ est
canoniquement scind\'ee, de sorte que $L$ est canoniquement
isomorphe au produit direct $T\times_kV$.
\endth

\rem Preuve
\endrem
Voir [SGA3] Expos\'e VI$_{\rm A}$ et [Se] Ch. III, \S 7, Prop. 11 et
Prop. 12.
\hfill\hfill\cqfd

Bien s\^ur, on a
$$
d_G=d_{G^0}=d_L+d_A=d_T+d_V+d_A
$$
et on a des isomorphismes non canoniques
$$
T\otimes_kk'\cong ({\bb G}_{{\rm m},k'})^{d_T}
$$
et
$$
V\cong ({\bb G}_{{\rm a},k})^{d_V}
$$
de groupes alg\'ebriques sur $k'$ et $k$ respectivement, o\`u
$k'$ est une extension finie de $k$ contenue dans $\overline k$
convenable.
\vskip 5mm

(4.2) On appellera simplement $k$-{\it groupes formels} les
$k$-sch\'emas formels affines en groupes ab\'eliens ${\cal G}$ (on
note additivement la loi de groupe) tels que l'action naturelle de
${\rm Gal}(\overline k/k)$ sur ${\cal G}(\overline k)$ se factorise
\`a travers le quotient ${\rm Gal}(k'/k)$ pour une extension finie
$k'$ de $k$ dans $\overline k$, que ${\cal G}(\overline k)={\cal
G}(k')$ soit un groupe ab\'elien de type fini et que la $k$-alg\`ebre
topologique ${\cal O}_{{\cal G},0}$ soit un quotient d'une
$k$-alg\`ebre de s\'eries formelles en un nombre fini
d'ind\'etermin\'ees (i.e. topologiquement de type fini). Puisque
$k$ est de caract\'eristique nulle, tout $k$-groupe formel est
formellement lisse et de dimension pure $d_{\cal G}$ pour un
certain entier $d_{\cal G}\geq 0$, i.e. ${\cal O}_{{\cal G},0}$ est
une $k$-alg\`ebre de s\'eries formelles en $d_{\cal G}$
ind\'etermin\'ees, et on a:

\th PROPOSITION (4.2.1)
\enonce
{\rm (i)} La cat\'egorie des $k$-groupes formels est ab\'elienne.
\decale{\rm (ii)} Tout $k$-groupe formel ${\cal G}$ admet le
d\'evissage canonique suivant:
$$
\matrix{
0\rightarrow {\cal G}^0\rightarrow {\cal G}\rightarrow {\cal
G}^{\rm \acute et}\rightarrow 0,\cr
\noalign{\medskip}
0\rightarrow {\cal G}_{\rm tors}^{\rm \acute et}\rightarrow {\cal
G}^{\rm \acute et}\rightarrow {\cal G}_{\rm lib}^{\rm \acute et}
\rightarrow 0,\cr}
$$
o\`u

\decale{$\bullet$} ${\cal G}^0$ est la composante neutre de ${\cal
G}$ et ${\cal G}^{\rm \acute et}$ est \'etale sur $k$,

\decale{$\bullet$} ${\cal G}_{\rm tors}^{\rm \acute et}$ est le plus
grand sous-$k$-sch\'ema en groupes de ${\cal G}^{\rm \acute et}$
dont le $k$-sch\'ema sous-jacent est fini et \'etale et ${\cal
G}_{\rm lib}^{\rm \acute et}(k')$ est un groupe ab\'elien libre de
rang fini.
\vskip 2mm

De plus, l'extension de ${\cal G}^{\rm \acute et}$ par ${\cal G}^0$
est canoniquement scind\'ee, de sorte que ${\cal G}$ est
canoniquement isomorphe au produit direct ${\cal G}^0
\times_k{\cal G}^{\rm \acute et}$.
\endth

\rem Preuve
\endrem
Voir [Fo] Ch. I, 6.6 et \S 7.
\hfill\hfill\cqfd

Bien s\^ur, on a
$$
d_{\cal G}=d_{{\cal G}^0}
$$
et, si on note $r_{\cal G}$ le rang du groupe ab\'elien de type fini
${\cal G}(k')$, i.e. $r_{\cal G}={\rm dim}_{\bb Q} ({\bb Q}\otimes
{\cal G}(k'))$, on a
$$
r_{\cal G}=r_{{\cal G}^{\rm \acute et}}= r_{{\cal G}_{\rm lib}^{\rm
\acute et}}.
$$
En outre, on a des isomorphismes non canoniques de $k$-groupes
formels
$$
{\cal G}^0\cong (\widehat {\bb G}_{{\rm a}, k}^0)^{d_{\cal G}}
$$
et de groupes ab\'eliens
$$
{\cal G}_{\rm lib}^{\rm \acute et}(k')\cong {\bb Z}^{r_{\cal G}},
$$
o\`u $\widehat {\bb G}_{{\rm a}, k}^0$ est le $k$-groupe formel
additif connexe (i.e. ${\rm Spf}(k[[x]])$ avec co-multiplication
$x\mapsto x\otimes 1+ 1\otimes x$).
\vskip 5mm

(4.3) Soit ${\bf Aff}/k$ la cat\'egorie des $k$-sch\'emas affines
munie de la topologie fppf. On notera ${\bf Ens}/k$ (resp. ${\bf
Ab}/k$) la cat\'egorie des faisceaux d'ensembles (resp. de groupes
ab\'eliens) sur le site ainsi d\'efini.

On identifie la cat\'egorie des $k$-sch\'emas \`a une
sous-cat\'egorie pleine de ${\bf Ens}/k$ de la mani\`ere habituelle:
au $k$-sch\'ema $X$ on associe le faisceau de valeur $X(R)$ sur
${\rm Spec}(R)$. En particulier, on identifiera la cat\'egorie des
$k$-groupes alg\'ebriques \`a une sous-cat\'egorie pleine de ${\bf
Ab}/k$, sous-cat\'egorie que l'on notera dans la suite ${\bf
Gr.alg}/k$.

De m\^eme on identifie la cat\'egorie des $k$-sch\'emas formels
affines \`a une sous-cat\'egorie pleine de ${\bf Ens}/k$: au
$k$-sch\'ema formel affine ${\cal X}={\rm Spf}(A)$, on associe le
faisceau  de valeur sur tout $k$-sch\'ema affine ${\rm Spec}(R)$
l'ensemble ${\cal X}(R)$ des homomorphismes continus de
$k$-alg\`ebres de $A$ dans $R$, la $k$-alg\`ebre $R$ \'etant munie
de la topologie discr\`ete. En particulier, on identifiera la
cat\'egorie des $k$-groupes formels \`a une sous-cat\'egorie pleine
${\bf Gr.form}/k$ de ${\bf Ab}/k$.

\th PROPOSITION (4.3.1)
\enonce
{\rm (i)} Si $f$ est une fl\`eche dans ${\bf Gr.alg}/k$ \(resp. ${\bf
Gr.form}/k$\), les noyau et conoyau de $f$ calcul\'es dans ${\bf
Gr.alg}/k$ \(resp. ${\bf Gr.form}/k$\) sont aussi des noyau et
conoyau de $f$ calcul\'es dans ${\bf Ab}/k$.

\decale{\rm (ii)} Si
$$
0\rightarrow F'\rightarrow F\rightarrow F''\rightarrow 0
$$
est une suite exacte dans ${\bf Ab}/k$ et si $F'$ et $F''$ sont
isomorphes \`a des objets de ${\bf Gr.alg}/k$ \(resp. ${\bf
Gr.form}/k$\), il en est de m\^eme de $F$.

\hfill\hfill\cqfd
\endth

\vskip 5mm

{\bf 5. La cat\'egorie des $1$-motifs g\'en\'eralis\'es; dualit\'e de
Cartier.}
\vskip 5mm

(5.1) On garde les notations du num\'ero pr\'ec\'edent. On  note
$C^{[-1,0]}({\bf Ab}/k)$ la cat\'egorie ab\'elienne des complexes
dans ${\bf Ab}/k$ concentr\'es en degr\'es $-1$ et $0$.

\th D\'EFINITION (5.1.1)
\enonce
Un $1$-motif g\'en\'eralis\'e sur $k$ est un objet $M$ de
$C^{[-1,0]}({\bf Ab}/k)$ de la forme $M=[{\cal G}\buildrel u\over
\longrightarrow G]$, o\`u ${\cal G}$ est un $k$-groupe formel sans
torsion \(${\cal G}_{\rm tors}^{\rm \acute et}=(0)$\) et  o\`u $G$
est un $k$-groupe alg\'ebrique connexe \($G=G^0$\).
\endth

La cat\'egorie ${\bf Mot}/k$ des $1$-motifs g\'en\'eralis\'es sur
$k$ est la sous-cat\'egorie pleine de $C^{[-1,0]}({\bf Ab}/k)$ dont
les objets sont les $1$-motifs g\'en\'eralis\'es sur $k$.

\vskip 2mm
{\pc REMARQUE} (5.1.2)\pointir Un $1$-motif sur $k$ au sens  de
[De] (10.1) est un $1$-motif g\'en\'eralis\'e sur $k$ tel que
${\cal G}^0=(0)$, ${\cal G}^{\rm \acute et}= {\cal G}_{\rm lib}^{\rm
\acute et}$ est constant (de valeur un ${\bb Z}$-module libre not\'e
$X$ dans loc. cit.) et o\`u $G$ est une extension d'un $k$-sch\'ema
ab\'elien $A$ par un $k$-tore d\'eploy\'e $T$.

\hfill\hfill\cqfd

\th PROPOSITION (5.1.3)
\enonce
La cat\'egorie ${\bf Mot}/k$ est une cat\'egorie exacte dans
laquelle toute fl\`eche admet un noyau et un conoyau.

Les suites exactes dans ${\bf Mot}/k$ sont les suites
$$
0\rightarrow M'\rightarrow M\rightarrow M''\rightarrow 0
$$
telles que les suites correspondantes
$$
0\rightarrow G'\rightarrow G\rightarrow G''\rightarrow 0
$$
dans ${\bf Gr.alg}/k$ et
$$
0\rightarrow {\cal G}'\rightarrow {\cal G}\rightarrow {\cal G}''
\rightarrow 0
$$
dans ${\bf Gr.form}/k$ soient toutes les deux exactes.
\endth

\rem Preuve
\endrem
Soit
$$
f=(f_G,f_{\cal G}):M_1=[{\cal G}_1\buildrel u_1\over
\longrightarrow G_1]\rightarrow [{\cal G}_2\buildrel u_2\over
\longrightarrow G_2]=M_2
$$
un morphisme dans ${\bf Mot}/k$. Posons
$$
\left\{\matrix{G_1'={\rm Ker}(f_G)^0,\hfill\cr
\noalign{\medskip}
{\cal G}_1'={\rm Ker}(f_{\cal G})\cap u_1^{-1}(G_1'),\hfill\cr
\noalign{\medskip}
u_1'=u_1|{\cal G}_1':{\cal G}_1'\rightarrow G_1',\hfill\cr}\right.
$$
et
$$
\left\{\matrix{G_2''={\rm Coker}(f_G)/\overline u_2 ({\rm Coker}
(f_{\cal G})_{\rm tors}^{\rm \acute et}),\hfill\cr
\noalign{\medskip}
{\cal G}_2''={\rm Coker}(f_{\cal G})/{\rm Coker}(f_{\cal G})_{\rm
tors}^{\rm \acute et},\hfill\cr
\noalign{\medskip}
u_2'':{\cal G}_2''\rightarrow G_2''\hbox{ induit par }\overline
u_2,\hfill\cr}\right.
$$
o\`u $\overline u_2:{\rm Coker}(f_{\cal G})\rightarrow {\rm
Coker}(f_G)$ est le morphisme induit par $u_2$. Alors
$M_1'=[{\cal G}_1'\buildrel u_1'\over\longrightarrow G_1']$ est un
noyau de $f$ et $M_2''=[{\cal G}_2''\buildrel u_2''\over
\longrightarrow G_2'']$ est un conoyau de $f$.

En outre, $f$ est un monomorphisme strict (resp. un \'epimorphisme
strict) si et seulement si $f_G$ et $f_{\cal G}$ sont des
monomorphismes (resp des \'epimorphismes) et  si ${\rm
Coker}(f_{\cal G})$ est sans torsion (resp. ${\rm Ker}(f_G)$ est
connexe).

Comme la classe des $k$-groupes alg\'ebriques connexes et  celle
des $k$-groupes formels sans torsion sont stables par extension
dans ${\bf Ab}/k$ (cf. (4.3.1)(ii)), on v\'erifie  facilement que:

\decale{$\bullet$} la classe des monomorphismes stricts (resp.
\'epimorphismes stricts) est stable par composition,

\decale{$\bullet$} tout ``push-out'' (resp. ``pull-back'') d'un
monomorphisme strict (resp. \'epimor-\break phisme strict) est
repr\'esentable dans ${\bf Mot}/k$ par un monomorphisme strict
(resp. \'epimorphisme strict),

\decale{$\bullet$} si un morphisme compos\'e de ${\bf Mot}/k$,
$$
M\buildrel f\over\longrightarrow M'\buildrel f'\over
\longrightarrow M''
$$
est un monomorphisme strict (resp. un \'epimorphisme strict), il en
est de m\^eme de $f'$ (resp. $f''$).

Ceci ach\`eve la preuve de la proposition.

\hfill\hfill\cqfd
\vskip 5mm

(5.2) Rappelons que la dualit\'e de Cartier (cf. [SGA3] VII$_B$,
(2.2.2)) induit une anti-\'equivalence de cat\'egories ab\'eliennes
entre ${\bf Gr.alg}/k$ et ${\bf Gr.form}/k$.

Plus pr\'ecis\'ement, si $G$ (resp. ${\cal G}$) est un $k$-groupe
alg\'ebrique affine (resp. un $k$-groupe formel), le faisceau de
groupes ab\'eliens sur ${\bf Ab}/k$
$$
\underline{\rm Hom}_{{\bf Ab}/k}(G,{\bb G}_{{\rm m},k})
$$
(resp.
$$
\underline{\rm Hom}_{{\bf Ab}/k}({\cal G},{\bb G}_{{\rm m},k})\,)
$$
est repr\'esentable par un $k$-groupe formel ${\cal G}'$ (resp. un
$k$-groupe alg\'ebrique affine $G'$), appel\'e le dual de Cartier de
$G$ (resp. ${\cal G}$). Les foncteurs contravariants $G\mapsto
{\cal G}'$ et ${\cal G}\mapsto G'$ sont exacts et quasi-inverses l'un
de l'autre. Si $G$ est connexe (resp. ${\cal G}$ est sans torsion),
${\cal G}'$ est sans torsion (resp. $G'$ est connexe). Si $G$ est un
$k$-tore (resp. ${\cal G}$ est \'etale sans torsion), ${\cal G}'$ est
\'etale sans torsion (resp. $G'$ est un $k$-tore). Si $G$ est un
$k$-vectoriel (resp. ${\cal G}$ est connexe), ${\cal G}'$ est connexe
(resp. $G'$ est un $k$-vectoriel). Enfin, on a
$$
d_G=d_{{\cal G}'}+r_{{\cal G}'}
$$
(resp.
$$
d_{\cal G}+r_{\cal G}=d_{G'}\,).
$$

On se propose maintenant de prolonger la dualit\'e de Cartier
ci-dessus et la dualit\'e pour les $k$-sch\'emas ab\'eliens
consid\'er\'ee dans la section 1 en une dualit\'e pour les
$1$-motifs g\'en\'eralis\'es. On veut bien entendu que cette
dualit\'e prolonge aussi la dualit\'e pour les $1$-motifs d\'efinie
par Deligne dans [De] (10.2.11).

Commen\c cons par montrer le lemme suivant.

\th LEMME (5.2.1)
\enonce
Soient $A$ un $k$-sch\'ema ab\'elien de dual $A'$ et ${\cal G}$ un
$k$-groupe formel sans torsion de dual de Cartier le $k$-groupe
alg\'ebrique affine connexe $L'$. Alors, les classes d'extensions
$G'$ de $A'$ par $L'$ \(dans ${\bf Ab}/k$ ou ce qui revient au
m\^eme dans ${\bf Gr.alg}/k$\) sont en bijection avec les fl\`eches
${\cal G}\rightarrow A$ dans ${\bf Ab}/k$.
\endth

\rem Preuve
\endrem
Comme on a les d\'ecompositions ${\cal G}={\cal G}^0\times_k
{\cal G}^{\rm \acute et}$ et $L'=T'\times_kV'$ en dualit\'e, il suffit
de traiter s\'epar\'ement les cas ${\cal G}$ \'etale et ${\cal G}$
connexe.

Si ${\cal G}$ est \'etale, $L'=T'$ est un $k$-tore. Par descente
galoisienne et par additivit\'e, on est ramen\'e au cas ${\cal G}
(\overline k)={\bb Z}$ avec action triviale de ${\rm Gal} (\overline
k/k)$ et $T'={\bb G}_{{\rm m},k}$. Alors la bijection cherch\'ee
associe la classe d'extension $\iota (a)\in A''(k)$ de $A'$ par ${\bb
G}_{{\rm m},k}$ \`a la fl\`eche ${\cal G}\rightarrow A$ qui envoie
$1\in {\bb Z}={\cal G}(k)$ sur $a\in A(k)$.

Si ${\cal G}$ est connexe, $L'=V'$ est aussi le $k$-vectoriel dual du
$k$-vectoriel ${\rm Lie}({\cal G})$ (alg\`ebre de Lie de ${\cal G}$).
Par suite la bijection cherch\'ee associe \`a une  fl\`eche $u:{\cal
G}\rightarrow A$ l'extension de $A'$ par $V'$ d\'eduite de
l'extension vectorielle universelle $A^\natural$ par ``push-out'' via
la fl\`eche $\Omega_A\rightarrow V'$ transpos\'ee de la fl\`eche
${\rm Lie}(u): {\rm Lie}({\cal G})\rightarrow {\rm Lie}(A)=
\epsilon^*{\cal T}_{A/k}$.

\hfill\hfill\cqfd

Soit maintenant $M=[{\cal G}\buildrel u\over\longrightarrow G]$ un
$1$-motif g\'en\'eralis\'e sur $k$. Le $k$-groupe alg\'ebrique $G$
admet le d\'evissage canonique
$$
1\rightarrow L\rightarrow G\rightarrow A\rightarrow 0
$$
o\`u $L$ est affine et $A$ est un $k$-sch\'ema ab\'elien (cf.
(4.1.1)(ii)). On note ${\cal G}\buildrel\overline u\over
\longrightarrow A$ la fl\`eche compos\'ee de $u$ et de la
projection $G\twoheadrightarrow A$. Soit ${\cal G}'$ le dual de
Cartier de $L$, $L'$ celui de ${\cal G}$ et $A'$ le $k$-sch\'ema
ab\'elien dual de $A$.

D'apr\`es le lemme ci-dessus, la fl\`eche $\overline u$ d\'efinit une
extension $G'$ de $A'$ par $L'$ et l'extension $G$ de $A$ par $L$
d\'efinit une fl\`eche ${\cal G}'\buildrel\overline
u'\over\longrightarrow A'$.

\th PROPOSITION (5.2.2)
\enonce
Soient ${\cal G}$, ${\cal G}'$, $G$, $G'$, $\overline u$ et $\overline
u'$ comme ci-dessus. Alors, la donn\'ee d'une fl\`eche ${\cal
G}\buildrel u\over\longrightarrow G$ relevant $\overline u$
\'equivaut \`a la donn\'ee d'une fl\`eche ${\cal G}'\buildrel
u'\over\longrightarrow G'$ relevant $\overline u'$.
\endth

\rem Preuve
\endrem
Soient $E$ l'extension de ${\cal G}$ par $L$ ``pull-back'' de
l'extension $G$ de $A$ par $L$ via $\overline u : {\cal G}\rightarrow
A$ et $E'$ l'extension de ${\cal G}'$ par $L'$ ``pull-back'' de
l'extension $G'$ de $A'$ par $L'$ via $\overline u':{\cal G}'
\rightarrow A'$. Les donn\'ees de $u$ et $u'$ sont \'equivalentes
\`a des scindages de $E$ et $E'$ respectivement. Il ne reste plus
qu'\`a remarquer que  les suites exactes
$$
1\rightarrow L\rightarrow E\rightarrow {\cal G}\rightarrow 0
$$
et
$$
1\rightarrow L'\rightarrow E'\rightarrow {\cal G}'\rightarrow 0
$$
dans ${\bf Ab}/k$ se d\'eduisent l'une de l'autre par application du
foncteur $\underline{\rm Hom}_{{\bf Ab}/k} (\,\cdot\,,{\bb
G}_{{\rm m},k})$, de sorte que scinder l'une  de ces suites exactes
\'equivaut \`a scinder l'autre.

\hfill\hfill\cqfd

\th D\'EFINITION (5.2.3)
\enonce
Le dual du $1$-motif g\'en\'eralis\'e $M$ sur $k$ est le $1$-motif
g\'en\'eralis\'e
$$
M'=[{\cal G}'\buildrel u'\over\longrightarrow G']
$$
construit ci-dessus.
\endth

Soit ${\cal Q}$ le ${\cal O}_{G'\times_kG}$-Module inversible
image r\'eciproque par la projection canonique $G'\times_kG
\twoheadrightarrow A'\times_kA$ du Module de Poincar\'e ${\cal
P}$. Alors ${\cal Q}$ est trivialis\'e le long de $(G'\times_k\{0\})
\cup (\{0'\}\times_kG)$ et satisfait le th\'eor\`eme du carr\'e
aussi bien pour le $G'$-sch\'ema en groupes $G'\times_kG$ que pour
le $G$-sch\'ema en groupes $G'\times_kG$. On laisse le soin au
lecteur de v\'erifier que:

\decale{$\bullet$} $({\rm id}_{G'}\times u)^*{\cal Q}$ (resp.
$(u'\times {\rm id}_G)^*{\cal Q}$) est canoniquement trivialis\'e en
tant que ${\cal O}_{G'\times_k{\cal G}}$-Module (resp. ${\cal
O}_{{\cal G'}\times_kG}$-Module) inversible rigidifi\'e et
satisfaisant le th\'eor\`eme du carr\'e pour le $G'$-sch\'ema
formel en groupes $G'\times_k{\cal G}$ (resp. pour le
$G$-sch\'ema formel en groupes ${\cal G'}\times_kG$),

\decale{$\bullet$} ces trivialisations de $({\rm id}_{G'}\times
u)^*{\cal Q}$ et $(u'\times {\rm id}_G)^*{\cal Q}$ induisent la
m\^eme trivilisation de
$$
(u'\times {\rm id}_{\cal G})^*({\rm id}_{G'}\times u)^*{\cal Q} \cong
({\rm id}_{{\cal G}'}\times u)^*(u'\times {\rm id}_G)^* {\cal Q}.
$$
\vskip 2mm

En d'autres termes, le ${\bb G}_{{\rm m},k}$-torseur sur
$G'\times_kG$ associ\'e \`a ${\cal Q}$ est muni d'une structure de
bi-extension de $M$ et $M'$ par ${\bb G}_{{\rm m},k}$ au sens de
Deligne (cf. [De] (10.2.1)). Le ${\cal O}_{G'\times_kG}$-Module ${\cal
Q}$ muni des structures suppl\'ementaires que l'on vient d'expliciter
sera encore dit  {\it de Poincar\'e}.

Il r\'esulte du caract\`ere sym\'etrique de la construction de $M'$
que l'on a un isomorphisme de bi-dualit\'e dans ${\bf Mot}/k$,
fonctoriel en $M$,
$$
\iota :M\buildrel\sim\over\longrightarrow M'',\leqno (5.2.4)
$$
gr\^ace auquel on peut identifier le Module de Poincar\'e ${\cal Q}'$
pour $M'$ au Module de Poincar\'e ${\cal Q}$ pour $M$ (\`a
permutation des facteurs pr\`es).

On laisse au lecteur le soin  de donner une description autoduale de
la cat\'egorie ${\bf Mot}/k$ qui g\'en\'eralise [De] (10.2.13).

Terminons cette section avec quelques exemples de paires de
$1$-motifs g\'en\'eralis\'es sur $k$ en dualit\'e.
\vskip 2mm

{\pc EXEMPLES} (5.2.5)\pointir {\rm (i)}
Soient $V$ et $W$ deux $k$-vectoriels et $W\buildrel u\over
\longrightarrow V$ un homomorphisme de $k$-vectoriels, de
transpos\'e $V'\buildrel u'\over\longrightarrow W'$. Notons
$\widehat W^0$ et $\widehat V'^0$ les compl\'et\'es formels \`a
l'origine de $W$ et $V'$ et notons encore $\widehat W^0\buildrel
u\over\longrightarrow V$ et $\widehat V'^0\buildrel u'\over
\longrightarrow W'$ les restrictions de $u$ et $u'$ \`a $\widehat
W^0$ et $\widehat V'^0$ respectivement. Alors, $M= [\widehat
W^0\buildrel u\over\longrightarrow V]$ et $M'=[\widehat
V'^0\buildrel u'\over\longrightarrow W']$ sont deux $1$-motifs
g\'en\'eralis\'es en dualit\'e.

\decale{\rm (ii)}
Soit $T$ un $k$-tore de dual de Cartier le $k$-groupe formel
\'etale sans torsion $X$. Notons simplement $\Omega_T$ le
$k$-vectoriel des $1$-formes diff\'erentielles (relatives
\`a $k$) invariantes sur $T$, i.e.
$$
\Omega_T=(\Omega_{T/k}^1)_{(1)},
$$
et notons $\widehat T^0$ le compl\'et\'e \`a l'origine de $T$ et
$\widehat T^0\buildrel {\rm can}\over\longrightarrow T$ la
fl\`eche canonique. Alors, $M=[\widehat T^0\buildrel {\rm can}
\over\longrightarrow T]$ et $M'=[X\buildrel {\rm can}'\over
\longrightarrow \Omega_T]$, o\`u
$$
{\rm can}'(x)=x^*({dt\over t})\qquad (\forall x\in X),
$$
sont deux $1$-motifs g\'en\'eralis\'es en dualit\'e (${dt/t}$ est la
$1$-forme diff\'erentielle invariante standard sur ${\bb G}_{{\rm
m},k}={\rm Spec}(k[t,t^{-1}])$).

\decale{\rm (iii)}
Soient $A$ un $k$-sch\'ema ab\'elien et $A^\natural$ l'extension
vectorielle universelle du sch\'ema ab\'elien dual de $A$. Alors,
$M=[\widehat A^0\buildrel {\rm can}\over\longrightarrow A]$, o\`u
$\widehat A^0$ est le compl\'et\'e \`a l'origine de $A$ et
$\widehat A^0\buildrel {\rm can}\over\longrightarrow A$ est la
fl\`eche canonique, et $M'=[0\rightarrow A^\natural ]$ sont deux
$1$-motifs g\'en\'eralis\'es en dualit\'e.
\vskip 5mm

{\bf 6. Modules quasi-coh\'erents et coh\'erents sur
les $1$-motifs g\'en\'eralis\'es sur $k$.}
\vskip 5mm

(6.1) Soit $M=[{\cal G}\buildrel u\over\longrightarrow G]$ un
$1$-motif g\'en\'eralis\'e sur $k$ (toujours suppos\'e de
caract\'eristique nulle). Posons
$$
M^0=[{\cal G}^0\buildrel u^0\over\longrightarrow G],
$$
o\`u $u^0$ est la restriction de $u$ \`a ${\cal G}^0$, notons $V'$ le
dual de Cartier de ${\cal G}^0$, de sorte que
$$
H^0(V',{\cal O}_{V'})={\rm Sym}_\bullet^k({\rm Lie}({\cal G}^0)),
$$
et introduisons la ${\cal O}_G$-Alg\`ebre quasi-coh\'erente
(associative, unitaire, mais non commutative en g\'en\'eral)
$$
{\cal D}_{M^0}={\cal O}_G\otimes_kH^0(V',{\cal O}_{V'})
$$
dont la multiplication est d\'etermin\'ee par la r\`egle:
$$
(\varphi_1\otimes\xi_1)\cdot(\varphi_2\otimes
(\xi_2\cdots\xi_n)) =(\varphi_1\varphi_2)\otimes (\xi_1\xi_2
\cdots\xi_n)+ (\varphi_1{\rm Lie}(u^0)(\xi_1)(\varphi_2))\otimes
(\xi_2\cdots\xi_n)
$$
pour toutes sections locales $\varphi_1,\varphi_2$ de ${\cal O}_G$
et tous $\xi_1,\xi_2,\ldots,\xi_n\in {\rm Lie}({\cal G}^0)$.

On rappelle que ${\cal G}$ se d\'ecompose canoniquement en le
produit ${\cal G}^0\times_k{\cal G}^{\rm \acute et}$ et que ${\cal
G}^{\rm \acute et}$ se d\'eploie sur une  extension finie $k'$ de $k$
dans $\overline k$, ce qui permet de voir ${\cal G}^{\rm \acute et}
(k')$ comme un sous-goupe de ${\cal G}(k')$.

\th D\'EFINITION (6.1.1)
\enonce
Un Module sur $M$ est un ${\cal D}_{M^0}$-Module ${\cal M}$ muni
d'une ${\cal G}^{\rm \acute et}$-structure, i.e. d'une famille
d'isomorphismes de $({\cal D}_{k'\otimes_kM^0})$-Modules
$$
a_x:\tau_{u(x)}^*(k'\otimes_k{\cal M})\buildrel\sim \over
\longrightarrow (k'\otimes_k{\cal M})\qquad (\forall x\in {\cal
G}^{\rm \acute et}(k')),
$$
o\`u on a not\'e $\tau_{u(x)}:k'\otimes_kG\rightarrow
k'\otimes_kG$ la translation par $u(x)\in G(k')$, satisfaisant \`a la
condition de co-cycle
$$
a_{x'+x}=a_{x'}\circ\tau_{u(x')}^*(a_x)\qquad (\forall x,x'\in
{\cal G}^{\rm \acute et}(k'))
$$
et \`a la condition de descente galoisienne
$$
({\rm Spec}(\sigma )\times {\rm id}_G)^*(a_x)=a_{\sigma (x)}
\qquad (\forall\sigma\in {\rm Gal}(k'/k),\ \forall x\in {\cal
G}^{\rm \acute et}(k')).
$$

Un tel Module est dit quasi-coh\'erent si le ${\cal O}_G$-Module
sous-jacent est quasi-coh\'erent.

\hfill\hfill\cqfd
\endth

Rappelons que, pour tout $k$-sch\'ema $X$ et tout
${\cal O}_X$-Module quasi-coh\'erent ${\cal M}$, la fl\`eche
d'adjonction
$$
{\cal M} \rightarrow (k'\otimes_k{\cal M})^{{\rm Gal}(k'/k)}
\buildrel{\rm dfn}\over{=\!=}(\alpha_*\alpha^* {\cal M})^{{\rm
Gal}(k'/k)},
$$
o\`u $\alpha :k'\otimes_kX \rightarrow X$ est la projection
canonique, est un isomorphisme.
\vskip 2mm

{\pc CAS PARTICULIERS} (6.1.2)\pointir (i)
Si ${\cal G}=(0)$, un Module (quasi-coh\'erent) sur $M$ n'est rien
d'autre qu'un ${\cal O}_G$-Module (quasi-coh\'erent).

\decale{(ii)} Si ${\cal G}$ est le compl\'et\'e formel $\widehat G^0$
de $G$ \`a l'origine et si $u$ est la fl\`eche canonique $\widehat
G^0\buildrel {\rm can}\over\longrightarrow G$, un Module
(quasi-coh\'erent) sur $M$ n'est rien d'autre qu'un ${\cal
D}_G$-Module (quasi-coh\'erent), o\`u ${\cal D}_G={\cal D}_{M^0}$
est la ${\cal O}_G$-Alg\`ebre des op\'erateurs diff\'erentiels
(relatifs \`a $k$) sur ${\cal O}_G$.

\decale{(iii)} Si ${\cal G}^0=(0)$ et si ${\cal G}^{\rm \acute et}$ est
constant de valeur le groupe ab\'elien libre de rang fini $X$, un
Module (quasi-coh\'erent) sur $M$ est un ${\cal O}_G$-Module
(quasi-coh\'erent) $X$-\'equivariant.
\vskip 2mm

On notera $M_{\rm qcoh}(M)$ la cat\'egorie ab\'elienne des Modules
quasi-coh\'erents sur $M$ et $D_{\rm qcoh}(M)$ sa cat\'egorie
d\'eriv\'ee. On a les sous-cat\'egories habituelles $D_{\rm
qcoh}^+(M)$, $D_{\rm qcoh}^{\rm b}(M)$, ... de $D_{\rm qcoh}(M)$.
\vskip 3mm

Pour tout morphisme $f=(f_{\cal G},f_G): M_1\rightarrow M_2$ de
$1$-motifs g\'en\'eralis\'es sur $k$, nous allons maintenant
d\'efinir des foncteurs
$$
f_*: D_{\rm qcoh}^{\rm b}(M_1)\rightarrow D_{\rm qcoh}^{\rm
b}(M_2)
$$
et
$$
f^!: D_{\rm qcoh}^{\rm b}(M_2)\rightarrow D_{\rm qcoh}^{\rm
b}(M_1).
$$

Le morphisme $f_{\cal G}$ \'etant le produit direct de sa
composante connexe $f_{{\cal G}^0}:{\cal G}_1^0\rightarrow {\cal
G}_2^0$ et de sa composante \'etale $f_{{\cal G}^{\rm \acute
et}}:{\cal G}_1^{\rm \acute et} \rightarrow {\cal G}_2^{\rm \acute
et}$, on a une factorisation canonique
$$
M_1\buildrel f'\over\longrightarrow M_2''\buildrel g'\over
\longrightarrow N_2'\buildrel h'\over\longrightarrow M_2'
\buildrel g\over \longrightarrow N_2\buildrel h\over
\longrightarrow M_2
$$
de $f:M_1\rightarrow M_2$, o\`u
$$\eqalign{
M_2''&=[{\cal G}_1\buildrel f_G\circ u_1\over {\hbox to
8mm{\rightarrowfill}} G_2]\cr f'&= ({\rm id},f_G)\cr
h\circ g\circ h'\circ g'&=(f_{\cal G},{\rm id}),\cr}
$$
$$\eqalign{
M_2'&=[{\cal G}_2^0\times_k{\cal G}_1^{\rm \acute et}\buildrel
u_2^0\cdot (f_G\circ u_1^{\rm \acute et})\over {\hbox to
16mm{\rightarrowfill}} G_2]\cr
h'\circ g'&= (f_{\cal G}^0\times {\rm id},{\rm id})\cr h\circ
g&=({\rm id}\times f_{\cal G}^{\rm\acute et}, {\rm id}),\cr}
$$
$$\eqalign{
N_2'&=[({\cal G}_1^0\times_k{\cal G}_2^0)\times_k{\cal G}_1^{\rm
\acute et}\buildrel 1\cdot u_2^0\cdot (f_G\circ u_1^{\rm \acute
et})\over {\hbox to 18mm{\rightarrowfill}} G_2]\cr
g'&= (({\rm id},f_{\cal G}^0)\times {\rm id},{\rm id})\cr
h'&=({\rm pr}_{{\cal G}_2^0}\times {\rm id},{\rm id}),\cr}
$$
et
$$\eqalign{
N_2&=[{\cal G}_2^0\times_k({\cal G}_1^{\rm \acute et}\times_k
{\cal G}_2^{\rm \acute et})\buildrel u_2^0\cdot 1\cdot u_2^{\rm
\acute et}\over {\hbox to 10mm{\rightarrowfill}} G_2]\cr
g&= ({\rm id}\times ({\rm id},f_{\cal G}^{\rm \acute et}), {\rm
id})\cr
h&=({\rm id}\times {\rm pr}_{{\cal G}_2^{\rm \acute et}},{\rm
id}),\cr}
$$
et, bien entendu, on veut que
$$
f_*=h_*\circ g_*\circ h_*'\circ g_*'\circ f_*'
$$
et que
$$
f^!=f'^!\circ g'^!\circ h'^!\circ g^!\circ h^!.
$$
Il suffit donc de d\'efinir $f_*{\cal M}_1$  pour chaque ${\cal M}_1
\in {\rm ob}\,D_{\rm qcoh}^+(M_1)$ et $f^!{\cal M}_2$  pour chaque
${\cal M}_2\in {\rm ob}\,D_{\rm qcoh}^+(M_2)$ dans les cinq cas
particuliers suivants:

\decale{(1)} $f_{\cal G}={\rm id}$,

\decale{(2)} $f_{{\cal G}^0}$ est un monomorphisme muni  d'une
r\'etraction, $f_{{\cal G}^{\rm \acute et}}={\rm id}$ et $f_G={\rm
id}$,

\decale{(3)} $f_{{\cal G}^0}$ est un \'epimorphisme muni  d'une
section, $f_{{\cal G}^{\rm \acute et}}={\rm id}$ et $f_G={\rm id}$,

\decale{(4)} $f_{{\cal G}^0}={\rm id}$, $f_{{\cal G}^{\rm \acute et}}$
est un monomorphisme muni  d'une r\'etraction  et $f_G={\rm id}$,

\decale{(5)} $f_{{\cal G}^0}={\rm id}$,$f_{{\cal G}^{\rm \acute et}}$
est un \'epimorphisme muni  d'une section  et $f_G={\rm id}$.
\vskip 2mm

Dans le cas particulier (1), l'image directe ordinaire $Rf_{G*} {\cal
M}_1$ et l'image inverse extraordinaire $f_G^!{\cal M}_2$ au sens
des ${\cal O}$-Modules sont trivialement munies de structures de
complexes de Modules quasi-coh\'erents sur $M_2$ et $M_1$
respectivement. Ces complexes sont par d\'efinition
$f_*{\cal M}_1$ et $f^!{\cal M}_2$.

Dans le cas particulier (2), on peut identifier $f_{{\cal G}^0}: {\cal
G}_1^0\rightarrow {\cal G}_2^0$ \`a
$$
({\rm id},0):{\cal G}_1^0\rightarrow {\cal G}_1^0\times_k{\cal H},
$$
o\`u ${\cal H}$ est le noyau de la r\'etraction, et on a alors
$$\eqalign{
{\cal D}_{M_2^0}&={\cal D}_{M_1^0}\otimes_k {\rm
Sym}_\bullet^k({\rm Lie}({\cal H}))\cr
&={\cal O}_G\otimes_k{\rm Sym}_\bullet^k ({\rm Lie}({\cal
G}_1^0))\otimes_k {\rm Sym}_\bullet^k({\rm Lie}({\cal H})).\cr}
$$
On pose d'une part
$$
\omega_{\cal H}=(\Omega_{{\cal H}/k}^{d_{\cal H}})_{(0)}
=\bigwedge_k^{d_{\cal H}}{\rm Hom}_k({\rm Lie}({\cal H}),k)
$$
et
$$
f_*{\cal M}_1={\cal M}_1\otimes_k {\rm Sym}_\bullet^k({\rm
Lie}({\cal H}))\otimes_k\omega_{\cal H}^{\otimes -1}
$$
avec la structure de ${\cal D}_{M_2^0}$-Module \`a gauche d\'efinie
par
$$\displaylines{
\qquad
(1\otimes 1\otimes\zeta_i)\cdot\bigl(m_1\otimes
(\zeta_1^{\alpha_1}\cdots\zeta_d^{\alpha_d})
\otimes (\zeta_1\wedge\cdots\wedge\zeta_d)\bigr)
\hfill\cr\hfill
=-m_1\otimes (\zeta_1^{\alpha_1}\cdots
\zeta_{i-1}^{\alpha_{i-1}}\zeta_i^{\alpha_i+1}
\zeta_{i+1}^{\alpha_{i+1}}\cdots\zeta_d^{\alpha_d})
\otimes(\zeta_1\wedge\cdots\wedge\zeta_d),\qquad}
$$
$$
(\varphi\otimes 1\otimes 1)\cdot\bigl(m_1\otimes 1\otimes
(\zeta_1\wedge\cdots\wedge\zeta_d)\bigr) =(\varphi\cdot
m_1)\otimes 1\otimes (\zeta_1\wedge\cdots\wedge\zeta_d)
$$
et
$$
(1\otimes\xi\otimes 1)\cdot\bigl(m_1\otimes 1\otimes
(\zeta_1\wedge\cdots\wedge\zeta_d)\bigr) = (\xi\cdot
m_1)\otimes 1\otimes (\zeta_1\wedge\cdots\wedge\zeta_d)
$$
pour toute base $(\zeta_1,\ldots ,\zeta_d)$ de ${\rm Lie}({\cal H})$,
tout $d$-uplet $(\alpha_1,\ldots ,\alpha_d)$ d'entiers $\geq 0$,
tout $i=1,\ldots ,d$, toute section locale $m_1$ de ${\cal M}_1$,
toute section locale $\varphi$ de ${\cal O}_G$ et tout $\xi\in {\rm
Lie}({\cal G}_1^0)$, et avec la ${\cal G}_2^{\rm \acute
et}$-structure induite par la ${\cal G}_1^{\rm \acute et}$-structure
de ${\cal M}_1$. On pose d'autre part
$$
f^!{\cal M}_2=\rho ({\cal M}_2)\otimes_k
\omega_{\cal H},
$$
o\`u $\rho ({\cal M}_2)$ est le Module sur $M_1$ d\'eduit du Module
${\cal M}_2$ par oubli de l'action de ${\rm Sym}_\bullet^k({\rm
Lie}({\cal H}))$, avec la ${\cal G}_1^{\rm \acute et}$-structure
induite par la ${\cal G}_2^{\rm \acute et}$-structure de ${\cal
M}_1$.

Dans la cas particulier (3), on peut identifier $f_{{\cal G}^0}: {\cal
G}_1^0\rightarrow {\cal G}_2^0$ \`a
$$
{\rm pr}_{{\cal G}_2^0}:{\cal G}_2^0\times_k{\cal H}\rightarrow
{\cal G}_2^0,
$$
o\`u ${\cal H}$ est le noyau de $f_{{\cal G}^0}$, et on a alors
$$
{\cal D}_{M_1^0}={\cal D}_{M_2^0}\otimes_k {\rm
Sym}_\bullet^k({\rm Lie}({\cal H}))
$$
(produit tensoriel de $k$-Alg\`ebres). On fait agir ${\rm
Sym}_\bullet^k({\rm Lie}({\cal H}))$ sur  le $k$-espace vectoriel
$\omega_{\cal H}$ \`a travers l'augmentation canonique ${\rm
Sym}_\bullet^k ({\rm Lie}({\cal H}))\twoheadrightarrow k$. On pose
d'une part
$$
f_*{\cal M}_1=\omega_{\cal H}\,\Lotimes_{{\rm Sym}_\bullet^k
({\rm Lie}({\cal H}))}\,{\cal M}_1,
$$
vu comme ${\cal D}_{M_2^0}$-Module, avec la ${\cal G}_2^{\rm
\acute et}$-structure induite par la ${\cal G}_1^{\rm \acute
et}$-structure de ${\cal M}_1$. On  pose d'autre part
$$
f^!{\cal M}_2={\cal M}_2\otimes_k\omega_{\cal H}^{\otimes -1}
$$
avec l'action de ${\cal D}_{M_1^0}$ produit tensoriel de l'action
donn\'ee de ${\cal D}_{M_2^0}$ sur ${\cal M}_2$ et de l'action de
${\rm Sym}_\bullet^k({\rm Lie}({\cal H}))$ sur $\omega_{\cal
H}^{\otimes -1}$ \`a travers l'augmentation canonique et avec la
${\cal G}_1^{\rm \acute et}$-structure induite par la ${\cal
G}_2^{\rm \acute et}$-structure de ${\cal M}_1$.

Dans le cas particulier (4), on peut identifier $f_{{\cal G}^{\rm
\acute et}}:{\cal G}_1^{\rm \acute et}\rightarrow {\cal G}_2^{\rm
\acute et}$ \`a
$$
({\rm id},0):{\cal G}_1^{\rm \acute et}\rightarrow {\cal G}_1^{\rm
\acute et}\times_k{\cal H},
$$
o\`u ${\cal H}$ est le noyau de la r\'etraction. On pose d'une part
$$
\omega_{\cal H}\buildrel{\rm dfn}\over{=\!=}
\Bigl(\bigwedge^{r_{\cal H}} {\rm Hom}_{k'{\rm
-gr}}(k'\otimes_k{\cal H}, {\bb G}_{{\rm a},k'})\Bigl)^{{\rm
Gal}(k'/k)}
$$
et
$$
f_*{\cal M}_1=\Bigl(\bigoplus_{y\in {\cal H}(k')}\tau_{u_2^{\rm
\acute et}(y)}^*(k'\otimes_k {\cal M}_1)\Bigl)^{{\rm
Gal}(k'/k)}\otimes_k\omega_{\cal H}^{\otimes -1},
$$
muni de la structure de ${\cal D}_{M_2^0}$-Module induite par la
structure de ${\cal D}_{M_1^0}$-Module de ${\cal M}_1$ et muni de
la ${\cal G}_2^{\rm \acute et}$-structure naturelle: pour tout
$x_2=(z,x_1)\in {\cal H}(k')\times {\cal G}_1^{\rm \acute
et}(k')={\cal G}_2^{\rm \acute et}(k')$, on a l'isomorphisme
$$\displaylines{
\qquad
\tau_{u_2^{\rm \acute et}(x_2)}^*\Bigl(\bigoplus_y\tau_{u_2^{\rm
\acute et}(y)}^*(k'\otimes_k  {\cal M}_1)\Bigr)\cong \bigoplus_y
\tau_{u_2^{\rm \acute et}(y)}^*\tau_{u_1^{\rm \acute
et}(x_1)}^*(k'\otimes_k {\cal M}_1)
\hfill\cr\hfill
\buildrel\bigoplus_y\tau_{u_2^{\rm \acute et}(y)}^*(a_{x_1})\over
{\hbox to 22mm{\rightarrowfill}}\bigoplus_y\tau_{u_2^{\rm \acute
et}(y)}^*(k'\otimes_k {\cal M}_1)\qquad}
$$
o\`u $y$ parcourt ${\cal H}(k')$. On pose d'autre part
$$
f^!{\cal M}_2=\rho ({\cal M}_2)\otimes_k\omega_{\cal H},
$$
o\`u $\rho ({\cal M}_2)$ est le Module sur $M_1$ d\'eduit du Module
${\cal M}_2$ sur $M_2$ par oubli de la ${\cal H}$-structure.

Dans le cas particulier (5), on peut identifier $f_{{\cal G}^{\rm
\acute et}}:{\cal G}_1^{\rm \acute et}\rightarrow {\cal G}_2^{\rm
\acute et}$ \`a
$$
{\rm pr}_{{\cal G}_2^{\rm \acute et}}:{\cal G}_2^{\rm \acute
et}\times_k{\cal H}\rightarrow {\cal G}_2^{\rm \acute et},
$$
o\`u ${\cal H}$ est le noyau de $f_{{\cal G}^{\rm \acute et}}$. On fait
agir $k[{\cal H}(k')]$ sur le $k$-espace vectoriel $\omega_{\cal H}$
\`a travers l'augmentation canonique $k[{\cal H}(k')]
\twoheadrightarrow k$. On pose d'une part
$$
f_*{\cal M}_1=\Bigl(\omega_{\cal H}\,\Lotimes_{k[{\cal H}(k')]}\,
(k'\otimes_k{\cal M}_1)\Bigr)^{{\rm Gal}(k'/k)}
$$
avec la structure de ${\cal D}_{M_2^0}$-Module induite par la
structure de ${\cal D}_{M_1^0}$-Module de ${\cal M}_1$ et la ${\cal
G}_2^{\rm \acute et}$-structure d\'eduite de la ${\cal G}_1^{\rm
\acute et}$-structure de ${\cal M}_1$ par oubli de l'action de ${\cal
H}$. On pose d'autre part
$$
f^!{\cal M}_2={\cal M}_2\otimes_k\omega_{\cal H}^{\otimes -1}
$$
avec la structure de ${\cal D}_{M_1^0}$-Module induite par la
structure de ${\cal D}_{M_2^0}$-Module de ${\cal M}_2$ et avec la
${\cal G}_1^{\rm \acute et}$-structure produit tensoriel de la
${\cal G}_2^{\rm \acute et}$-structure donn\'ee sur ${\cal M}_2$ et
de la ${\cal H}$-structure triviale sur $\omega_{\cal H}^{\otimes
-1}$.

On laissera au lecteur le soin de v\'erifier que, pour tout
morphisme compos\'e
$$
M_1\buildrel f\over\longrightarrow M_2\buildrel g\over
\longrightarrow M_3
$$
dans la cat\'egorie des $1$-motifs g\'en\'eralis\'es, on a des
isomorphismes de transitivit\'e
$$
(g\circ f)_*\cong g_*\circ f_*
$$
et
$$
(g\circ f)^!\cong f^!\circ g^!.
$$
\vskip 2mm

{\pc EXEMPLES} (6.1.3)\pointir (i) Pour $f=(0,f_G):[0\rightarrow
G_1] \rightarrow [0 \rightarrow G_2]$, $f_*$ et $f^!$ sont
respectivement l'image directe ordinaire $Rf_{G*}$ et l'image
inverse extraordinaire $f_G^!$ pour les ${\cal O}$-Modules.

\decale{\rm (ii)} Pour $f=(\widehat f_G,f_G):[\widehat G_1^0
\buildrel{\rm can}\over\longrightarrow G_1] \rightarrow
[\widehat G_2^0\buildrel{\rm can}\over\longrightarrow G_2]$,
$f_*$ et $f^!$ sont respectivement l'image directe ordinaire
$f_{G+}$ et l'image inverse extraordinaire $f_G^!$ pour les
${\cal D}$-Modules \`a gauche (cf. [Bo] \S 4 et \S 5 ou  les rappels
(3.3.1)).

\decale{\rm (iii)} Pour $f=(0,{\rm id}_G):[0\rightarrow  G]
\rightarrow [\widehat G^0\buildrel{\rm can}\over\longrightarrow
G]$, on a
$$
f_*(\cdot )=(\cdot )\otimes_{{\cal O}_G}{\cal D}_G^\Omega ,
$$
o\`u
$$
{\cal D}_G^\Omega ={\cal D}_G\otimes_{{\cal O}_G}
(\Omega_{G/k}^{d_G})^{\otimes -1}
$$
est muni de sa structure de $({\cal O}_G,{\cal D}_G)$-bi-Module
gauche-gauche naturelle (cf. [Bo] \S 3),
et on a
$$
f^!(\cdot )=\rho (\cdot )\otimes_{{\cal O}_G}\Omega_{G/k}^{d_G},
$$
o\`u $\rho :D_{\rm qcoh}^{\rm b}({\cal D}_G)\rightarrow D_{\rm
qcoh}^{\rm b}({\cal O}_G)$ est le foncteur d'oubli de la structure de
${\cal D}_G$-Module.

\decale{\rm (iv)} Pour $M=[X\rightarrow G]$, o\`u $X$ est un
${\bb Z}$-module libre de rang fini $r$, avec action triviale de
${\rm Gal}(k'/k)$, et pour $f:M\rightarrow [0\rightarrow 0]={\rm
Spec}(k)$ la projection canonique, $f_*$ est le foncteur de
cohomologie \'equivariante
$$
R\Gamma^X(G,\,\cdot\,)[r]=R\Gamma (X,R\Gamma (G,\,\cdot\,))[r]
$$
puisque l'on a la r\'esolution canonique
$$\displaylines{
\qquad k\hookrightarrow k[X] \rightarrow k[X]\otimes_{\bb Z}
X^\vee\rightarrow k[X]\otimes_{\bb Z}\bigwedge^2X^\vee
\rightarrow \cdots
\hfill\cr\hfill
\cdots\rightarrow
k[X]\otimes_{\bb Z}\bigwedge^{r-1}X^\vee\rightarrow k[X]
\otimes_{\bb Z}\bigwedge^rX^\vee\twoheadrightarrow k
\otimes_{\bb Z}\bigwedge^rX^\vee=\omega_X,\qquad}
$$
o\`u la diff\'erentielle
$$
k[X]\otimes_{\bb Z}\bigwedge^iX^\vee\rightarrow k[X] \otimes_{\bb
Z}\bigwedge^{i+1}X^\vee
$$
est d\'efinie par
$$
[x]\otimes (x_1^\vee\wedge\cdots x_i^\vee )\mapsto
\sum_{j=1}^r\langle e_j^\vee ,x\rangle [x-e_j]\otimes
(e_j^\vee\wedge x_1^\vee\wedge\cdots x_i^\vee )
$$
si $(e_j)_{j=1,\ldots ,r}$ est une base de $X$ et si
$\sum_{j=1}^re_j\otimes e_j^\vee\in X\otimes X^\vee$ est
l'\'el\'ement qui est envoy\'e sur l'identit\'e par l'isomorphisme
canonique de $X\otimes_{\bb Z}X^\vee$ sur ${\rm End}_{\bb Z}(X)$.

\decale{\rm (v)} Pour $M=[{\cal G}\rightarrow G]$ arbitraire et pour
$f:M\rightarrow [0\rightarrow 0]={\rm Spec}(k)$ la projection
canonique, on a
$$
f^!k={\cal O}_G\otimes_k\omega_M[d_G]
$$
avec sa structure naturelle de Module sur $M$, o\`u on a pos\'e
$$
\omega_M=\omega_G\otimes\omega_{\cal G}^{\otimes -1}
$$
avec
$$
\omega_G=(\Omega_{G/k}^{d_G})_{(1)}
$$
et
$$
\omega_{\cal G}=\omega_{{\cal G}^0}\otimes_k
\omega_{{\cal G}^{\rm \acute et}} =(\Omega_{{\cal G}^0/k}^{d_{\cal
G}})_{(1)}\otimes_k\Bigl(\bigwedge^{r_{\cal H}} {\rm Hom}_{k'{\rm
-gr}}(k'\otimes_k {\cal H},{\bb G}_{{\rm a},k'})\Bigl)^{{\rm
Gal}(k'/k)}.
$$

\hfill\hfill\cqfd
\vskip 2mm

Si $M'=[{\cal G}'\buildrel u'\over\longrightarrow G']$ et $M''=[{\cal
G}''\buildrel u''\over\longrightarrow G'']$ sont deux $1$-motifs
g\'en\'eralis\'es, on dispose aussi d'un produit tensoriel externe
$$
(\cdot )\boxtimes (\cdot ):D_{\rm qcoh}^{\rm b}(M')\times D_{\rm
qcoh}^{\rm b}(M'')\rightarrow  D_{\rm qcoh}^{\rm b}(M'\times M''),
$$
o\`u $M'\times M''=[{\cal G}'\times_k {\cal G}''\buildrel u'\times
u''\over{\hbox to 12mm{\rightarrowfill}}  G'\times_kG'']$. Par
suite, si $M=[{\cal G}\buildrel u\over\longrightarrow G]$ est un
$1$-motif g\'en\'eralis\'e, on peut d\'efinir un produit tensoriel
interne
$$
(\cdot )\otimes^! (\cdot ):D_{\rm qcoh}^{\rm b}(M)\times D_{\rm
qcoh}^{\rm b}(M)\rightarrow  D_{\rm qcoh}^{\rm b}(M)
$$
par
$$
(\cdot )\otimes^!(\cdot )=\Delta_M^!((\cdot )\boxtimes (\cdot )),
$$
o\`u $\Delta_M=(\Delta_{\cal G},\Delta_G):M \rightarrow M\times
M$ est le morphisme diagonal.

On notera simplement
$$
(\cdot )\otimes_k(\cdot ):D_{\rm qcoh}^{\rm b}(M)\times D^{\rm
b}(k) \rightarrow D_{\rm qcoh}^{\rm b}(M)
$$
le foncteur $(\cdot )\otimes^!\pi^!(\cdot )$ pour $\pi :M\rightarrow
{\rm Spec}(k)$ le morphisme canonique (ce foncteur n'est autre que
le foncteur $(\cdot )\boxtimes (\cdot )$ pour $M_1=M$ et
$M_2=[0\rightarrow 0]$ puisque $M\times_k[0\rightarrow 0]$
s'identifie canoniquement \`a $M$).

On laisse au lecteur le soin de d\'eduire la proposition suivante des
r\'esultats analogues pour les ${\cal O}$-Modules (cf. [Ha]).

\th PROPOSITION (6.1.4)
\enonce
{\rm (i)} {\rm (Changement de base)} Pour tout carr\'e cart\'esien
$$
\matrix{M_1'&\maprightover{\alpha_1} &M_1\cr
\noalign{\smallskip}
\mapdownleft{f'}&&\mapdownright{f}\cr
\noalign{\smallskip}
M_2'&\maprightunder{\alpha_2} &M_2\cr}
$$
dans ${\bf Mot}/k$ tel que $f_G$ et $f_{\cal G}$ soient des
\'epimorphismes, les foncteurs $\alpha_2^!\circ f_*$ et $f_*'\circ
\alpha_1^!$ de $D_{\rm qcoh}^{\rm b}(M_1)$ dans $D_{\rm qcoh}^{\rm
b}(M_2')$ sont naturellement isomorphes.

\decale{\rm (ii)} {\rm (Formule des projections)} Pour tout
morphisme $f:M_1 \rightarrow M_2$ dans ${\bf Mot}/k$, les
foncteurs $f_*((\cdot )\otimes^!f^!(\cdot ))$ et $f_*(\cdot )
\otimes^!(\cdot )$ de $D_{\rm qcoh}^{\rm b}(M_1)\times D_{\rm
qcoh}^{\rm b}(M_2)$ dans $D_{\rm qcoh}^{\rm b} (M_2)$ sont
naturellement isomorphes.

\decale{\rm (iii)} {\rm (Compatibilit\'e au produit externe)}
Pour tous morphismes $f':M_1' \rightarrow M_2'$ et $f'':M_1''
\rightarrow M_2''$ dans ${\bf Mot}/k$, les foncteurs $f_*'(\cdot
)\boxtimes f_*''(\cdot )$ et $(f'\times f'')_*((\cdot )\boxtimes
(\cdot ))$ de $D_{\rm qcoh}^{\rm b}(M_1')\times D_{\rm qcoh}^{\rm
b}(M_1'')$ dans $D_{\rm qcoh}^{\rm b}(M_2'\times M_2'')$ d'une part
et les foncteurs $f'^!(\cdot )\boxtimes f''^!(\cdot )$ et $(f'\times
f'')^!((\cdot )\boxtimes (\cdot ))$ de $D_{\rm qcoh}^{\rm
b}(M_2')\times D_{\rm qcoh}^{\rm b} (M_2'')$ dans $D_{\rm
qcoh}^{\rm b}(M_1'\times M_1'')$ d'autre part sont naturellement
isomorphes.

\hfill\hfill\cqfd
\endth

{\pc REMARQUE} 6.1.5\pointir L'hypoth\`ese de l'\'enonc\'e de
changement de base n'est certainement pas l'hypoth\`ese optimale.
On laisse au lecteur le soin d'introduire une notion de
``Tor-ind\'ependance'' pour deux morphismes $f: M_1\rightarrow
M_2$ et $\alpha_2: M_2'\rightarrow M_2$ qui g\'en\'eralise la
notion de Tor-ind\'ependance pour les morphismes de sch\'emas et
qui soit une hypoth\`ese suffisante pour le changement de base (cf.
[SGA6] Ch. IV, Prop. 3.1.0). Bien entendu, cette hypoth\`ese de
Tor-ind\'ependance doit \^etre automatiquement v\'erifi\'ee dans le
cas o\`u $f_G$ et $f_{\cal G}$ sont tous les deux des
\'epimorphismes.

\hfill\hfill\cqfd
\vskip 5mm

(6.2) Soit $M=[{\cal G}\buildrel u\over\longrightarrow G]$ un
$1$-motif g\'en\'eralis\'e sur $k$. Pour tout ${\cal O}_G$-Module
quasi-coh\'erent ${\cal F}$, on notera
$$
{\rm Ind}_G^M({\cal F})
$$
le Module quasi-coh\'erent sur le $M$ {\it induit par} ${\cal F}$. Par
d\'efinition, ${\rm Ind}_G^M({\cal F})$ est le ${\cal
D}_{M^0}$-Module
$$
\Bigl(\bigoplus_{x\in {\cal G}^{\rm \acute et}(k')}
\tau_{u(x)}^*(k'\otimes_k{\cal D}_{M^0}
\otimes_{{\cal O}_G}{\cal F})\Bigr)^{{\rm Gal}(k'/k)}
$$
muni de la ${\cal G}^{\rm \acute et}$-structure \'evidente. On a
encore
$$
{\rm Ind}_G^M({\cal F})={\rm Ind}_{M^0}^M ({\rm Ind}_G^{M^0}({\cal
F})),
$$
o\`u
$$
{\rm Ind}_G^{M^0}({\cal F})={\cal D}_{M^0}\otimes_{{\cal O}_G}{\cal
F}
$$
et
$$
{\rm Ind}_{M^0}^M({\cal M}^0)=\Bigl(\bigoplus_{x\in {\cal G}^{\rm
\acute et}(k')}\tau_{u(x)}^*(k'\otimes_k{\cal M}^0)\Bigr)^{{\rm
Gal}(k'/k)}
$$
pour tout ${\cal D}_{M^0}$-Module quasi-coh\'erent ${\cal M}^0$.

Pour tout Module quasi-coh\'erent ${\cal M}$ sur $M$, de ${\cal
O}_G$-Module sous-jacent $\rho ({\cal M})$, on a une fl\`eche
d'adjonction
$$
{\rm Ind}_G^M(\rho ({\cal M}))\rightarrow  (k'\otimes_k{\cal
M})^{{\rm Gal}(k'/k)}= {\cal M}
$$
qui est clairement un \'epimorphisme de Modules sur $M$.

\th D\'EFINITION (6.2.1)
\enonce
Un Module coh\'erent sur $M$ est un Module quasi-coh\'erent ${\cal
M}$ sur $M$ pour lequel on puisse trouver au moins un ${\cal
O}_G$-Module coh\'erent ${\cal F}$ et un \'epimorphisme
$$
{\rm Ind}_G^M({\cal F})\twoheadrightarrow {\cal M}
$$
de Modules quasi-coh\'erents sur $M$.
\endth

\th LEMME (6.2.2)
\enonce
{\rm (i)} Pour tout ouvert affine $U$ de $G$, l'anneau
${\cal D}_{M^0}(U)$ est noeth\'erien.

\decale{\rm (ii)} L'Anneau ${\cal D}_{M^0}$ est coh\'erent et un
Module coh\'erent sur $M^0$ n'est rien d'autre qu'un
${\cal D}_{M^0}$-Module \`a gauche coh\'erent.

\endth

\rem Preuve
\endrem
On peut munir l'Anneau ${\cal D}_{M^0}$ de la filtration croissante
et exhaustive
$$
{\cal D}_{M^0,i}={\cal O}_G\otimes_k\Bigl(\bigoplus_{j=1}^i{\rm
Sym}_j^k({\rm Lie}({\cal G}^0))\Bigr)
$$
dont le gradu\'e associ\'e est l'Anneau (maintenant commutatif)
$$
{\rm gr}_\bullet{\cal D}_{M^0}={\rm Sym}_\bullet^k ({\cal
O}_G\otimes_k{\rm Lie}({\cal G}^0)).
$$
On conclut alors comme dans [Bo] Ch. II, \S 3.

\hfill\hfill\cqfd
\vskip 2mm

On notera $M_{\rm coh}(M)$ la sous-cat\'egorie (strictement)
pleine de $M_{\rm qcoh}(M)$ dont les objets sont les Modules
coh\'erents.

\th PROPOSITION (6.2.3)
\enonce
La cat\'egorie ab\'elienne $M_{\rm qcoh}(M)$ est localement
noeth\'erienne et les objets de $M_{\rm coh}(M)$ sont  exactement
les objets noeth\'eriens de $M_{\rm qcoh}(M)$. En particulier,
$M_{\rm coh}(M)$ est une sous-cat\'egorie ab\'elienne de $M_{\rm
qcoh}(M)$ stable par sous-objets, quotients et extensions.
\endth

\rem Preuve {\rm (Thomasson)}
\endrem
Il suffit de montrer que tout Module quasi-coh\'erent sur $M$ est
limite inductive de ses sous-objets coh\'erents et que les Modules
coh\'erents sur $M$ sont exactement les objets noeth\'eriens de
$M_{\rm qcoh}(M)$.

Mais tout Module quasi-coh\'erent ${\cal M}$ sur $M$ est un
quotient de ${\rm Ind}_G^M(\rho ({\cal M}))$ et $\rho ({\cal M})$ est
la limite inductive de ses sous-${\cal O}_G$-Modules coh\'erents.
Par suite ${\cal M}$ est bien la limite inductive de ses
sous-Modules coh\'erents. En particulier, tout objet noeth\'erien de
$M_{\rm qcoh}(M)$ est n\'ecessairement coh\'erent.

Inversement, montrons que tout Module coh\'erent ${\cal M}$  sur
$M$ est noeth\'erien. L'action de ${\rm Gal}(k'/k)$ sur ${\cal
G}^{\rm \acute et}(k')$ se factorise \`a travers un quotient fini
${\rm Gal}(k'/k)$ et donc, quitte \`a \'etendre les scalaires de $k$
\`a $k'$, on peut supposer que ${\cal G}^{\rm \acute et}$ est
constant de valeur un groupe ab\'elien $X$ libre de rang fini $r$. On
proc\`ede alors par r\'ecurrence sur $r$. Si $r=0$, l'assertion
r\'esulte imm\'ediatement du lemme (6.2.2). Si $r>0$, on
d\'ecompose $X$ en ${\bb Z}\times X'$ et on suppose par
r\'ecurrence que tout Module coh\'erent ${\cal M}'$ sur
$M'=[{\cal G}^0\times X'\buildrel u'\over\longrightarrow G]$ est
noeth\'erien ($u'$ est la restriction de $u$ \`a ${\cal G}^0\times
X'\subset {\cal G}^0\times X$). Puisque ${\cal M}$ est coh\'erent, on
peut trouver un Module coh\'erent ${\cal M}'$ sur $M'$ et un
\'epimorphisme
$$
\pi :{\rm Ind}_{M'}^M({\cal M}')\buildrel{\rm dfn}\over{=\!=}
\bigoplus_{i\in {\bb Z}}\tau_{u(i,0)}^*{\cal M}'\twoheadrightarrow
{\cal M}.
$$

Dans un premier temps, montrons que, pour tout sous-Module
quasi-coh\'erent ${\cal N}$ de ${\cal M}$, il existe un entier $i_0$
tel que
$$
{\cal N}=\sum_{j\in {\bb Z}}\tau_{u(j,0)}^*\bigl({\cal N}\cap \pi
({\cal M}'\oplus\tau_{u(1,0)}^*{\cal M}'\oplus\cdots
\oplus\tau_{u(i_0,0)}^*{\cal M}')\bigr).
$$
Pour cela, posons
$$
{\cal N}_i'=\tau_{u(-i,0)}^*\biggl({\pi^{-1}({\cal N})\cap
\bigoplus_{j=0}^i\tau_{u(j,0)}^*{\cal M}'\over\pi^{-1} ({\cal
N})\cap\bigoplus_{j=0}^{i-1}\tau_{u(j,0)}^*{\cal M}'}\biggr)\subset
\tau_{u(0,0)}^*{\cal M}'={\cal M}'
$$
pour chaque entier $i\geq 0$ et remarquons que
$$
{\cal N}_0'\subset{\cal N}_1'\subset {\cal N}_2'\subset\cdots
\subset {\cal M}'.
$$
Alors, comme ${\cal M}'$ est noeth\'erien, il existe un entier
$i_0\geq 0$ tel que ${\cal N}_i'={\cal N}_{i_0}'$ pour tout $i\geq
i_0$. L'entier $i_0$ r\'epond \`a la question. En effet, si $n$ est une
section de ${\cal N}$ image par $\pi$ de
$$
\sum_{i=a}^b\tau_{u(i,0)}^*(m_i')=\tau_{u(a,0)}^*
\Bigl(\sum_{i=0}^{b-a}\tau_{u(i,0)}^*(m_{i+a}')\Bigr),
$$
soit $b-a$ est inf\'erieur ou \'egal \`a $i_0$ et $n$ est une section
de  $\tau_{u(a,0)}^*\bigl({\cal N}\cap \pi ({\cal M}'\oplus
\tau_{u(1,0)}^*{\cal M}'\oplus\cdots\oplus\tau_{u(i_0,0)}^* {\cal
M}') \bigr)$, soit $b-a$ est strictement sup\'erieur \`a $i_0$ et,
comme $m_b'$ est une section de ${\cal N}_{b-a}'$, on peut modifier
$n$ par une section de $\tau_{u(b-i_0,0)}^*\bigl({\cal N}\cap \pi
({\cal M}'\oplus\tau_{u(1,0)}^*{\cal M}'\oplus\cdots\oplus
\tau_{u(i_0,0)}^* {\cal M}')\bigr)$ de mani\`ere \`a diminuer la
diff\'erence $b-a$ d'une unit\'e sans changer $a$.

Maintenant, soit
$$
{\cal N}_0\subset {\cal N}_1\subset {\cal N}_2\subset\cdots
\subset {\cal M}
$$
une tour de sous-Modules quasi-coh\'erents. Pour chaque entier
$j\geq 0$, on consid\`ere les sous-Modules ${\cal N}_i'$ de ${\cal
M}'$ et l'entier $i_0\geq 0$ correspondants \`a ${\cal N}_j$ que l'on
a d\'efinis ci-dessus. On les notera ${\cal N}_{ji}'$ et $i_j$
respectivement. Clairement,  on a
$$
{\cal N}_{0i}'\subset {\cal N}_{1i}'\subset {\cal N}_{2i}'\subset
\cdots\subset {\cal M}'
$$
pour chaque $i\geq 0$ et on peut supposer que $i_0\leq i_1\leq
i_2\leq \cdots$, de sorte que
$$
{\cal N}_{0,i_0}'\subset {\cal N}_{1,i_1}'\subset {\cal N}_{2,i_2}'
\subset\cdots\subset {\cal M}'.
$$
Utilisant de nouveau le fait que ${\cal M}'$ est noeth\'erien, on voit
qu'il existe un entier $j_0$ tel que ${\cal N}_{j,i_j}' ={\cal
N}_{j_0,i_{j_0}}'$ pour tout $j\geq j_0$. Il s'en suit facilement que
l'on peut supposer que $i_j=i_{j_0}$ pour tout $j\geq j_0$.

Finallement, le Module ${\cal M}'\oplus\tau_{u(1,0)}^* {\cal
M}'\oplus\cdots\oplus\tau_{u(i_{j_0},0)}^*{\cal M}'$ sur $M'$ est
lui aussi noeth\'erien et la suite croissante de ses sous-Modules
$$
\pi^{-1}({\cal N}_j)\cap ({\cal M}'\oplus \tau_{u(1,0)}^*{\cal M}'
\oplus\cdots\oplus\tau_{u(i_{j_0},0)}^* {\cal M}')
$$
est donc stationnaire, de sorte que la suite des ${\cal N}_j$
elle-m\^eme est stationnaire, ce que l'on voulait d\'emontrer.

\hfill\hfill\cqfd

On dispose d'une dualit\'e parfaite
$$
D: D_{\rm coh}^{\rm b}(M)^{\rm opp}\rightarrow D_{\rm coh}^{\rm
b}(M)
$$
d\'efinie comme suit. On pose
$$
{\cal D}_{M^0}^\Omega ={\cal D}_{M^0}\otimes_k
\omega_{{\cal G}^0}^{\otimes -1}.
$$
On munit ${\cal D}_{M^0}^\Omega$ d'une structure de bi-${\cal
D}_{M^0}$-Module gauche-gauche de la fa\c con suivante. D'une part,
on a
$$
{\cal D}_{M^0}^\Omega =f_*{\cal O}_G
$$
o\`u $f=(0,{\rm id}):[0\rightarrow G] \rightarrow [{\cal
G}^0\buildrel u^0\over\longrightarrow G]=M^0$ et on a vu comment
munir $f_*{\cal O}_G$ d'une structure de ${\cal D}_{M^0}$-Module
\`a gauche, c'est la structure (1). D'autre part, ${\cal D}_{M^0}$ agit
par multiplication \`a gauche sur ${\cal D}_{M^0}$, ce qui induit une
structure de ${\cal D}_{M^0}$-Module \`a gauche sur ${\cal D}_{M^0}
\otimes_k\omega_{{\cal G}^0}^{\otimes -1}$, c'est la structure (2).
Cette structure de bi-${\cal D}_{M^0}$-Module gauche-gauche sur
${\cal D}_{M^0}^\Omega $ en induit bien entendu une sur
$$
{\cal D}_M^\Omega\buildrel{\rm dfn}\over{=\!=}
\Bigl(\bigoplus_{x\in {\cal G}^{\rm \acute et}(k')}\tau_{u^{\rm
\acute et}(x)}^*(k'\otimes_k {\cal D}_{M^0}^\Omega )\Bigr)^{{\rm
Gal}(k'/k)}\otimes_k\omega_{{\cal G}^{\rm \acute et}}^{\otimes -1}.
$$
On munit ${\cal D}_M^\Omega$ d'une structure de bi-Module sur $M$
de la fa\c con suivante. D'une part, on a
$$
{\cal D}_M^\Omega =g_*{\cal D}_{M^0}
$$
o\`u $g=(({\rm id},0),{\rm id}):[{\cal G}^0\rightarrow G]
\rightarrow [{\cal G}^0\times_k{\cal G}^{\rm \acute et}
\buildrel u\over\longrightarrow G]=M$ et on a vu comment munir le
${\cal D}_{M^0}$-Module \`a gauche $g_*{\cal D}_{M^0}^\Omega$
(structure (1)) d'une  ${\cal G}^{\rm \acute et}$-structure et donc
d'une structure de Module sur $M$, c'est la structure (1). D'autre
part, le ${\cal D}_{M^0}$-Module \`a gauche ${\cal D}_{M^0}^\Omega$
(structure (2))  admet lui-m\^eme une ${\cal G}^{\rm \acute
et}$-structure naturelle, d'o\`u une autre structure de Module sur
$M$, c'est la structure (2).

Alors, pour tout Module coh\'erent ${\cal M}$ sur $M$, $D({\cal M})$
est le complexe de Modules sur $M$,
$$
R\underline {\rm Hom}_{D_{\rm qcoh}(M)}({\cal M},{\cal
D}_M^\Omega)\otimes_k\omega_M[d_G]
$$
(le $R\underline {\rm Hom}$ est form\'e en utilisant la structure
(1) de Module sur $M$ de ${\cal D}_M^\Omega$ et la structure de
Module sur $M$ de ce $R\underline {\rm Hom}$ provient elle de la
structure (2) de Module sur $M$ de ${\cal D}_M^\Omega$. On
v\'erifie que ce complexe est \`a cohomologie coh\'erente en
r\'esolvant ${\cal M}$ par des Modules induits: on a
$$
D({\rm Ind}_G^M({\cal F})) ={\rm Ind}_G^M(R\underline{\rm
Hom}_{{\cal O}_G}({\cal F},{\cal O}_G))\otimes_k\omega_M
\otimes_k\omega_{\cal G}^{\otimes -1}[d_G].
$$

Bien entendu, dans le cas $M=[0\rightarrow G]$, $D$ n'est autre que
le foncteur de dualit\'e pour les ${\cal O}_G$-Modules, $R\underline
{\rm Hom}_{{\cal O}_G} (\,\cdot\,,{\cal O}_G))\otimes_k\omega_G
[d_G]$, et dans le cas $M=[\widehat G^0\rightarrow G]$, $D$ n'est
autre que le foncteur de dualit\'e pour les ${\cal D}_G$-Modules \`a
gauche.

\th PROPOSITION (6.2.4)
\enonce
{\rm (i)} Si la composante $f_G$ de $f$ est propre, on a
$$
f_*({\rm ob}\,D_{\rm coh}^{\rm b}(M_1))\subset {\rm ob}\, D_{\rm
coh}^{\rm b}(M_2)
$$
et
$$
D\circ f_*\cong f_*\circ D.
$$
\decale{\rm (ii)} Si $f_{\cal G}^0$ est un \'epimorphisme et si le
conoyau de $f_{\cal G}^{\rm \acute et}$ est fini sur $k$, on a
$$
f^!({\rm ob}\,D_{\rm coh}^{\rm b}(M_2))\subset {\rm ob}\,D_{\rm
coh}^{\rm b}(M_1)
$$
et
$$
D(f^!(\cdot ))\cong f^!(D(\cdot ))\otimes_k\omega_{M_2}
\otimes_k\omega_{M_1}^{\otimes -1}[d_{G_2}-d_{G_1}+
d_{{\cal G}_2}-d_{{\cal G}_1}+r_{{\cal G}_2}-r_{{\cal G}_1}].
$$

\decale{\rm (iii)} On a
$$
{\rm ob}\,D_{\rm coh}^{\rm b}(M_1)\boxtimes {\rm ob}\,D_{\rm
coh}^{\rm b}(M_2)\subset {\rm ob}\,D_{\rm coh}^{\rm b}(M_1\times
M_2)
$$
et
$$
D((\cdot )\boxtimes (\cdot ))=D(\cdot )\boxtimes D(\cdot ).
$$
\endth

\rem Preuve
\endrem
L'assertion (iii) est triviale.

Pour les assertions (i) et (ii) nous nous contenterons de v\'erifier
quelques cas particuliers.

Tout d'abord, si $f=(0,f_G):M_1=[0 \rightarrow G_1]\rightarrow [0
\rightarrow G_2]=M_2$, la dualit\'e est celle des ${\cal
O}$-Modules et les \'enonc\'es sont ceux de Grothendieck (cf. [Ha]
Ch. III, \S 5, Thm. 11.1, pour $f_*$ et [Ha] Ch. III, \S 8, Prop. 8.8, \S
7, Cor. 7.3, et \S 2, Def.).

Supposons dor\'enavant que $G_1=G_2=G$ et que $f_G={\rm id}$.

Alors, l'hypoth\`ese de la partie (i) est automatiquement
v\'erifi\'ee. Il suffit donc de montrer les assertions de cette partie
dans chacun des cas particuliers (2) \`a (5) consid\'er\'e dans la
section (6.1). Faisons le par exemple dans le cas (2), en supposant
plus ${\cal G}_1^{\rm \acute et} ={\cal G}_2^{\rm \acute et}=(0)$,
et dans le cas (5), en supposant de plus que ${\cal H}$ est
d\'eploy\'e sur $k$. Dans le premier de ces cas, on a
$$
\displaylines{
\qquad R\underline {\rm Hom}_{D_{\rm qcoh}(M_2)}(f_*{\cal
M}_1,{\cal D}_{M_2}^\Omega)
\hfill\cr\hfill
\eqalign{&\cong R\underline {\rm Hom}_{{\cal D}_{M_1^0}
\otimes_k{\rm Sym}_\bullet^k({\rm Lie}({\cal H}))}\bigl({\cal
M}_1\otimes_k{\rm Sym}_\bullet^k({\rm Lie} ({\cal
H}))\otimes_k\omega_{\cal H}^{\otimes -1},\cr &\kern 5cm {\cal
D}_{M_1^0}\otimes_k{\rm Sym}_\bullet^k({\rm Lie} ({\cal
H}))\otimes_k\omega_{{\cal G}_1^0}^{\otimes -1}\otimes
\omega_{\cal H}^{\otimes -1}\bigr)\cr
&\cong R\underline {\rm Hom}_{{\cal D}_{M_1^0}}\bigl({\cal
M}_1,{\cal D}_{M_1^0}\otimes_k\omega_{{\cal G}_1^0}^{\otimes
-1}\bigr)\otimes_k{\rm Sym}_\bullet^k({\rm Lie}({\cal H}))\cr
&\cong f_*R\underline {\rm Hom}_{D_{\rm qcoh}(M_1)}\bigl({\cal
M}_1,{\cal D}_{M_1}^\Omega\bigr)\otimes_k
\omega_{\cal H}.\cr}
\qquad}
$$
Dans le second, on a
$$
\eqalign{
R\underline {\rm Hom}_{D_{\rm qcoh}(M_2)}(f_*{\cal M}_1,{\cal
D}_{M_2}^\Omega) &\cong R\underline {\rm Hom}_{D_{\rm
qcoh}(M_2)}\bigl(\omega_{\cal H}\,\Lotimes_{k[{\cal H}]}\,{\cal
M}_1, {\cal D}_{M_2}^\Omega\bigr)\cr
&\cong R\underline {\rm Hom}_{D_{\rm qcoh}(M_1)}\bigl({\cal
M}_1,{\cal D}_{M_2}^\Omega\otimes_k k[{\cal H}]\bigr)\cr
&\cong f_*R\underline {\rm Hom}_{D_{\rm qcoh}(M_1)}
\bigl({\cal M}_1,{\cal D}_{M_1}^\Omega\bigr)\otimes_k
\omega_{\cal H}\cr}
$$
(si $S={\rm Spec}(k[{\cal H}])$, si $\pi :S \rightarrow {\rm
Spec}(k)$ est la projection canonique et si $\epsilon :{\rm Spec}(k)
\hookrightarrow S$ est l'origine du tore
$S$, on a
$$\eqalign{
R{\rm Hom}_k(\epsilon^*{\cal F},{\cal G})&\cong
\epsilon^!R\underline {\rm Hom}_{{\cal O}_S} ({\cal F},\pi^!{\cal
G})\cr
&\cong \epsilon^*R\underline {\rm Hom}_{{\cal O}_T}({\cal
F},\pi^*{\cal G})\cr}
$$
pour tout ${\cal O}_S$-Module coh\'erent ${\cal F}$ et tout
$k$-espace vectoriel ${\cal G}$ puisque ${\cal G}\cong
\epsilon^!\pi^!{\cal G}$, que $\pi^!(\cdot )\cong \pi^*(\cdot
)\otimes_k\omega_S[d_S]$ et que $\epsilon^!(\cdot
)\cong\epsilon^*(\cdot )\otimes_k\omega_S^{\otimes -1}[-d_S]$).

Passons maintenant \`a la partie (ii), toujours en supposant
$f_G={\rm id}$. On ne peut plus utiliser la factorisation de
$f$ de la section (6.1) puisque les monomorphismes $g_{\cal G}$ et
$g_{\cal G}'$ n'ont pas un conoyau fini. Cependant, il suffit de
v\'erifier les assertions de la partie (ii) dans les cas particuliers
(3) et (5) de (6.1) et dans le cas particulier o\`u $f_{\cal G}^0={\rm
id}$ et $f_{\cal G}^{\rm \acute et}$ est une isog\'enie, i.e. un
monomorphisme \`a conoyau fini.  Faisons le par exemple dans le cas
particulier (3), en supposant plus ${\cal G}_1^{\rm \acute et} = {\cal
G}_2^{\rm\acute et}=(0)$, et dans le cas particulier o\`u $f_{\cal
G}^0={\rm id}$ et $f_{\cal G}^{\rm \acute et}$ est une isog\'enie,
en supposant de plus que ${\cal G}_1^0={\cal G}_2^0=(0)$ et que
${\cal G}_1^{\rm \acute et}$ et ${\cal G}_2^{\rm \acute et}$ sont
d\'eploy\'es sur $k$.

Dans le premier de ces cas, on a
$$
\displaylines{
\qquad R\underline {\rm Hom}_{D_{\rm qcoh}(M_1)}(f^!{\cal
M}_2,{\cal D}_{M_1}^\Omega)
\hfill\cr\hfill
\eqalign{&\cong R\underline {\rm Hom}_{{\cal D}_{M_2^0}
\otimes_k{\rm Sym}_\bullet^k({\rm Lie}({\cal H}))}\bigl({\cal
M}_2\otimes_k\omega_{\cal H}^{\otimes -1},\cr &\kern 5cm {\cal
D}_{M_1^0}^\Omega\otimes_k {\rm Sym}_\bullet^k({\rm Lie}({\cal
H}))\otimes_k\omega_{\cal H}^{\otimes -1}\bigr)\cr
&\cong R\underline {\rm Hom}_{D_{\rm qcoh}(M_2)}
\bigl({\cal M}_2,{\cal D}_{M_2}^\Omega\bigr)
\otimes_kR{\rm Hom}_{{\rm Sym}_\bullet^k({\rm Lie}({\cal
H}))}\bigl(k,{\rm Sym}_\bullet^k({\rm Lie}({\cal H}))\bigr)\cr
&\cong R\underline {\rm Hom}_{D_{\rm qcoh}(M_2)}\bigl({\cal
M}_2,{\cal D}_{M_2}^\Omega\bigr)\otimes_k\omega_{\cal
H}[-d_{\cal H}]\cr}\qquad}
$$
(si $V={\rm Spec}({\rm Sym}_\bullet^k({\rm Lie}({\cal H}))$ et si $
\epsilon :{\rm Spec}(k) \hookrightarrow V$ est l'origine, on a
$$
R\underline {\rm Hom}_{{\cal O}_V}\bigl(\epsilon_*k,{\cal
O}_V)\cong\epsilon_*R{\rm Hom}_k\bigl(k,\epsilon^!{\cal
O}_V)\cong\epsilon_*\omega_V^{\otimes -1}[-d_V]
$$
avec $\omega_V=\omega_{\cal H}^{\otimes -1}$ et $d_V=d_{\cal
H}$).

Dans le second, posons $X_i={\cal G}_i^{\rm \acute et}(k)$ pour
$i=1,2$ et identifions $X_1$ \`a un sous-module de $X_2$ ($X_1$ et
$X_2$ sont tous les deux des ${\bb Z}$-modules libres de m\^eme
rang fini $r$ et $X_2/X_1$ est un groupe ab\'elien fini). Alors, on a
$$
f^!{\cal M}_2=\rho ({\cal M}_2)\otimes_k \omega_2 \otimes_k
\omega_1^{\otimes -1},
$$
o\`u $\rho ({\cal M}_2)$ est le ${\cal O}_G$-Module
$X_1$-\'equivariant sous-jacent au ${\cal O}_G$-Module
$X_2$-\'equivariant ${\cal M}_2$ et o\`u on a pos\'e $\omega_i=
\bigwedge^r{\rm Hom}_{\bb Z}(X_i,k)$ pour $i=1,2$. Par suite, si
l'on note simplement ${\cal O}_G[X_i]$ le ${\cal O}_G$-Module
$X_i$-\'equivariant $\bigoplus_{x_i\in X_i}\tau_{u_i^{\rm \acute
et}(x_i)}^*{\cal O}_G$ pour la $X_i$-structure induite par la
translation de $X_i$ sur lui-m\^eme, on a
$$
R\underline {\rm Hom}_{D_{\rm qcoh}(M_1)}(f^!{\cal M}_2, {\cal
D}_{M_1}^\Omega )\cong R\underline {\rm Hom}_{{\cal O}_G,X_1}
\bigl(\rho ({\cal M}_2),{\cal O}_G[X_1]\bigr)\otimes_k
\omega_2^{\otimes -1}
$$
et
$$
f^!R\underline {\rm Hom}_{D_{\rm qcoh}(M_2)}({\cal M}_2,
{\cal D}_{M_2}^\Omega )\cong\rho\bigl(R\underline {\rm
Hom}_{{\cal O}_G,X_2}({\cal M}_2, {\cal O}_G[X_2])\bigr) \otimes_k
\omega_1^{\otimes -1},
$$
o\`u
$$
R\underline {\rm Hom}_{{\cal O}_G,X_i}(\,\cdot\,,{\cal O}_G[X_i])
$$
est muni de la $X_i$-structure induite par celle de ${\cal O}_G[X_i]$
qui provient de la $X_i$-structure naturelle de ${\cal O}_G$. Il ne
reste plus qu'\`a montrer que
$$
R\underline {\rm Hom}_{{\cal O}_G,X_1}\bigl(\rho ({\cal M}_2),{\cal
O}_G[X_1]\bigr)
$$
est isomorphe \`a
$$
\rho\bigl(R\underline {\rm Hom}_{{\cal O}_G,X_2} ({\cal M}_2,{\cal
O}_G[X_2])\bigr)\otimes_k\omega_1\otimes_k\omega_2^{\otimes
-1}.
$$
Mais, la $X_2$-structure de ${\cal M}_2$ et celle par translation de
${\cal O}_G[X_2]$ induisent une action du groupe fini $X_2/X_1$ sur
$R\underline {\rm Hom}_{{\cal O}_G,X_1}\bigl(\rho ({\cal M}_2),
{\cal O}_G[X_2]\bigr)$ et on a un quasi-isomorphisme naturel de
complexes de ${\cal O}_G$-Modules
$$
R\underline {\rm Hom}_{{\cal O}_G,X_2}({\cal M}_2,{\cal O}_G[X_2])
\cong\bigl(R\underline {\rm Hom}_{{\cal O}_G,X_1}(\rho ({\cal
M}_2), {\cal O}_G[X_2])\bigr)^{X_2/X_1}
$$
et notre assertion en r\'esulte puisque ${\cal O}_G[X_2]\cong
\bigoplus_{y\in X_2/X_1}\tau_{u_2^{\rm \acute et}(y)}^*{\cal
O}_G[X_1]$.

\hfill\hfill\cqfd
\vskip 5mm

(6.3) On est maintenant en mesure de d\'efinir une transformation
de Fourier pour les $1$-motifs g\'en\'eralis\'es.

Soit $M$ un objet de ${\bf Mot}/k$, de dual $M'$. Le ${\cal
O}_{G'\times_kG}$-Module inversible ${\cal Q}$ de Poincar\'e avec
les structures suppl\'ementaires introduites en (5.2) est un Module
inversible sur $M'\times M$. On d\'efinit alors un  foncteur
$$
{\cal F}_M:D_{\rm qcoh}^{\rm b}(M)\rightarrow D_{\rm qcoh}^{\rm
b}(M')\leqno (6.3.1)
$$
par
$$
{\cal F}_M(\cdot )={\rm pr}_*'({\cal Q}\otimes^!{\rm pr}^!(\cdot ))
$$
o\`u ${\rm pr}$, ${\rm pr}'$ sont les projections canoniques de
$M'\times M$ sur $M$ et $M'$ respectivement.
\vskip 2mm

{\pc CAS PARTICULIERS} (6.3.2)\pointir {\rm (i)} Si $M=
[0\rightarrow A]$ avec $A$ un $k$-sch\'ema ab\'elien de dimension
$g$, on a
$$
{\cal F}_M(\cdot )={\cal F}(\cdot )\otimes_k\omega_A^{\otimes
-1}[-g].
$$

\decale{\rm (ii)} Si $M=[\widehat A^0\rightarrow A]$ avec $A$ un
$k$-sch\'ema ab\'elien de dimension $g$, on a
$$
{\cal F}_M(\cdot )=\widetilde{\cal F}(\cdot )[-g]
$$
et si $M=[0\rightarrow A^\natural]$ avec $A$ un $k$-sch\'ema
ab\'elien de dimension $g$, on a
$$
{\cal F}_M(\cdot )=\widetilde{\cal F}^\natural (\cdot )[-2g].
$$

\decale{\rm (iii)} Si $M=[0\rightarrow V]$ avec $V$ un
$k$-vectoriel de dimension $d_V$, on a $M'=[\widehat
V'^0\rightarrow 0]$, o\`u $V'$ est le $k$-vectoriel dual de $V$, et
$$
{\cal F}_M(\cdot )=R\Gamma (V,\,\cdot\,)\otimes_k
\omega_V^{\otimes -1}[-d_V]
$$
en tant que complexe de modules sur $H^0(V,{\cal O}_V)={\rm
Sym}_\bullet^k(V')$. Si $M=[\widehat V^0\rightarrow 0]$ avec $V$ un
$k$-vectoriel de dimension $d_V$, de dual not\'e $V'$, on a
$M'=[0\rightarrow V']$ et
$$
{\cal F}_M(\cdot )=(\cdot )\widetilde{\,\,}\,\otimes_k
\omega_V^{\otimes 2},
$$
o\`u $(\cdot )\widetilde{\,\,}$ est le ${\cal O}_{V'}$-Module
quasi-coh\'erent d'espace des sections globales le ${\rm
Sym}_\bullet^k(V)$-module $(\cdot )$.

\decale{\rm (iv)} Si $M=[\widehat V^0\buildrel {\rm can}\over
\longrightarrow V]$ avec $V$ un $k$-vectoriel de dimension $d_V$,
de dual not\'e $V'$, on a $M'=[\widehat V'^0\buildrel {\rm
can}'\over\longrightarrow V']$ et, pour tout
${\cal D}_V$-Module quasi-coh\'erent ${\cal M}$, on a
$$
R\Gamma (V',{\cal F}_M({\cal M}))=H^0(V,{\cal M})[-d_V]
$$
en tant que complexe de $k$-vectoriels et l'action de $H^0(V',{\cal
D}_{V'})$ est induite par celle de $H^0(V,{\cal D}_V)$, compte tenu
de l'identification habituelle de
$$
H^0(V',{\cal D}_{V'})={\rm Sym}_\bullet^k(V)\otimes_k{\rm
Sym}_\bullet^k(V')
$$
avec
$$
H^0(V,{\cal D}_V)={\rm Sym}_\bullet^k(V')\otimes_k{\rm
Sym}_\bullet^k(V).
$$

\decale{\rm (v)} Si $M=[0\rightarrow T]$ avec $T$ un $k$-tore
d\'eploy\'e de dimension $d_T$, de groupe des caract\`eres not\'e
$X$, on a $M'=[X\rightarrow 0]$ et
$$
{\cal F}_M(\cdot )=R\Gamma (T,\,\cdot\,)\otimes_k
\omega_T^{\otimes -1}[-d_T]
$$
en tant que complexe de modules sur $H^0(T,{\cal O}_V)=k[X]$. Si
$M=[X\rightarrow 0]$ avec $X$ un ${\bb Z}$-module libre de rang
fini, on a $M'=[0\rightarrow T]$, o\`u $T$ est le $k$-tore dont le
groupe des caract\`eres est $X$, et
$$
{\cal F}_M(\cdot )=(\cdot )\widetilde{\,\,}\,\otimes_k
\omega_X^{\otimes 2},
$$
o\`u $(\cdot )\widetilde{\,\,}$ est le ${\cal O}_T$-Module
quasi-coh\'erent d'espace des sections globales le $k[X]$-module
$(\cdot )$.

\decale{\rm (vi)} Si $M=[\widehat T^0\buildrel {\rm can}
\over\longrightarrow T]$ avec $T$ un $k$-tore d\'eploy\'e de
dimension $d_T$, on a $M'=[X\buildrel {\rm can}'\over
\longrightarrow\Omega ]$ avec les notations de (5.2.5)(ii) et, pour
tout ${\cal D}_T$-Module quasi-coh\'erent ${\cal M}$, on a
$$
R\Gamma (V,{\cal F}_M({\cal M}))=H^0(T,{\cal M})[-d_T]
$$
en tant que complexe de $k$-vectoriels et la structure de Module
sur $M'$ est induite par la structure de $H^0(T,{\cal D}_T)$-module,
compte tenu de l'identification
$$
H^0(T,{\cal D}_T)=k[X]\otimes_kH^0(\Omega ,{\cal O}_\Omega ).
$$

\hfill\hfill\cqfd
\vskip 2mm

Les propri\'et\'es les plus importantes de la transformation de
Fourier pour les $1$-motifs g\'en\'eralis\'es sont les suivantes.

\th TH\'EOR\`EME (6.3.3)
\enonce
{\rm (i)} Les foncteurs exacts ${\cal F}_M$ et ${\cal F}_{M'}$ sont
des \'equivalences de cat\'egories triangul\'ees. Plus
pr\'ecis\'ement,  on a les isomorphismes de foncteurs
$$
{\cal F}_{M'}\circ {\cal F}_M\cong \langle -1\rangle^!(\cdot )
\otimes_k\omega_M^{\otimes -2}\otimes_k
\omega_{M'}^{\otimes -1}[-2d_G+d_{{\cal G}'}+r_{{\cal G}'}-d_{G'}].
$$

\decale{\rm (ii)} Pour tout morphisme $f:M_1\rightarrow M_2$
dans ${\bf Mot}/k$ de transpos\'e $f':M_2' \rightarrow M_1'$, on a
l'isomorphisme de foncteurs
$$
{\cal F}_{M_2}\circ f_*(\cdot)\otimes_k\omega_{M_2'}\otimes_k
\omega_{M_2}[d_{G_2'}+d_{G_2}]\cong f'^!\circ {\cal F}_{M_1}
(\cdot)\otimes_k\omega_{M_1'}\otimes_k\omega_{M_1}[d_{G_1'}+
d_{G_1}]
$$
de $D_{\rm qcoh}^{\rm b}(M_1)$ dans $D_{\rm qcoh}^{\rm b}(M_2')$.

\decale{\rm (iii)} Le foncteur ${\cal F}_M$ envoie $D_{\rm
coh}^{\rm b}(M)$ dans $D_{\rm coh}^{\rm b}(M')$.

\decale{\rm (iv)} On a l'isomorphisme de foncteurs
$$
D'\circ {\cal F}_M(\cdot )\cong \langle -1\rangle^!\circ {\cal
F}_M\circ D\otimes_k\omega_M^{\otimes 3}[d_G-d_{\cal G}-r_{\cal
G}+ 2d_{G'}]
$$
de $D_{\rm coh}^{\rm b}(M)$ dans $D_{\rm coh}^{\rm b}(M')$, o\`u
$D$ et $D'$ sont les foncteurs de dualit\'e pour $M$ et $M'$
respectivement.

\endth

On remarquera que $\omega_{M'}=\omega_M$ puisque $\omega_A'=
\omega_A$, $\omega_{V'}=\omega_{{\cal G}^0}^{\otimes -1}$,
$\omega_{T'}=\omega_{{\cal G}^{\rm \acute et}}^{\otimes -1}$,
$\omega_{{\cal G}'^0}=\omega_V$ et $\omega_{{\cal G}'^{\rm \acute
et}}=\omega_T$ et que $-2d_{G'}+d_{\cal G}+r_{\cal G}-d_G=
-2d_G+d_{{\cal G}'}+r_{{\cal G}'}-d_{G'}$ puisque $d_{A'}=d_A$,
$d_{V'}=d_{\cal G}$, $d_{T'}=r_{\cal G}$, $d_{{\cal G}'}=d_V$ et
$r_{{\cal G}'}=d_T$.

\rem Preuve
\endrem
L'assertion (ii) r\'esulte facilement du th\'eor\`eme de
changement de base, de la formule des projections et de
l'isomorphisme canonique
$$
({\rm id}\times f)^!{\cal Q}_2\otimes_k
\omega_{M_2'}\otimes_k\omega_{M_2}[d_{G_2'}+d_{G_2}]
\cong (f'\times {\rm id})^!{\cal Q}_1\otimes_k
\omega_{M_1'}\otimes_k\omega_{M_1}[d_{G_1'}+d_{G_1}]
$$
sur $M_2'\times M_1$.

L'assertion (iii) est une cons\'equence de l'assertion (ii) appliqu\'ee
au morphisme $(0,{\rm id}): [0\rightarrow G]\rightarrow [{\cal
G}\rightarrow G]$.

On se contentera de v\'erifier les assertions (i) et (iv) dans les cas
particuliers (6.3.2).

Si $M=[0\rightarrow A]$, ces assertions r\'esultent
aussit\^ot du th\'eor\`eme (1.2.4) et de la proposition (1.3.4).

Si $M=[\widehat A^0\rightarrow A]$ ou $M=[0\rightarrow A^\natural
]$, ces assertions r\'esultent aussit\^ot du th\'eor\`eme (3.2.1) et
de la proposition (3.3.4).

Si $M=[0\rightarrow V]$ et $M'=[\widehat V^0\rightarrow 0]$, on a
$$
D(\cdot )=R\underline{\rm Hom}_{{\cal O}_V}(\,\cdot\,,{\cal O}_V)
\otimes_k\omega_V[d_V]
$$
et
$$
D'(\cdot )=R{\rm Hom}_{{\rm Sym}_\bullet^k(V')}(\,\cdot\,,{\rm
Sym}_\bullet^k(V'))\otimes_k \omega_{V'}^{\otimes -2}
$$
et les assertions (i) et (iv) du th\'eor\`eme sont triviales.

Si $M=[\widehat V^0\rightarrow V]$, on a ${\cal F}_M={\cal F}(\cdot
) [-d_V]$ o\`u ${\cal F}$ est la transformation de Fourier usuelle
pour les ${\cal D}_V$-Modules quasi-coh\'erents, transformation de
Fourier qui v\'erifie trivialement ${\cal F}'\circ {\cal F} \cong
\langle -1\rangle^!$ et $D'\circ {\cal F}\cong \langle -1 \rangle^!
\circ {\cal F}\circ D$, d'o\`u les assertions (i) et (iv) du
th\'eor\`eme dans ce cas.

Les cas particuliers $M=[0\rightarrow T]$ et $M=[\widehat T^0
\rightarrow T]$ se traitent de m\^eme.

\hfill\hfill\cqfd

\vfill\eject
\centerline{\bf Bibliographie}
\vskip 20mm
\parindent=10mm

\newtoks\ref
\newtoks\auteur
\newtoks\titre
\newtoks\editeur
\newtoks\annee
\newtoks\revue
\newtoks\tome
\newtoks\pages
\newtoks\reste
\newtoks\autre

\def\livre{\leavevmode
\llap{[\the\ref]\enspace}%
\the\auteur\pointir
{\sl \the\titre},
\the\editeur,
{\the\annee}.
\smallskip
\filbreak}

\def\article{\leavevmode
\llap{[\the\ref]\enspace}%
\the\auteur\pointir
\the\titre,
{\sl\the\revue}
{\bf\the\tome},
({\the\annee}),
\the\pages.
\smallskip
\filbreak}

\def\autre{\leavevmode
\llap{[\the\ref]\enspace}%
\the\auteur\pointir
\the\reste.
\smallskip
\filbreak}


\ref={Bo}
\auteur={A. {\pc BOREL} et al.}
\titre={Algebraic $D$-Modules}
\editeur={Academic Press}
\annee={1987}
\livre

\ref={De}
\auteur={P. {\pc DELIGNE}}
\titre={Th\'eorie de Hodge, III}
\revue={Publ. Math. IHES}
\tome={44}
\annee={1974}
\pages={5-78}
\article

\ref={Fa-Ch}
\auteur={G. {\pc FALTINGS} and C.-L. {\pc CHAI}}
\titre={Degeneration of Abelian Varieties}
\editeur={Springer-Verlag}
\annee={1990}
\livre

\ref={Fo}
\auteur={J.-M. {\pc FONTAINE}}
\titre={groupes $p$-divisibles sur les corps locaux}
\editeur={Ast\'erisque 47-48}
\annee={1977}
\livre

\ref={Ha}
\auteur={R. {\pc HARTSHORNE}}
\titre={Residues ans Duality}
\editeur={Lectures Notes in Mathematics 20, Springer-Verlag}
\annee={1966}
\livre

\ref={Gr}
\auteur={A. {\pc GROTHENDIECK}}
\titre={Technique de descente et th\'eor\`emes d'existence
en g\'eom\'e-\break trie alg\'ebrique, V. Les sch\'emas de Picard.
Th\'eor\`emes d'existence}
\editeur={S\'eminaire Bourbaki 232}
\annee={1961/62}
\livre

\ref={Ma-Me}
\auteur={B. {\pc MAZUR} and W. {\pc MESSING}}
\titre={Universal Extensions and One Dimensional Crystalline
Cohomology}
\editeur={Lectures Notes in Mathematics 370, Springer-Verlag}
\annee={1974}
\livre

\ref={Muk}
\auteur={S. {\pc MUKAI}}
\titre={Duality between $D(X)$ and $D(\widehat X)$ with its
applications to Picard sheaves}
\revue={Nagoya Math. J.}
\tome={81}
\annee={1981}
\pages={153-175}
\article

\ref={Mum}
\auteur={D. {\pc MUMFORD}}
\titre={Abelian Varieties}
\editeur={Tata Institute of Fundamental Research, Bombay,
et Oxford University Press}
\annee={1985}
\livre

\ref={Ro}
\auteur={M. {\pc ROTHSTEIN}}
\reste={{\it Sheaves with connection on abelian varieties},
Preprint, alg-geom/9602023}
\autre

\ref={Se}
\auteur={J.-P. {\pc SERRE}}
\titre={Groupes alg\'ebriques et corps de classes}
\editeur={Hermann}
\annee={1959}
\livre

\ref={SGA3}
\auteur={M. {\pc DEMAZURE} et A. {\pc GROTHENDIECK}}
\titre={Propri\'et\'es G\'en\'erales des Sch\'emas en Groupes}
\editeur={Lectures Notes in Mathematics 151, Springer-Verlag}
\annee={1959}
\livre

\ref={SGA6}
\auteur={P. {\pc BERTHELOT}, A. {\pc GROTHENDIECK} et L. {\pc
ILLUSIE}}
\titre={Th\'eorie des Intersections et Th\'eor\`eme de
Riemann-Roch}
\editeur={Lectures Notes in Mathematics 225, Springer-Verlag}
\annee={1971}
\livre

\vskip 10mm
\let\+=\tabalign

\line{\hbox{\kern 5mm{\vtop{\+ G\'erard LAUMON\cr
\+ URA 752 du CNRS\cr
\+ Universit\'e de Paris-Sud\cr
\+ Math\'ematiques, b\^at. 425\cr
\+ 91405 ORSAY C\'edex (France)\cr}}}}

\bye